\documentclass[12pt,reqno]{amsart}
\usepackage{amsmath,amssymb,geometry,color,tikz-cd,extarrows}
\usepackage[backref,pagebackref,pdftex,hyperindex]{hyperref}
\usepackage[alphabetic]{amsrefs}
\usepackage{amssymb}
\usepackage{bbm}
\usepackage{enumitem}
\usepackage{upgreek}

\newtheorem{proposition}{Proposition}[section]
\newtheorem{theorem}[proposition]{Theorem}
\newtheorem{lemma}[proposition]{Lemma}
\newtheorem{corollary}[proposition]{Corollary}
\newtheorem{conjecture}[proposition]{Conjecture}

\theoremstyle{definition}
\newtheorem{remark}[proposition]{Remark}
\newtheorem{definition}[proposition]{Definition}
\newtheorem{example}[proposition]{Example}

\title{Non-proportional wall crossing for K-stability}

\author{Yuchen Liu}
\address{Department of Mathematics, Northwestern University, Evanston, IL 60208, USA.}
\email{yuchenl@northwestern.edu}

\author{Chuyu Zhou}
\address{School of Mathematical Sciences, Xiamen University, Siming South Road 422, Xiamen, Fujian 361005, China.}
\email{chuyuzhou@xmu.edu.cn, chuyuzhou1@gmail.com}

\date{} 

\thanks{
        \newline
		\indent
		\emph{Competing interests}:  The authors have no conflict of interest to declare.
		}

\newcommand{\Fut}{{\rm{Fut}}}

\newcommand{\ord}{{\rm {ord}}}
\newcommand{\tc}{{\rm {tc}}}
\newcommand{\vol}{{\rm {vol}}}

\newcommand{\red}{{\rm {red}}}
\newcommand{\Spec}{{\rm {Spec}}}

\newcommand{\Exc}{{\rm {Exc}}}
\newcommand{\Val}{{\rm {Val}}}
\newcommand{\wt}{{\rm {wt}}}

\newcommand{\Pic}{{\rm {Pic}}}

\newcommand{\dt}{{\rm {dt}}}

\newcommand{\Supp}{{\rm {Supp}}}

\newcommand{\ac}{{\rm {ac}}}

\newcommand{\SL}{{\rm {SL}}}
\newcommand{\Hilb}{{\rm {Hilb}}}

\newcommand{\PGL}{{\rm {PGL}}}

\newcommand{\Aut}{{\rm {Aut}}}

\newcommand{\CM}{{\rm {CM}}}

\newcommand{\Kss}{{\rm {Kss}}}

\newcommand{\WLF}{{\rm {WLF}}}

\newcommand{\LF}{{\rm {LF}}}
\newcommand{\Nef}{{\rm {Nef}}}
\newcommand{\NS}{{\rm {NS}}}

\newcommand{\OD}{{\rm {OD}}}
\newcommand{\Kps}{{\rm {Kps}}}
\newcommand{\Hom}{{\rm {Hom}}}
\newcommand{\GIT}{{\rm {GIT}}}
\newcommand{\Isom}{{\rm {Isom}}}
\newcommand{\RC}{{\rm {RC}}}

\newcommand{\bA}{\mathbb{A}}

\newcommand{\bC}{\mathbb{C}}

\newcommand{\bN}{\mathbb{N}}

\newcommand{\bP}{\mathbb{P}}
\newcommand{\bQ}{\mathbb{Q}}
\newcommand{\bR}{\mathbb{R}}

\newcommand{\bT}{\mathbb{T}}

\newcommand{\bZ}{\mathbb{Z}}

\newcommand{\mA}{\mathcal{A}}
\newcommand{\mB}{\mathcal{B}}

\newcommand{\mD}{\mathcal{D}}
\newcommand{\mE}{\mathcal{E}}
\newcommand{\mF}{\mathcal{F}}
\newcommand{\mG}{\mathcal{G}}
\newcommand{\mH}{\mathcal{H}}

\newcommand{\mL}{\mathcal{L}}
\newcommand{\mM}{\mathcal{M}}

\newcommand{\mO}{\mathcal{O}}
\newcommand{\mP}{\mathcal{P}}

\newcommand{\mX}{\mathcal{X}}
\newcommand{\mY}{\mathcal{Y}}

\newcommand{\fD}{\mathbf{D}}
\newcommand{\fE}{\mathbf{E}}

\newcommand{\fX}{\mathbf{X}}
\newcommand{\fY}{\mathbf{Y}}

\newcommand{\fa}{\mathbf{a}}

\newcommand{\vc}{\vec{c}}
\newcommand{\vx}{\vec{x}}
\newcommand{\va}{\vec{a}}

\newcommand{\kA}{\mathfrak{A}}

\newcommand{\YL}[1]{{\textcolor{blue}{[Yuchen: #1]}}}

\begin{document}

\begin{abstract}
In this paper, we present a general wall crossing theory for K-stability and K-moduli of log Fano pairs whose boundary divisors can be non-proportional to the anti-canonical divisor. Along the way, we prove 
that there are only finitely many K-semistable domains associated to the fibers of a log bounded family of couples. Under the additional assumption of volume bounded from below, we show that K-semistable domains are semi-algebraic sets (although not necessarily polytopes).  As a consequence, we obtain a finite semi-algebraic chamber decomposition for wall crossing of K-moduli spaces. 
In the case of one boundary divisor, this decomposition is an expected  finite interval chamber decomposition. As an application of the theory, we prove a comparison theorem between GIT-stability and K-stability in non-proportional setting when the coefficient of the boundary is sufficiently small.

\end{abstract}

\maketitle

\setcounter{tocdepth}{1}

\tableofcontents

We work over the field of complex numbers  $\bC$ throughout the article.

\section{Introduction}

The wall crossing phenomenon is a novel characteristic for moduli spaces of higher dimensional algebraic varieties, and it plays an important role in the study of birational geometry of various moduli spaces of algebraic varieties. In the canonically polarized case, the wall crossing framework for KSBA moduli spaces was established \cite{ABIP21, MZ23}. In the log Fano case, from the theoretical viewpoint, there has been a rather satisfactory theory for wall crossing of K-moduli spaces in the proportional case, i.e. where the boundary divisors are proportional to the anti-canonical divisor (e.g. \cite{ADL19, Zhou23, Zhou23b}). Moreover, many explicit wall crossing examples are studied in recent years. For instance, an explicit birational geometry picture of moduli spaces of K3 surfaces of degree four is worked out in \cite{ADL23} by wall crossing of K-stability. However, it is still a little mysterious in the non-proportional case, i.e. where the boundary divisors are not necessarily proportional to the anti-canonical divisor. In this paper, we attempt to find the suitable setting to formulate this problem and present the non-proportional wall crossing phenomenon.

\subsection{Wall crossing of K-moduli}\label{subsec: wall-crossing}
To study the wall crossing phenomenon, similar as the approach in \cite{Zhou23b}, we a priori construct a set of couples which is large enough to contain all the objects we are interested. 
Fixing two positive integers $d, k$ and a rational polytope $P\subset [0,1]^k$, we say a couple $(X, \sum_{j=1}^kD_j)$ belongs to $\mG:=\mG(d,k,P)$ if and only if the following conditions are satisfied:
\begin{enumerate}
\item $X$ is a normal projective variety of dimension $d$;
\item $D_j, 1\leq j\leq k,$ are effective Weil divisors on $X$;
\item there exists a vector $(c_1,...,c_k)\in P$ such that $(X, \sum_{j=1}^kc_jD_j)$ is a K-semistable log Fano pair.\footnote{In this paper, log Fano pairs are allowed to have $\bR$-coefficients, and K-semi(poly)stability is defined for them in Section \ref{sec: preliminaries}.}
\end{enumerate}
Note that in the definition of the set $\mG$, we do not assume that $X$ is Fano and $D_j$'s are proportional to $-K_X$. This is the essential difference compared to the setup in \cite{Zhou23b}.
On the other hand, as $\mG$ is not necessarily log bounded, we just fix a log bounded subset $\widetilde{\mG}$ satisfying the following volume condition $(\clubsuit)$: \textit{there exists a positive real number $\epsilon>0$ such that $\vol(-K_X-\sum_{j=1}^k x_jD_j)\geq \epsilon$ for any $(X, \sum_{j=1}^k D_j)\in \widetilde{\mG}$ and any $(x_1,...,x_k)\in P$.} For such a subset $\widetilde{\mG}$, we have the following result on the chamber decomposition of $P$ to control the variation of K-semi(poly)stability.

\begin{theorem}\label{mainthm: kps fcd}{\rm (Theorem \ref{thm: fcd}, Theorem \ref{thm: kps fcd})}
Suppose $\widetilde{\mG}\subset \mG$ is a subset which is log bounded and satisfies the volume condition $(\clubsuit)$. Then there exists a  decomposition of $P$ into finite disjoint subsets (depending only on $\widetilde{\mG}$), denoted by
$P=\amalg_{i=1}^{m} A_{i}, $
such that: (1)  each $A_{i}$ is a semi-algebraic set homeomorphic to $(0,1)^r$ for some $r\in \bN$ (with $(0,1)^0$ being a point);
(2) the K-semi(poly)stability of $(X, \sum_{j=1}^kx_jD_j)$ does not change as $(x_1,...,x_k)$ varies in $A_{i}$ for any $0\leq i\leq m$ and any $(X, \sum_{j=1}^kD_j)\in \widetilde{\mG}$.
\end{theorem}

We refer to Definition \ref{def: semi-algebraic} and Proposition \ref{prop: semi-algebraic property} for the definition and properties of semi-algebraic sets. Theorem \ref{mainthm: kps fcd} enables us to characterize the wall crossing phenomenon of K-moduli when the coefficients vary. When we consider the K-moduli parametrizing certain couples in $\widetilde{\mG}$, probably we will not get a proper one since $\widetilde{\mG}$ may not be complete in moduli sense. So we consider the moduli completion $\widehat{\mG}$ \footnote{We should take the notation $\widehat{\widetilde{\mG}}$ rather than $\widehat{\mG}$, since we are considering the completion of $\widetilde{\mG}$. As there should not be a confusion, we just take the notation $\widehat{\mG}$ for simplicity.} of $\widetilde{\mG}$ as in the following definition.

\begin{definition}
We say a couple $(Y, \sum_{j=1}^kB_j)\in \mG$ is contained in $\widehat{\mG}$ if and only if there exists a family of couples $(\mY, \sum_{j=1}^k\mB_j)\to C\ni 0$ over a smooth pointed curve such that the following conditions are satisfied:
\begin{enumerate}
\item $(\mY_t, \sum_{j=1}^k \mB_{j,t})\in \widetilde{\mG}$ for $t\ne 0$;
\item $(\mY_0, \sum_{j=1}^k \mB_{j,0})\cong (Y, \sum_{j=1}^kB_j)$;
\item for some $(c_1,...,c_k)\in P$, $(\mY, \sum_{j=1}^k c_j\mB_j)\to C$ is a family of K-semistable log Fano pairs.
\end{enumerate}
\end{definition}
In other words, $\widehat{\mG}$ is obtained by adding to $\widetilde{\mG}$ all the K-semistable degenerations of K-semistable log Fano pairs arising from $\widetilde{\mG}$. It is not hard to check that $\widehat{\mG}$ is still log bounded satisfying the volume condition $(\clubsuit)$ (e.g. Proposition \ref{prop: G-complete}). \textit{Hence we can freely replace $\widetilde{\mG}$ in Theorem \ref{mainthm: kps fcd} with its moduli completion $\widehat{\mG}$. }
For each vector $\fa:=(a_1,...,a_k)\in P$, we consider the following subset of $\widehat{\mG}$:
$$\widehat{\mG}_\fa:=\{(Y, \sum_{j=1}^kB_j)\in \widehat{\mG}\ |\ \textit{$(Y, \sum_{j=1}^ka_jB_j)$ is a K-semistable log Fano pair}\}. $$
By \cite[Theorem 1.3]{LXZ22} complemented by \cite{LZ24}, there exists a finite type Artin stack, denoted by $\mM^{\Kss}_{\widehat{\mG}, \fa}$, parametrizing the following set of K-semistable log Fano pairs
$$\{(Y, \sum_{j=1}^ka_jB_j)\ |\ \textit{$(Y, \sum_{j=1}^kB_j)\in \widehat{\mG}_\fa$}\}, $$
and $\mM^{\Kss}_{\widehat{\mG}, \fa}$ admits a proper good moduli space $M_{\widehat{\mG},\fa}^{\Kps}$. Let $\mM^{\Kss, \red}_{\widehat{\mG},\vec{x}}$ be the moduli functor associated to reduced schemes and  $M^{\Kps,\red}_{\widehat{\mG}, \vec{x}}$ the corresponding good moduli space.  Denote $\phi_\fa: \mM^{\Kss,\red}_{\widehat{\mG}, \fa}\to
M^{\Kps,\red}_{\widehat{\mG}, \fa}$. We have the following result.

\begin{theorem}\label{mainthm: moduli-invariance}{\rm (Theorem \ref{thm: moduli invariance})}
Let $A_i$ be a chamber in the chamber decomposition in Theorem \ref{mainthm: kps fcd}. Then $\phi_\fa: \mM^{\Kss,\red}_{\widehat{\mG}, \fa}\to M^{\Kps,\red}_{\widehat{\mG}, \fa}$ does not change as $\fa$ varies in $A_i$.
\end{theorem}

To characterize the wall crossing phenomenon, we fix a semi-algebraic chamber $A_i$ in Theorem \ref{mainthm: kps fcd} and denote $\partial \bar{A_i}$ to be the boundary of the closure of $A_i$. Fix a vector $\fa_0:=(a_{01},...,a_{0k})\in \partial \bar{A_i}$ and a vector $\fa=(a_1,...,a_k)\in A_i$, which are not necessarily rational. Then we have the following:

\begin{theorem}\label{thm: wall crossing}{\rm (Theorem \ref{thm: diagram})}
Notation as above, there exists a wall crossing diagram as follows:
\begin{center}
	\begin{tikzcd}[column sep = 2em, row sep = 1.5em]
	 \mM^{\Kss,\red}_{\widehat{\mG}, \fa} \arrow[d,"\phi_\fa", swap] \arrow[rr,"\Psi_{\fa, \fa_0}"]&& \mM^{\Kss,\red}_{\widehat{\mG}, \fa_0} \arrow[d,"\phi_{\fa_0}",swap]\\
	 M^{\Kps,\red}_{\widehat{\mG}, \fa}\arrow[rr,"\psi_{\fa, \fa_0}"]&& M^{\Kps,\red}_{\widehat{\mG}, \fa_0}
	\end{tikzcd}
\end{center}
where 
$\psi_{\fa, \fa_0}$ is an induced proper morphism by the universal property of good moduli spaces. If the $\bQ$-span of $\fa_0$ (i.e. the smallest affine space defined over $\bQ$ containing $\fa_0$) contains the chamber $A_i$, 
then $\Psi_{\fa, \fa_0}$ is an open embedding.
\end{theorem}

As a corollary, in the case of one boundary divisor, we get the following finite polytope chamber decomposition to control the variation of K-semi(poly)stability. 

\begin{corollary}\label{maincor: kps fcd}{\rm (Corollary \ref{cor: fcd} and Corollary \ref{cor: kps fcd})}
Suppose $\widetilde{\mG}\subset \mG$ is a subset which is log bounded and satisfies the volume condition $(\clubsuit)$.  Suppose $k=1$ and $P=[a,b]\subset [0,1]$. Then there exist finitely many algebraic numbers (depending only on $\widetilde{\mG}$)
$$a=w_0<w_1<w_2<...<w_p<w_{p+1}=b $$
so that the K-semi(poly)stability of $(X, xD)$ does not change as $x$ varies in $(w_i, w_{i+1})$ for any $0\leq i\leq p$ and any $(X, D)\in \widetilde{\mG}$. Consequently,  we have the following wall crossing diagram for any $1\leq i\leq p$:

\begin{center}
	\begin{tikzcd}[column sep = 2em, row sep = 2em]
	\mM^{\Kss,\red}_{\widehat{\mG}, (w_{i-1}, w_i)} \arrow[d,"\phi_i^-",swap] \arrow[rr,"\Psi_i^-"]&& \mM^{\Kss,\red}_{\widehat{\mG}, w_i} \arrow[d,"\phi_i",swap] &&\mM^{\Kss,\red}_{\widehat{\mG}, (w_i, w_{i+1})}\arrow[d,"\phi_i^+",swap]\arrow[ll,"\Psi_i^+",swap]\\
	M^{\Kps,\red}_{\widehat{\mG}, (w_{i-1}, w_i)}\arrow[rr,"\psi_i^-"]&& M^{\Kps,\red}_{\widehat{\mG}, w_i}&&M^{\Kps,\red}_{\widehat{\mG}, (w_i, w_{i+1})}\arrow[ll,"\psi_i^+",swap].	 	
	\end{tikzcd}
\end{center}

\end{corollary}

We note that the walls $w_i$ can be irrational in the non-proportional case as opposed to the proportional case, where an example is provided in \cite{DJ+24}. Nevertheless,  $w_i$ is always an algebraic number.

\subsection{Properties of the K-semistable domain}\label{subsec: kss domain}
The strategy towards Theorem \ref{mainthm: kps fcd} is similar to \cite{Zhou23b}, i.e. we define the K-semistable domains associated to the couples in $\mG$ and study their properties, i.e. finiteness and good shape. 
For any $(X, \sum_{j=1}^kD_j)\in \mG$, we define the K-semistable domain (restricted to $P$) as follows:
$$\Kss(X, \sum_{j=1}^k D_j)_P:=\{(x_1,...,x_k)\in P\mid \textit{$(X, \sum_{j=1}^k x_jD_j)$ is a K-semistable log Fano pair}\}. $$
Note that the definition above of the K-semistable domain is also a little different compared with \cite{Zhou23b}, since we do not take the closure\footnote{Since we also take into account of $\bR$-coefficients (as have mentioned in the previous footnote), we do not take the closure of the domain. This means the K-semistable domain defined in this paper is the precise one.}. Under this definition, the K-semistable domain is not necessarily closed, but only locally closed under Euclidean topology (see Lemma \ref{lem: locally closed}). 
The first property of the K-semistable domain is the following finiteness criterion.

\begin{theorem}{\rm (= Theorem \ref{thm: finiteness})}\label{mainthm: finiteness}
Let $\pi: \mX\to T$ be a flat family of normal projective varieties over a smooth base $T$, and $\mD_j\to T$ is a family of effective Weil divisors\footnote{By a family of effective Weil divisors, we mean a relative effective Mumford $\bZ$-divisor in the sense of Koll\'ar, see \cite[Definition 4.68]{Kollar23}.} for each $1\leq j\leq k$.

Suppose for any closed point $t\in T$, there is a vector $(c_1,...,c_k)\in \bR_{\geq 0}^{k}$ such that $(\mX_t, \sum_{j=1}^kc_j\mD_{j,t})$ is a K-semistable log Fano pair. Define
$$\Kss(\mX_t, \sum_{j=1}^k \mD_{j,t}):=\{(x_1,...,x_k)\in \bR_{\geq 0}^k |\ \textit{$(\mX_t, \sum_{j=1}^k x_j\mD_{j,t})$ is a K-semistable log Fano pair}\}.$$
Then the following set is finite:
$$\{\Kss(\mX_t, \sum_{j=1}^k\mD_{j, t})\ |\ \text{$t\in T$}\}.$$
\end{theorem}

The above finiteness property enables us to cook up a finite chamber decomposition to control the variation of K-semistability as the coefficients vary. However, the shape of each chamber may be pathological (e.g. Corollary \ref{cor: chamber}).  In the proportional setting, each chamber turns out to be a rational polytope (e.g. \cite{Zhou23b}), which is quite satisfactory. But in non-proportional setting, there are  examples where K-semistable domains are not polytopes (e.g. \cite{Loginov23}). The following result gives a satisfactory characterization on the shape of the K-semistable domain.

\begin{theorem}\label{mainthm: shape}{\rm (= Corollary \ref{cor: semi-algebraic kss})}
Let $X$ be a normal projective variety of dimension $d$, and $D_j$'s ($1\leq j\leq k$) are effective Weil divisors on $X$. Given a rational polytope $Q\subset [0,1]^k$ and assume $Q\subset \overline{\LF(X, \sum_{j=1}^kD_j)}$ (see Section \ref{sec: constructibility}). Suppose there exists a positive real number $\epsilon$ such that $\vol(-K_X-\sum_{j=1}^kx_jD_j)\geq \epsilon$ for any $(x_1,...,x_k)\in Q$. 
Then $\Kss(X, \sum_{j=1}^kD_j)_Q$ is a semi-algebraic set. 
\end{theorem}

Theorem \ref{mainthm: shape} says that the K-semistable domain is semi-algebraic, which leads to the semi-algebraic chamber decomposition as in Theorem \ref{mainthm: kps fcd}.

\subsection{Constructible properties}\label{subsec: constructible}

In the proportional situation, the corresponding results in Section \ref{subsec: wall-crossing} and Section \ref{subsec: kss domain} are worked out in \cite{ADL19, Zhou23, Zhou23b}, where we could easily find a suitable family parametrizing the required log pairs.  However, it is much trickier in non-proportional setting, since a priori we need to show various constructible properties to select the required fibers in a family of couples. We overcome this difficulty by proving the following propositions.

\begin{proposition}{\rm (= Proposition \ref{prop: QC stratification})}\label{mainprop: QC}
Let $\pi: \mX\to T$ be a flat family of normal projective varieties over a smooth base $T$, and $\mD_j\to T$ is a family of effective Weil divisors for each $1\leq j\leq k$.
For each $t\in T$, we define the $\bR$-Cartier domain as follows:
$$\RC(\mX_t, \sum_{j=1}^k\mD_{j,t}):=\{(x_1,...,x_k)\in \bR_{\geq 0}^k\ |\ \text{$(\mX_t,\sum_{j=1}^kx_j\mD_{j,t})$ is a log pair}\}.$$
Suppose every fiber of $\pi$ is of Fano type, then the following map is constructible:
$$t\in T\mapsto \RC(\mX_t, \sum_{j=1}^k\mD_{j,t}). $$
More precisely, there is a constructible stratification of $T$, denoted by $T=\amalg_i T_i$, such that the fibers over each fixed piece $T_i$ admit the same $\bR$-Cartier domain.
\end{proposition} 

In Section \ref{sec: algebraic locus}, we will prove a little stronger result than Proposition \ref{mainprop: QC} by posing a weaker condition (see Lemma \ref{lem: family of RC}).

\begin{proposition}{\rm (= Proposition \ref{prop: LF stratification})}\label{mainprop: LF}
Let $\pi: \mX\to T$ be a flat family of normal projective varieties over a smooth base $T$, and $\mD_j\to T$ is a family of effective Weil divisors for each $1\leq j\leq k$. For each $t\in T$, we define the log Fano domain and the weak log Fano domain as follows:
$$\LF(\mX_t, \sum_{j=1}^k\mD_{j,t}):=\{(x_1,...,x_k)\in \bR_{\geq 0}^k\ |\ \text{$(\mX_t,\sum_{j=1}^kx_j\mD_{j,t})$ is a log Fano pair}\},$$
$$\WLF(\mX_t, \sum_{j=1}^k\mD_{j,t}):=\{(x_1,...,x_k)\in \bR_{\geq 0}^k\ |\ \text{$(\mX_t,\sum_{j=1}^kx_j\mD_{j,t})$ is a weak log Fano pair}\}.$$
Suppose for any closed point $t\in T$, there is a vector $(c_1,...,c_k)\in \bR_{\geq 0}^{k}$ such that $(\mX_t, \sum_{j=1}^kc_j\mD_{j,t})$ is a log Fano pair (resp. a weak log Fano pair). Then the following map is constructible:
$$t\in T\mapsto \LF(\mX_t, \sum_{j=1}^k\mD_{j,t}),\ \ \ ({\rm resp.}\ \  t\in T\mapsto \WLF(\mX_t, \sum_{j=1}^k\mD_{j,t})).$$
More precisely, there is a constructible stratification of $T$, denoted by $T=\amalg_i T_i$, such that the fibers over each fixed piece $T_i$ admit the same log Fano domain (resp. weak log Fano domain).
\end{proposition}

\begin{proposition}{\rm (= Proposition \ref{prop: algebraic locus})}\label{mainprop: locus}
Let $\pi: \mX\to T$ be a flat family of normal projective varieties over a smooth base $T$, and $\mD_j\to T$ is a family of effective Weil divisors for each $1\leq j\leq k$. Then the following two sets are algebraic subsets of $T$:
$$T_0:=\{t\in T|\text{$(\mX_t, \sum_{j=1}^kc_j\mD_{j,t})$ is log Fano for some numbers $c_j$'s}\}, $$
$$T'_0:=\{t\in T|\text{$(\mX_t, \sum_{j=1}^kc_j\mD_{j,t})$ is weak log Fano for some numbers $c_j$'s}\}. $$
\end{proposition} 

Armed by these propositions, we could find a suitable family of couples parametrizing couples in $\widetilde{\mG}$. This would be crucial to the proof of Theorem \ref{mainthm: finiteness} and Theorem \ref{mainthm: shape} (see Section \ref{sec: finiteness} and Section \ref{sec: shape}).

\subsection{GIT-stability vs K-stability}

As an application, in Section \ref{sec: GIT=K}, we will also establish the following comparison between K-stability and GIT-stability. It could be viewed as an analogue result of \cite[Theorem 1.1]{Zhou21a} in non-proportional setting. 

\begin{theorem}{\rm (Theorem \ref{thm: GIT=K})}\label{mainthm: GIT=K}
Let $X$ be a K-polystable Fano variety with  $\bQ$-factorial terminal singularities.
Fix $|L|$ to be a linear system of an effective line bundle $L$ on $X$. Suppose $(X, \epsilon D_0)$ is a K-semistable log Fano pair for some sufficiently small rational number $0<\epsilon\ll 1$ and some $D_0\in |L|$, then we have the following conclusion:  the K-(semi/poly)stability of $(X, c D)$ is identical to the GIT-(semi/poly)stability of $D\in |L|$ under $\Aut(X)$-action for any sufficiently small rational number $0<c\ll 1$.
\end{theorem}

Let $\mD\subset X\times |L|$ be the universal divisor and $ (X\times |L|, \mD)\to |L|$ the universal family. The condition in the above theorem just means that $(X\times |L|, \epsilon\mD)\to |L|$ admits a fiber which is a K-semistable log Fano pair for some sufficiently small rational number $0< \epsilon \ll 1$. Under this condition, Theorem \ref{mainthm: GIT=K} essentially tells that GIT-moduli coincides with K-moduli (see Remark \ref{rmk: GIT=K}).

\subsection{Conjectures}\label{subsec: conj}

The results on wall crossing in Section \ref{subsec: wall-crossing} are worked out under a volume condition $(\clubsuit)$ for $\widetilde{\mG}$. In fact, we mainly use this condition to derive various boundedness results, and it is essentially applied in the proof of Theorem \ref{mainthm: shape}, where the semi-algebraicity of the K-semistable domain is obtained (see Section \ref{sec: shape}). Thus the semi-algebraic chamber decomposition in Theorem \ref{mainthm: kps fcd} essentially relies on this condition.
However, as $\widetilde{\mG}$ is already assumed to be log bounded, one may naturally ask whether it still works without the volume condition $(\clubsuit)$. This leads to the following conjectures.

\begin{conjecture}\label{conj: 1}
Theorem \ref{mainthm: kps fcd} still holds without the volume condition $(\clubsuit)$.
\end{conjecture}

An important ingredient towards Theorem \ref{mainthm: kps fcd} without the volume condition $(\clubsuit)$ is the following conjectural characterization of the K-semistable domain:

\begin{conjecture}\label{conj: 1.5}
Let $X$ be a normal projective variety, and $D_j$'s ($1\leq j\leq k$) are effective Weil divisors on $X$. Then $\Kss(X, \sum_{j=1}^kD_j)$ is a semi-algebraic set. 
\end{conjecture}

In Theorem \ref{mainthm: shape}, Conjecture \ref{conj: 1.5} is confirmed when restricted to a rational polytope for which the volume condition $(\clubsuit)$ is satisfied. 
Without the volume condition $(\clubsuit)$, we immediately encounter a difficulty on boundedness of $\widehat{\mG}$. We give an example to explain the subtle point. Let $f: (\mX, \mD)\to T$ be a flat family of log smooth couples over a smooth base $T$, where $\mX$ is smooth, $\mD$ is a smooth Weil divisor on $\mX$, and $T$ is not necessarily proper.   Suppose $f$ is of relative dimension $d$ and $(\mX_t, c\mD_t)$ is a K-semistable log Fano pair for any $c\in(0,1)$ and any $t\in T$.  Take $P:=[0,1]$ and $\widetilde{\mG}$ to be the set of fibers of $f$. Then we clearly see $\widetilde{\mG}\subset \mG(d,1, P)$. Note that $v_0:=\vol(-K_{\mX_t}-\mD_t)$ does not depend on $t\in T$. If $v_0>0$ (i.e. $\widetilde{\mG}$ satisfies the volume condition $(\clubsuit)$), then $\widehat{\mG}$ is log bounded (e.g. Proposition \ref{prop: G-complete}). In this case, Corollary \ref{maincor: kps fcd} tells us that there are only finitely many compactifications of $T$ in the K-moduli spaces (see Example \ref{example: p1p2}). The subtlety appears when $v_0=0$.

\begin{conjecture}\label{conj: 2}
If $v_0:=\vol(-K_{\mX_t}-\mD_t)=0$, then $\widehat{\mG}$ is still log bounded.
\end{conjecture}

We note that all of the above conjectures are known in the proportional case, see \cite{Zhou23b}.

\subsection{The relation with \cite{LZ24}}

The readers may find that this article and \cite{LZ24} mutually refer to each other at several places, and then suspect whether the two works stand on their own. We make a clarification in this subsection. 

Going through \cite{LZ24}, one can easily make sure that \cite{LZ24} essentially only depend on Propositions \ref{prop: polytope QC}, \ref{prop: polytope LF}, \ref{prop: QC stratification}, \ref{prop: LF stratification}, \ref{prop: algebraic locus}, Proposition \ref{prop: N-complement}, Theorem \ref{thm: finiteness}, Propositions \ref{prop: Abdd}, \ref{prop: perfect fcd}, \ref{prop: uniform polarization}, Corollary \ref{prop: perfect tc}. Let us first take a look at Propositions \ref{prop: polytope QC}, \ref{prop: polytope LF}, \ref{prop: QC stratification}, \ref{prop: LF stratification}, \ref{prop: algebraic locus}. These results are all about the properties of Fano type varieties, which can be proved quite independently (of the both two works). We then look at Propositions \ref{prop: Abdd}, \ref{prop: perfect fcd}, \ref{prop: uniform polarization}, Corollary \ref{prop: perfect tc}. Proposition \ref{prop: Abdd} is about the boundedness of log Fano pairs or weak log Fano pairs; Propositions \ref{prop: perfect fcd}, \ref{prop: uniform polarization}, Corollary \ref{prop: perfect tc} are also properties of Fano type varieties, whose proofs only rely on Propositions \ref{prop: polytope QC}, \ref{prop: polytope LF}, \ref{prop: QC stratification}, \ref{prop: LF stratification}, \ref{prop: algebraic locus}. One may also note that some arguments in \cite[Proposition 4.1, Theorem 5.2]{LZ24} rely on the proofs of Proposition \ref{prop: N-complement} and Theorem \ref{thm: finiteness}. But the ideas of these proofs are actually standard and well-known to experts before, thus being quite independent of both two works. 

In summary, we conclude that \cite{LZ24} stands independently of this article and there is no circularity between the two works. Actually, the core techniques regarding wall crossing results in Section \ref{subsec: wall-crossing} have nothing to do with the core techniques of the construction of the K-moduli in \cite{LZ24}.

\subsection{Organization and convention}
In Section \ref{sec: preliminaries}, we recall some basic knowledge on K-stability and boundedness of log Fano pairs. In Section \ref{sec: constructibility} and Section \ref{sec: algebraic locus}, we prove several constructible properties (listed in Section \ref{subsec: constructible}) that will be applied through the article. In Section \ref{sec: finiteness}, we prove the finiteness criterion for K-semistable domains, i.e. Theorem \ref{mainthm: finiteness}. In Section \ref{sec: shape}, we prove the semi-algebraicity of the K-semistable domain, i.e. Theorem \ref{mainthm: shape}.  Section \ref{sec: bdd set} contains several useful boundedness results, which serve as the preparations for the final proof of the technical Theorem \ref{mainthm: shape} in Section \ref{sec: shape}. As applications, we prove the results (listed in Section \ref{subsec: wall-crossing}) on wall crossing in Section \ref{sec: fcd}, and Theorem \ref{mainthm: GIT=K} in Section \ref{sec: GIT=K}.
Throughout the article, we make the following convention:
\begin{itemize} 
\item a polytope and a face of a polytope are always convex and compact; 
\item a polyhedron $P\subset \bR^k$ is always closed but not necessarily bounded as the intersection of finitely many half spaces;
\item for a polytope $P\subset \bR^k$, we use $P^\circ$  to denote the interior part of $P$ under relative topology; for any subset $A\subset \bR^k$, we use $A(\bQ)$ to denote the set of rational points in $A$;
\item when we mention a point on the base of a family of varieties, it always means a closed point;
\item a Weil divisor means a $\bZ$-divisor which is not necessarily irreducible;
\item a constructible stratification of a variety $T$ means decomposing $T$ into a finite disjoint union where each piece is locally closed.
\end{itemize}

\noindent
\subsection*{Acknowledgement}
Part of this work was completed while CZ was visiting SCMS in Fudan university during Aug-Sep in 2023. CZ would like to thank Chen Jiang for the hospitality and helpful communications. 
YL is partially supported by NSF CAREER Grant DMS-2237139 and an AT\&T Research Fellowship from Northwestern University. CZ is 
partially supported by the NSFC grant (No. 12501058), and a grant from Xiamen University (No. X2450214), and Samsung Science and Technology Foundation under Project Number SSTF-BA2302-03.

\section{K-stability of Fano varieties and boundedness}\label{sec: preliminaries}

In this section, we recall some basic knowledge about K-stability and boundedness of Fano varieties. 

We say $(X, \Delta)$ is a \textit{couple} if $X$ is a normal projective variety and $\Delta$ is an effective $\bR$-divisor on $X$. 

We say a couple $(X, \Delta)$ is a \textit{log pair} if $K_X+\Delta$ is $\bR$-Cartier; if $\Delta$ is a $\bQ$-divisor, we say $(X, \Delta)$ is a $\bQ$-pair.

For a log pair $(X, \Delta)$, we say it is \textit{log Fano} (resp. \textit{weak log Fano}) if $(X, \Delta)$ admits klt singularities and $-K_X-\Delta$ is ample (resp. big and nef); in the case $\Delta=0$, we just say \textit{Fano} (resp. \textit{weak Fano}) instead of log Fano (resp. weak log Fano). 

We say $X$ is a \textit{Fano type} variety if there is an $\bR$-divisor $\Delta$ on $X$ such that $(X, \Delta)$ is a log Fano pair. 

For the definition of many types of singularities in minimal model program such as $klt$, $lc$, $plt$, $dlt$, etc., we refer to the standard textbooks \cite{KM98, Kollar13}.

After the works by Fujita and Li (\cite{Fuj19, Li17}), we can define the K-semistability of a log Fano pair as follows.
\begin{definition} (Fujita-Li)\label{def: Fujita-Li}
For a given weak log Fano pair $(X, \Delta)$\footnote{In the original works of K.Fujita and C. Li, they mainly talk about the K-semistability of log Fano pairs rather than weak log Fano pairs. However, one can just define the K-semistability of weak log Fano pairs as above. Up to passing to the anti-canonical model, one easily finds the compatibility between the K-semistability of a weak log Fano pair and that of its anti-canonical model.}, we say it is \textit{K-semistable} if 
$$\beta_{X, \Delta}(E):=A_{X, \Delta}(E)-S_{X, \Delta}(E)\geq 0$$ for any prime divisor $E$ over $X$, where 
$$A_{X, \Delta}(E):=\ord_{E}\left(K_Y-f^*(K_X+\Delta)\right)+1 $$
and
$$S_{X, \Delta}(E):=\frac{1}{\vol(-K_X-\Delta)}\int_0^{\infty}\vol(-f^*(K_X+\Delta)-tE){\rm d} t .$$
Here $f: Y\to (X, \Delta)$ is a proper log resolution such that $E$ is a prime divisor on $Y$.
\end{definition}

Given a weak log Fano pair $(X, \Delta)$, we define $\delta(X, \Delta):=\inf_{E}\frac{A_{X, \Delta}(E)}{S_{X, \Delta}(E)}$, where $E$ runs through all prime divisors over $X$. The following result is clear given Definition \ref{def: Fujita-Li}.

\begin{theorem}{\rm (\cite{Fuj19, Li17, FO18, BJ20})}
Let $(X, \Delta)$ be a weak log Fano pair. Then $(X, \Delta)$ is K-semistable if and only if $\delta(X, \Delta)\geq 1$.
\end{theorem}

\begin{remark}
By \cite{FO18, BJ20}, there is another way to reformulate $\delta(X, \Delta)$, i.e. to define it as the limit of a sequence of real numbers $\delta_m(X, \Delta)$, where $m$ runs through all sufficiently divisible positive integers. The number $\delta_m(X, \Delta)$ is defined as follows:
$$\delta_m(X, \Delta):=\inf_{E}\frac{A_{X, \Delta}(E)}{\sup_{D_m}\ord_{E}D_m}, $$
where $E$ runs through all prime divisors over $X$ and $D_m$ runs through all $m$-basis type divisors. Note that we say $D_m$ is an $m$-basis type divisor if it is of the following form
$$D_m=\frac{\sum_{j=1}^{N_m}{\rm div}(s_j)}{mN_m}\sim_\bR -(K_X+\Delta),  $$
where $N_m=\dim H^0(X,- m(K_X+\Delta))$ and $\{s_j\}_{j=1}^{N_m}$ is a basis of the vector space $H^0(X, - m(K_X+\Delta))$. Note the following
$$ H^0(X, - m(K_X+\Delta)):=H^0(X, \lfloor- m(K_X+\Delta)\rfloor)+\lceil m\Delta\rceil -m\Delta.$$
By the similar arguments as in \cite{BJ20, AZ22}, we still have $\lim_m \delta_m(X, \Delta)=\delta(X, \Delta)$ (e.g. \cite[Section 2.2]{LZ24}).
\end{remark}

We also recall the concept of K-polystability.

\begin{definition}{\rm (Test configuration)}\label{def: tc}
Let $(X,\Delta)$ be a log pair of dimension $d$ and $L$ an ample $\bR$-line bundle on $X$. A \textit{test configuration} (of $(X, \Delta; L)$) $\pi: (\mX,\Delta_\tc;\mL)\to \bA^1$ is a degenerating family over $\bA^1$ consisting of the following data:
\begin{enumerate}
\item $\pi: \mX\to \bA^1$ is a projective flat morphism of normal varieties, $\Delta_\tc$ is an effective $\bQ$-divisor on $\mX$, and $\mL$ is a relatively ample $\bR$-line bundle on $\mX$;
\item the family $\pi$ admits a $\bC^*$-action which lifts the natural $\bC^*$-action on $\bA^1$ such that $(\mX,\Delta_\tc; \mL)\times_{\bA^1}\bC^*$ is $\bC^*$-equivariantly isomorphic to $(X, \Delta; L)\times_{\bA^1}\bC^*$.
\end{enumerate}
We say the test configuration $(\mX, \Delta_\tc;\mL)$ is \textit{product type} if it is induced by a one parameter subgroup of $\Aut(X, \Delta; L)$. 

We denote $(\bar{\mX}, \bar{\Delta}_\tc;\bar{\mL})\to \bP^1$ to be the natural compactification of the original test configuration, which is obtained by glueing $(\mX, \Delta_\tc;\mL)$ and $(X,\Delta;L)\times (\bP^1\setminus 0)$ along their common open subset $(X, \Delta;L)\times \bC^*$.
\textit{The generalized Futaki invariant} associated to $(\mX, \Delta_\tc; \mL)$ is defined as follows:
$$\Fut(\mX,\Delta_\tc;\mL):=\frac{(K_{\bar{\mX}/\bP^1}+\bar{\Delta}_\tc).\bar{\mL}^d}{L^d} -\frac{d}{d+1}\frac{\left((K_X+\Delta)L^{d-1}\right)\bar{\mL}^{d+1}}{(L^d)^2}.$$
\end{definition}

\begin{definition}{\rm (Special test configuration)}\label{def: special tc}
Suppose $(X,\Delta)$ is a log Fano pair of dimension $d$ and $L=-K_X-\Delta$. Let  $(\mX,\Delta_\tc; \mL)$ be a test configuration of $(X, \Delta;L)$ such that $\mL=-K_{\mX/\bA^1}-\Delta_\tc$. We call it a \textit{special test configuration} if  $(\mX, \mX_0+\Delta_{\tc})$ admits plt singularities (or equivalently, the central fiber $(\mX_0, \Delta_{\tc,0})$ is a log Fano pair). The generalized Futaki invariant of a special test configuration is simplified as follows:
$$\Fut(\mX,\Delta_\tc;-K_\mX-\Delta_\tc)=-\frac{1}{d+1}\cdot\frac{(-K_{\bar{\mX}/\bP^1}-\bar{\Delta}_\tc)^{d+1}}{(-K_X-\Delta)^d}.$$
\end{definition}

The following definition is inspired by \cite{LWX21}.

\begin{definition}\label{def: kps}
Let $(X, \Delta)$ be a K-semistable log Fano pair. We say it is  \textit{K-polystable} if any special test configuration of $(X, \Delta; -K_X-\Delta)$ with a K-semistable central fiber is product type.
\end{definition}

These definitions (i.e. \ref{def: Fujita-Li}, \ref{def: tc}, \ref{def: special tc}, \ref{def: kps}) are standard in literature for log Fano pairs with $\bQ$-coefficients (e.g. \cite{Xu21}). However, in the case of $\bR$-coefficients, many standard theorems in algebraic K-stability theory (which are proven for $\bQ$-coefficients case) under these definitions need to be checked again. We list some of them below that will be applied in this article.

\begin{theorem}\label{thm: R-copy}
Given a log Fano pair $(X, \Delta)$. We have the following conclusions:
\begin{enumerate}
\item if $\delta(X, \Delta)\leq 1$, then there exists a prime divisor $E$ over $X$ such that $\delta(X, \Delta)=\frac{A_{X, \Delta}(E)}{S_{X, \Delta}(E)}$, and $E$ induces a special test configuration $(\mX, \Delta_{\tc}; -K_\mX-\Delta_{\tc})\to \bA^1$ of $(X, \Delta; -K_X-\Delta)$ such that 
$\delta(\mX_0, \Delta_{\tc, 0})=\delta(X, \Delta)$ and 
$\Fut(\mX, \Delta_{\tc}; -K_\mX-\Delta_{\tc})$ is equal to $A_{X, \Delta}(E)-S_{X, \Delta}(E)$ up to a positive multiple;
\item if $(X, \Delta)$ is K-semistable, then there exists a unique K-polystable degeneration via a special test configuration;
\item $(X, \Delta)$ being K-semistable is equivalent to the non-negativity of the generalized Futaki invariants of all special test configurations of $(X, \Delta; -K_X-\Delta)$.
\end{enumerate}
\end{theorem}

\begin{proof}
In the case of $\bQ$-coefficients, statement (1) is proved in \cite{LXZ22} (combined with \cite{BLZ22}); statement (2) is proved in \cite{LWX21}; statement (3) is proved in \cite{LX14}. 

In the case of $\bR$-coefficients, these statements are proved in \cite[Theorem 4.5 and Theorem 4.11]{LZ24}, \cite[Theorem 6.3]{LZ24}, \cite[Theorem 3.2]{LZ24}.
\end{proof}

Next we recall some results on boundedness of log Fano pairs. 

\begin{definition}
Let $\mP$ be a set of projective varieties of dimension $d$. We say $\mP$ is bounded if there exists a projective morphism $ \mY\to T$ between schemes of finite type such that for any $X\in \mP$, there exists a closed point $t\in T$ such that $X\cong \mY_t$. Let $\mP'$ be a set of couples of dimension $d$, we say 
$\mP'$ is log bounded if there exist a projective morphism $\mY\to T$ between schemes of finite type and a reduced Weil divisor $\mD$ on $\mY$ such that for any $(X, \Delta)\in \mP$, there exists a closed point $t\in T$ such that $(X, \red(\Delta))\cong (\mY_t, \mD_t)$. Here $\red(\Delta)$ means taking all the coefficients of components in $\Delta$ to be one. 
\end{definition}

\begin{theorem}{\rm (\cite[Theorem 1.1]{Birkar21})}\label{thm: BAB}
Fix a positive integer $d$ and positive real number $\epsilon$. The following set of algebraic varieties is bounded:
$$\{X\ | \ \text{$(X, \Delta)$ is an $\epsilon$-lc weak log Fano pair of dimension $d$ for some $\bR$-divisor $\Delta$ on $X$}\} .$$
\end{theorem}

Based on the above boundedness, Jiang (\cite{Jiang20}) proves the boundedness of K-semistable Fano varieties with fixed dimension and volume, and this is an important step to construct the K-moduli space (e.g. \cite{Xu21}).

\begin{theorem}{\rm (\cite{Jiang20})}\label{thm: bddkss}
Fix a positive integer $d$ and two positive real numbers $v$ and $a$. The following set of algebraic varieties is bounded:
$$\{X\ |\ \text{$X$ is a  Fano variety of dimension $d$ with $\vol(-K_X)\geq v$ and $\delta(X)\geq a$}\}. $$
\end{theorem}

The following result is a generalization of Theorem \ref{thm: bddkss} to log Fano pairs.
\begin{theorem}\label{cor: bdd}
Fix two positive real numbers $v$ and $a$. The following set of algebraic varieties is bounded:
$$\{X\ | \ \text{$(X, \Delta)$ is a log Fano pair of dimension $d$ with $\vol(-K_X-\Delta)\geq v$ and $\delta(X, \Delta)\geq a$}\}. $$
\end{theorem}

\begin{proof}
See \cite[Corollary 2.22]{LZ24}.
\end{proof}

The above boundedness result still holds if one replaces \textit{log Fano pair} with \textit{weak  log Fano pair} by applying \cite[Theorem 1.1]{Birkar21}.

\section{Constructible properties}\label{sec: constructibility}

We aim to prove several constructible properties in this section.

\begin{definition}
Let $X$ be a normal projective variety and $D_j, 1\leq j\leq k,$ are effective Weil divisors on $X$. We define several domains associated to the couple $(X, \sum_{j=1}^kD_j)$. The $\bR$-Cartier domain is defined as follows:
$$\RC(X, \sum_{j=1}^kD_j):=\{(x_1,...,x_k)\in \bR_{\geq 0}^k\ |\ \text{$(X,\sum_{j=1}^kx_jD_j)$ is a log pair}\};$$
The log Fano domain is defined as follows:
$$\LF(X, \sum_{j=1}^kD_j):=\{(x_1,...,x_k)\in \bR_{\geq 0}^k\ |\ \text{$(X, \sum_{j=1}^kx_jD_j)$ is a log Fano pair}\}; $$
The weak log Fano domain is defined as follows:
$$\WLF(X, \sum_{j=1}^kD_j):=\{(x_1,...,x_k)\in \bR_{\geq 0}^k\ |\ \text{$(X, \sum_{j=1}^kx_jD_j)$ is a weak log Fano pair}\}. $$
\end{definition}

It is clear that $\LF(X, \sum_{j=1}^kD_j)$ and $\WLF(X, \sum_{j=1}^kD_j)$ are both convex and not necessarily closed in $\bR_{\geq 0}^k$. We denote
$\overline{\LF(X, \sum_{j=1}^kD_j)}$ and $\overline{\WLF(X, \sum_{j=1}^kD_j)}$ to be their closures.

\begin{proposition}\label{prop: polytope QC}
Let $X$ be a normal projective variety and $D_j, 1\leq j\leq k,$ are effective Weil divisors on $X$. Then $\RC(X, \sum_{j=1}^kD_j)$ is a rational polyhedron in $\bR^k_{\geq 0}$.
\end{proposition}
\begin{proof}
We denote $V$ to be the $k+1$-dimensional vector space (defined over $\bR$) generated by $k+1$ divisors: $K_X, D_1,...,D_k$. In other words, an element in $V$ is of the form $x_0K_X+\sum_{j=1}^kx_jD_j$, which is determined by a vector $(x_0,x_1,...,x_k)\in \bR^{k+1}$. Then there is a linear subspace $V'\subset V$ defined over $\bQ$ such that $V'$ consists of $\bR$-Cartier divisors of the form $x_0K_X+\sum_{j=1}^k x_jD_j$. It is then not hard to see the following formulation:
$$\RC(X, \sum_{j=1}^k D_j)= \{(x_1,...,x_k)\mid\text{$(x_0,x_1,...,x_k)\in V'\cap \{x_0=1\}\cap_{j=1}^k\{x_j\geq 0\}$}\}, $$
which is a rational polyhedron in $\bR^k_{\geq 0}$.
\end{proof}

If $\LF(X, \sum_{j=1}^kD_j)$ is non-empty, then $X$ is a Fano type variety. In this case, it is well-known by \cite{BCHM10} that $X$ is a Mori dream space. More precisely, $X$ satisfies the conditions that $H^1(X, \mO_X)=0$ and the Cox ring is finitely generated. In particular, the nef cone of $X$, denoted by $\Nef(X)$, is a rational polytope (e.g. \cite{HK00}). Note that $\Nef(X)$ is defined as a cone in N\'eron-Severi  group $\NS(X)_\bR$ generated by nef divisors. Here $\NS(X)_\bR$ is the $\bR$-vector space generated by Cartier divisors modulo numerical equivalences.
The following proposition is well-known to experts by \cite{BCHM10}.

\begin{proposition}\label{prop: polytope LF}
Let $X$ be a normal projective variety and $D_j, 1\leq j\leq k,$ are  effective Weil divisors on $X$. Then $\overline{\LF(X, \sum_{j=1}^kD_j)}$ and $\overline{\WLF(X, \sum_{j=1}^kD_j)}$ are both rational polytopes in $\bR_{\geq 0}^k$.
\end{proposition}

\begin{proof}
We first show that $\overline{\WLF(X, \sum_{j=1}^kD_j)}$ is a rational polytope. Suppose it is not empty. Write
$$P_1:= \RC(X, \sum_{j=1}^k D_j). $$
Let $X'\to X$ be a small $\bQ$-factorization of $X$, and $D_j'$ the birational transformation of $D_j$. We define the following domain
$$P_2:=\{(x_1,...,x_k)\in \bR^k_{\geq 0}\ |\ \text{$(X', \sum_{j=1}^kx_jD'_j)$ is log canonical}\}. $$
It is clear that $P_2$ is a rational polytope. Let $\NS(X')_\bR$ be the N\'eron-Severi group of $X'$,  we define the following map:
$$\phi: P_2\to \NS(X')_\bR, \ \ \ \ \ \ (x_1,...,x_k)\mapsto -K_{X'}-\sum_{j=1}^kx_jD'_j. $$
Note that $X'$ is a Mori dream space by \cite{BCHM10}. Thus the nef cone $\Nef(X')$ is a rational polytope in $\NS(X')_\bR$ and $\phi^{-1}\Nef(X')$ is a rational polytope in $P_2$. Therefore
$$\overline{\WLF(X, \sum_{j=1}^kD_j)}=P_1\cap \phi^{-1}\Nef(X') $$
is a rational polytope in $\bR_{\geq 0}^k$.

If $\LF(X, \sum_{j=1}^kD_j)$ is not empty, then $\overline{\LF(X, \sum_{j=1}^kD_j)}$ is also a rational polytope as it coincides with $\overline{\WLF(X, \sum_{j=1}^kD_j)}$. 
The proof is complete.
\end{proof}

It is a simple observation that $\overline{\LF(X, \sum_{j=1}^kD_j)}\setminus \LF(X, \sum_{j=1}^kD_j)$ is an union of some faces of $\overline{\LF(X, \sum_{j=1}^kD_j)}$, and the weak log Fano cone also admits the same property.

\begin{proposition}\label{prop: QC stratification}
Let $\pi: \mX\to T$ be a flat family of normal projective varieties over a smooth base $T$, and $\mD_j\to T$ is a family of effective Weil divisors for each $1\leq j\leq k$. 
Suppose every fiber of $\pi$ is of Fano type. Then the following map is constructible:
$$t\in T\mapsto \RC(\mX_t, \sum_{j=1}^k\mD_{j,t}). $$
More precisely, there is a constructible stratification of $T$, denoted by $T=\amalg_i T_i$, such that the fibers over each fixed piece $T_i$ admit the same $\bR$-Cartier domain.
\end{proposition}

\begin{proof}
Up to a constructible stratification of $T$, we may assume $(\mX, \sum_{j=1}^k \mD_j)\to T$ admits a fiberwise log resolution, i.e., there exists a proper birational morphism over $T$:
$$f: (\mY, \sum_{j=1}^k\widetilde{\mD}_j+\sum_l \mE_l)\to (\mX, \sum_{j=1}^k\mD_j), $$
where $\widetilde{\mD}_j$ is the birational transformation of $\mD_j$, $\mE_l$ is the exceptional divisors of $f$, and $(\mY, \sum_{j=1}^k\widetilde{\mD}_j+\sum_l \mE_l)$ is log smooth over $T$. 

Since $R^2\pi_*\bZ$ is a constructible sheaf over $T$, up to a further constructible stratification, we may assume $R^2\pi_*\bZ$ is a locally constant local system under the analytic topology on $T$. Note that every fiber of $\pi$ is of Fano type,  we have following conclusions:
\begin{enumerate}
\item $\pi\circ f: \mY\to T$ is a family of smooth rationally connected varieties (see \cite[Corollary 1.5]{HM07}); 
\item $\Pic(\mX_t)=H^2(\mX_t, \bZ)$ and $\Pic(\mY_t)\cong H^2(\mY_t, \bZ)$ for any $t\in T$;
\item Up to a finite base change of $T$, the Picard functors $Pic_{\mX/T}$ and $Pic_{\mY/T}$ are constant sheaves on $T$ (e.g. \cite[Theorem 1.3]{CLZ25}).
\end{enumerate}
Fix a closed point $t_0\in T$. Suppose $K_{\mX_{t_0}}+\sum_{j=1}^kx_j\mD_{j,t_0}$ is $\bQ$-Cartier and $m(K_{\mX_{t_0}}+\sum_{j=1}^kx_j\mD_{j,t_0})$ is Cartier for some integer $m\in \bZ^+$. Then 
$$m\left(K_{\mY_{t_0}}+\sum_{j=1}^kx_j\widetilde{\mD}_{j,t_0}+\sum_l (a_l-1) \mE_{l, t_0}\right)=f^*m(K_{\mX_{t_0}}+\sum_{j=1}^kx_j\mD_{j,t_0}), $$
where $a_l=A_{\mX_{t_0}, \sum_{j=1}^kx_j\mD_{j, t_0}}(\mE_{l,t_0})$. By the conclusions listed above, we see that 
$$m\left(K_{\mY}+\sum_{j=1}^kx_j\widetilde{\mD}_{j}+\sum_l (a_l-1) \mE_{l}\right)=f^*\mL$$
for some Cartier divisor $\mL$ on $\mX$. Therefore,
$$\mL=f_*f^*\mL=f_*m\left(K_{\mY}+\sum_{j=1}^kx_j\widetilde{\mD}_{j}+\sum_l (a_l-1) \mE_{l}\right)=m(K_\mX+\sum_{j=1}^kx_j\mD_j). $$
This implies that $(\mX_t, \sum_{j=1}^kx_j\mD_{j,t})$ is a log pair for any $t\in T$ and $\RC(\mX_t, \sum_{j=1}^k\mD_{j,t})$ does not depend on $t\in T$. The proof is complete.
\end{proof}

\begin{proposition}\label{prop: LF stratification}
Let $\pi: \mX\to T$ be a flat family of normal projective varieties over a smooth base $T$, and $\mD_j\to T$ is a family of effective Weil divisors for each $1\leq j\leq k$.

Suppose for any closed point $t\in T$, there is a vector $(c_1,...,c_k)\in \bR_{\geq 0}^{k}$ such that $(\mX_t, \sum_{j=1}^kc_j\mD_{j,t})$ is a log Fano pair (resp. a weak log Fano pair). Then the following map is constructible:
$$t\in T\mapsto \LF(\mX_t, \sum_{j=1}^k\mD_{j,t}) \ \ \ \left({\rm resp.}\ \  t\in T\mapsto \WLF(\mX_t, \sum_{j=1}^k\mD_{j,t})\right). $$
More precisely, there is a constructible stratification of $T$, denoted by $T=\amalg_i T_i$, such that the fibers over each fixed piece $T_i$ admit the same log Fano domain (resp. weak log Fano domain).
\end{proposition}

\begin{proof}
We first prove the result for the log Fano domain. Up to a constructible stratification of $T$, we may assume $(\mX, \sum_{j=1}^k \mD_j)\to T$ admits a fiberwise log resolution. By Proposition \ref{prop: QC stratification}, we may assume 
$$\RC(\mX_t, \sum_{j=1}^k\mD_{j,t})$$ 
does not depend on $t\in T$ up to a further constructible stratification. Thus there exists a rational polytope $\tilde{P}$ in $\bR^k_{\geq 0}$ such that $(\mX_t, \sum_{j=1}^k x_j\mD_{j,t})$ is a log canonical pair for any $(x_1,...,x_k)\in \tilde{P}$ and any $t\in T$. Here $\tilde{P}$ is the log canonical polytope. 
Let $(c_1,...,c_k)$ be a vertex of $\tilde{P}$, which is a rational vector since $\tilde{P}$ is rational, then 
$$(\mX, \sum_{j=1}^kc_j\mD_j)\to T$$ 
is a family of log canonical pairs. Up to shrinking the base, we may assume $K_{\mX/T}+\sum_{j=1}^kc_j\mD_j$ is $\bQ$-Cartier (e.g. \cite[Lemma 29]{MST20}). Since there are only finitely many vertices, up to further shrinking the base we may assume $-K_{\mX}-\sum_{j=1}^k x_j \mD_j$ is $\bR$-Cartier for any $(x_1,...,x_k)\in \tilde{P}$. Up to shrinking the base again and by the openness of ampleness, we may assume $(\mX, \sum_{j=1}^ka_j\mD_j)\to T$ is a family of log Fano pairs for some $(a_1,...,a_k)\in \tilde{P}$.

Fix a closed point $t_0\in T$. Let $(c_1,...,c_k)$ be a vertex of $\overline{\LF(\mX_{t_0}, \sum_{j=1}^kc_j\mD_{j,t_0})}$, which is a rational vector. Then $-K_{\mX_{t_0}}-\sum_{j=1}^kc_j\mD_{j,t_0}$ is nef. By Lemma \ref{lem: open nef}, there exists an open neighborhood $t_0\in U_{t_0}$
such that $-K_{\mX_{t}}-\sum_{j=1}^kc_j\mD_{j,t}$ is nef for any $t\in U_{t_0}$. Thus we see that the map 
$$t\mapsto \overline{\LF(\mX_t, \sum_{j=1}^k\mD_{j,t})}$$ 
is lower semi-continuous, i.e., given a polytope $P\subset \tilde{P}$, the set
$$\{t\in T\mid  \text{$P\subset \overline{\LF(\mX_t, \sum_{j=1}^k\mD_{j,t})}$}\} $$
is an open subset of $T$. 

Next we show that there exists an open subset $U\subset T$ such that $\overline{\LF(\mX_t, \sum_{j=1}^k\mD_{j,t})}$ does not depend on $t\in U$, which implies that the map $t\in T\mapsto \overline{\LF(\mX_t, \sum_{j=1}^k\mD_{j,t})}$ is constructible (by Noetherian induction).  We apply the induction on dimension of $T$. First we assume $\dim T=1$. Fixing a closed point $t_1\in T$ and writing 
$P_1:=\overline{\LF(\mX_{t_1}, \sum_{j=1}^k\mD_{j,t_1})},$
by lower semi-continuity as above, there exists an open subset $t_1\in U_1\subset T$ such that $P_1\subset \overline{\LF(\mX_t, \sum_{j=1}^k\mD_{j,t})}$ for any $t\in U_1$. Suppose there exists a closed point $t_2\in U_1$ such that the inclusion 
$$P_1\subset P_2:=\overline{\LF(\mX_{t_2}, \sum_{j=1}^k\mD_{j,t_2})}$$ 
is strict, then there is an open subset $t_2\in U_2$ strictly contained in $U_1$ such that $P_2\subset \overline{\LF(\mX_{t}, \sum_{j=1}^k\mD_{j,t})}$ for any $t\in U_2$. Continuing the process, we obtain  a sequence of open subsets of $T$, denoted by
$$ U_1\supset U_2\supset...\supset U_m\supset... ,$$
such that each inclusion is strict, and there exists a sequence of closed points $t_1,t_2,...,t_m,...$ satisfying the following conditions:
\begin{enumerate}
\item $t_i\in U_i$ and $t_i\notin U_{i+1}$;
\item writing $P_i:=\overline{\LF(\mX_{t_i}, \sum_{j=1}^k\mD_{j,t_i})}$, then $P_i\subset \overline{\LF(\mX_t, \sum_{j=1}^k\mD_{j,t})}$ for any $t\in U_i$;
\item 
$P_1\subset P_2 \subset...\subset P_m\subset..., $
where each inclusion is strict.
\end{enumerate}
It suffices to show the above process terminates after finitely many steps. Suppose not, we aim to derive a contradiction.
By Proposition \ref{prop: polytope LF}, every $P_i$ is a rational polytope. Since $P_i\subset \tilde{P}$, the following set
$$P_\infty:= \overline{\{\cup_{i\geq 1} P_i\}}$$
is a closed convex set in $\tilde{P}$. Write $U_\infty:=\cap_{i\geq 1} U_i$. We first show the following equality
$$P_\infty=\overline{\LF(\mX_{t}, \sum_{j=1}^k\mD_{j,t})}$$
for any $t\in U_\infty$. By the second condition listed above, we see that 
$$\cup_i P_i\subset\overline{\LF(\mX_{t}, \sum_{j=1}^k\mD_{j,t})}.$$ 
This implies that $P_\infty\subset\overline{\LF(\mX_{t}, \sum_{j=1}^k\mD_{j,t})}$ since the right hand side is closed. We show that the inclusion is not strict. Suppose the inclusion is strict, then there exists some open neighborhood $t\in V$ such that the inclusion $P_\infty\subset \overline{\LF(\mX_{\tau}, \sum_{j=1}^k\mD_{j,\tau})}$ is strict for any $\tau\in V$. Recall that we assume $\dim T=1$, thus we see $\{t_i\}$ is dense in $T$ and $V\cap \{t_i\}_i$ is non-empty. This implies that the inclusion $P_\infty\subset P_i$ is strict for some $i$, contradiction. By induction, we assume there always exists an open subset $U\subset T$ such that $\overline{\LF(\mX_t, \sum_{j=1}^k\mD_{j,t})}$ does not depend on $t\in U$ if $\dim T\leq d$. Now we assume $\dim T=d+1$. Let $\{Z_l\}$ be countable smooth hypersurfaces on $T$ such that $\cup Z_l$ is dense in $T$. Let $V_l\subset T$ be the open subset with $V_l\cap Z_l$ non-empty such that 
$$P_l':=\overline{\LF(\mX_t, \sum_{j=1}^k\mD_{j,t})}$$
does not depend on $t\in V_l\cap Z_l$ (by induction). Up to shrinking $V_l$, we may assume 
$$P_l'\subset \overline{\LF(\mX_t, \sum_{j=1}^k\mD_{j,t})}$$
for any $t\in V_l$. Denote by $P'_\infty:=\overline{\{\cup_l P'_l\}}$. We show $P'_\infty$ is achieved by any $t\in \cap_l V_l$. We first see $P_\infty'\subset \overline{\LF(\mX_t, \sum_{j=1}^k\mD_{j,t})}$ similar as before. Suppose the inclusion is strict, then there exists an open subset $t\in V'$ such that the inclusion $P_\infty'\subset \overline{\LF(\mX_\tau, \sum_{j=1}^k\mD_{j,\tau})}$ is strict for any $\tau\in V'$. Note that $V'\cap Z_k$ non-empty for some $k$ (by the density of $\cup_l Z_l$ in $T$), we see that the inclusion $P_\infty'\subset P_k'$ is strict, contradiction. The contraction implies that $P'_\infty=\overline{\LF(\mX_t, \sum_{j=1}^k\mD_{j,t})}$ for any $t\in \cap V_l$. This clearly implies $P'_\infty=\overline{\LF(\mX_t, \sum_{j=1}^k\mD_{j,t})}$ for general $t\in T$ (rather than very general $t\in T$). Thus there exists an open subset $U\subset T$ such that $\overline{\LF(\mX_t, \sum_{j=1}^k\mD_{j,t})}$ does not depend on $t\in U$. By now, we complete the induction and confirm that general fibers admit the same closure of the log Fano domain. 
We turn to the map $t\in T\mapsto \LF(\mX_t, \sum_{j=1}^k\mD_{j,t}). $

Up to a constructible stratification, we assume the fibers of $(\mX, \sum_{j=1}^k\mD_j)\to T$ share the same closure of the log Fano domain, denoted by $Q$. 
Up to shrinking the base we may assume $-K_{\mX}-\sum_{j=1}^k x_j \mD_j$ is $\bR$-Cartier for any $(x_1,...,x_k)\in Q$. Let $\{F_l\}_{l=1}^m$ be the set of all proper faces of $Q$. For any $t\in T$, we have the following formulation
$$\LF(\mX_t, \sum_{j=1}^k\mD_{j,t})=Q^\circ\ \cup\ \cup_{i=1}^s F^\circ_{r_i},  $$
where $\{r_1, r_2,..., r_s\}$ is a subset of $\{1,2,...,m\}$ depending on the choice of $t\in T$. Thus the following set is finite:
$$S:=\{\LF(\mX_t, \sum_{j=1}^k\mD_{j,t})| \textit{$t\in T$}\}. $$
Suppose the log Fano domain $\LF(\mX_{t_0}, \sum_{j=1}^k\mD_{j,t_0})$ associated to $t_0\in T$ is maximal under the inclusion among elements of $S$. By the openness of ampleness,  it is not hard to see that there exists an open subset $t_0\in U$ such that 
$$\LF(\mX_{t_0}, \sum_{j=1}^k\mD_{j,t_0})\subset \LF(\mX_{t}, \sum_{j=1}^k\mD_{j,t})$$
for any $t\in U$, which implies the following equality (by the maximal condition on $\LF(\mX_{t_0}, \sum_{j=1}^k\mD_{j,t_0})$):
$$\LF(\mX_{t_0}, \sum_{j=1}^k\mD_{j,t_0})= \LF(\mX_{t}, \sum_{j=1}^k\mD_{j,t}).$$
By Noether induction, we see the following map is also constructible:
$$t\in T\mapsto \LF(\mX_t, \sum_{j=1}^k\mD_{j,t}).$$

The corresponding statement for the weak log Fano domain is obtained by replacing the openness of ampleness (as applied in the previous paragraph) with the openness of weak log Fano condition (e.g. \cite[Lemma 3.9]{dFH11} and \cite[Lemma 2.15]{zhuang19}).
\end{proof}

\begin{lemma}\label{lem: open nef}
Let $f: (X, \Delta)\to T$ be a family of log Fano pairs over a smooth base $T$ such that $X\to T$ is flat and $K_{X/T}+\Delta$ is $\bR$-Cartier. Let $D$ be an $\bR$-Cartier $\bR$-divisor on $X$. Then the following set is Zariski open:
$$\{t\in T\ |\ \textit{$D_t$ is nef\ }\}. $$
\end{lemma}

\begin{proof}
Suppose $\{t\in T\ |\ \textit{$D_t$ is nef\ }\}$ is not empty. By the openness of ampleness, we see that $D_t$ is nef for very general $t\in T$. We claim that $D_t$ is nef for general $t\in T$ (not just for very general $t\in T$). To see this, we run a $D$-MMP over $T$ on $X$.
Since $X$ is Fano type over $T$, we may assume the MMP sequence terminates with a good minimal model $\phi: (X, D)\dashrightarrow (X', D')/T$, where $D':=\phi_*D$ is semi-ample over $T$. It is clear that the extremal rays contracted in the course of MMP do not intersect very general fibers. Therefore they do not intersect the generic fiber of $X\to T$, which implies that $X\dashrightarrow X'/T$ is an isomorphism over the generic point of $T$. Thus we see that $D_t$ is nef for general $t\in T$ and we prove the claim. The claim says that $\{t\in T\ |\ \textit{$D_t$ is nef\ }\}$ contains an open subset (which is not empty) $U\subset T$. 

Denote $T_1:=T\setminus U$. We may write $T_1$ as a disjoint union of smooth varieties, denoted by $T_1=\coprod_i Z_i$. Let $(X_i, D_i)\to Z_i$ be the pullback of $(X, D)\to T$ along $Z_i\to T$. If $D_{i,t}$ is nef for some $t\in Z_i$, by the same argument as before, we see that $\{t\in Z_i\ |\ \textit{$D_{i,t}$ is nef\ }\}$ contains an open subset (which is not empty) of $Z_i$. By Noether induction and \cite[Chapter 2, Exercise 3.18 (c)]{Har77}, the set $\{t\in T\ |\ \textit{$D_t$ is nef\ }\}$ is Zariski open. The proof is complete.
\end{proof}

\section{Algebraic log Fano locus}\label{sec: algebraic locus}

Given a family of couples $\pi: (\mX, \sum_{j=1}^k\mD_j)\to T$, i.e. every fiber $(\mX_t, \sum_{j=1}^k\mD_{j,t})$ is a couple, we want to select the fibers satisfying the following condition:
$$\textit{there exist numbers $c_j\in \bR_{\geq 0}$ such that $(\mX_t, \sum_{j=1}^kc_j\mD_{j,t})$ is a log Fano pair}. $$
Similarly, we also want to select the fibers satisfying the following condition:
$$\textit{there exist numbers $c_j\in \bR_{\geq 0}$ such that $(\mX_t, \sum_{j=1}^kc_j\mD_{j,t})$ is a weak log Fano pair}. $$

In this section, we prove the above conditions are in fact algebraic. More precisely, we have the following result.

\begin{proposition}\label{prop: algebraic locus}
Let $\pi: \mX\to T$ be a flat family of normal projective varieties over a smooth base $T$, and $\mD_j\to T$ is a family of effective Weil divisors for each $1\leq j\leq k$. Then the following two sets are algebraic subsets of $T$:
$$T_0:=\{t\in T|\text{$(\mX_t, \sum_{j=1}^kc_j\mD_{j,t})$ is log Fano for some numbers $c_j$'s}\}, $$
$$T'_0:=\{t\in T|\text{$(\mX_t, \sum_{j=1}^kc_j\mD_{j,t})$ is weak log Fano for some numbers $c_j$'s}\}. $$
\end{proposition}

\begin{proof}
We first show that $T_0$ is algebraic.  Assume $T_0$ is non-empty, otherwise there is nothing to prove. We may also assume $T_0$ is dense in $T$ by taking the closure of $T_0$ in the sense of Zariski topology. Up to a constructible stratification of $T$, we assume $(\mX, \sum_{j=1}^k\mD_j)\to T$ admits a fiberwise log resolution.

By Lemma \ref{lem: family of RC}, up to a further constructible stratification, we may assume 
$$\RC(\mX_t, \sum_{j=1}^k\mD_{j,t})$$ 
does not depend on $t\in T$. By the first paragraph in the proof of Prop \ref{prop: LF stratification}, there exists a rational polytope $P$ consisting of all vectors $(x_1,...,x_k)$ such that $(\mX_t, \sum_{j=1}^kx_j\mD_{j,t})$ is log canonical, and up to shrinking $T$ we could assume $K_{\mX/T}+\sum_{j=1}^kx_j\mD_{j}$ is $\bR$-Cartier for any $(x_1,...,x_k)\in P$.
Fix $t_0\in T_0$.  By Prop \ref{prop: polytope LF}, there exists a rational vector $(c_1,...,c_k)\in P$ such that $(\mX_{t_0}, \sum_{j=1}^kc_j\mD_{j,t_0})$ is log Fano. By the openness of ampleness, there exists an open subset $t_0\in U_0$ such that $(\mX_{t}, \sum_{j=1}^kc_j\mD_{j,t})$ is also log Fano for any $t\in U_0$. This implies $U_0\subset T_0$. Hence $T_0$ is an algebraic subset of $T$.

For $T_0'$, we just replace the openness of ampleness applied in the second paragraph with the openness of weak log Fano condition as applied in \cite[Lemma 3.9]{dFH11} (see also \cite[Lemma 2.15]{zhuang19}).
\end{proof}

\begin{lemma}\label{lem: family of RC}
Let $\pi: \mX\to T$ be a flat family of normal projective varieties over a smooth base $T$, and $\mD_j\to T$ is a family of effective Weil divisors for each $1\leq j\leq k$. Suppose the set of closed points which correspond to Fano type fibers is dense in $T$ under Zariski topology, then the following map is constructible:
$$t\mapsto \RC(\mX_t, \sum_{j=1}^k\mD_{j,t}).$$ 
More precisely, there is a constructible stratification of $T$, denoted by $T=\amalg_i T_i$, such that the fibers over each fixed piece $T_i$ admit the same $\bR$-Cartier domain.
\end{lemma}

In Proposition \ref{prop: QC stratification}, we assume all fibers are Fano type, while in Lemma \ref{lem: family of RC} we pose a weaker condition that Fano type fibers are dense.

\begin{proof}
Up to a constructible stratification of $T$, we may assume $(\mX, \sum_{j=1}^k \mD_j)\to T$ admits a fiberwise log resolution, i.e., there exists a morphism over $T$:
$$f: (\mY, \sum_{j=1}^k\widetilde{\mD}_j+\sum_l \mE_l)\to (\mX, \sum_{j=1}^k\mD_j), $$
where $\widetilde{\mD}_j$ is the birational transformation of $\mD_j$, and $(\mY, \sum_{j=1}^k\widetilde{\mD}_j+\sum_l \mE_l)$ is log smooth over $T$. By deformation invariance of rationally connectedness (e.g. \cite[Theorem IV 3.11]{Kollar96}), the induced morphism $\pi\circ f: \mY\to T$ is a family of rationally connected projective manifolds. Thus for any $t\in T$ we have the following by \cite[Corollary IV 3.8]{Kollar96}
$$H^1(\mY_t, \mO_{\mY_t})= H^2(\mY_t, \mO_{\mY_t})=0.$$
Since Fano type fibers all admit rational singularities, by deformation invariance of rational singularities (e.g. \cite{Elkik78}), we may assume all fibers of $\pi$ have rational singularities (up to shrinking the base). Applying Leray spectral sequence we obtain
$$H^1(\mX_t, \mO_{\mX_t})= H^2(\mX_t, \mO_{\mX_t})=0 $$
for any $t\in T$. This implies $\Pic(\mX_t)=H^2(\mX_t, \bZ)$ for any $t\in T$. The rest is given by the similar argument as in the proof of Proposition \ref{prop: QC stratification}, since we only rely on $\Pic(\mX_t)=H^2(\mX_t, \bZ)$ therein.
\end{proof}

\section{Finiteness of K-semistable domains}\label{sec: finiteness}

In this section, we present a criterion to show the finiteness of K-semistable domains. This approach is inspired by \cite{Zhou23b}.

Let $X$ be a Fano type variety, then $-K_X$ is big. We say a log pair $(X, B)$ is an $N$-complement of $K_X$ if $(X, B)$ has lc singularities with $N(K_X+B)\sim 0$. 
The following result may be well-known to experts, e.g. \cite{BLX22, Zhou23b}.

\begin{proposition}\label{prop: N-complement}
Fixing a positive integer $d$, then there exists a sufficiently divisible positive integer $N$ depending only on $d$ such that the following statement holds.

Let $X$ be a Fano type variety of dimension $d$. Suppose there exists an effective $\bR$-divisor $\Delta$ on $X$ such that $(X, \Delta)$ is a weak log Fano pair with $\delta(X, \Delta)\leq 1$, then $\delta(X, \Delta)$ can be approximated by lc places of $N$-complements of $K_X$. More precisely, $\delta(X, \Delta)=\inf_{E} \frac{A_{X, \Delta}(E)}{S_{X, \Delta}(E)}$, where $E$ runs over all lc places of $N$-complements of $K_X$.
\end{proposition}

\begin{proof}
The idea is the same as \cite[Lemma 6.1]{Zhou23b}. We provide a complete proof for the readers' convenience. We may replace $(X, \Delta)$ with its small $\bQ$-factorization.

We first assume $\delta(X,\Delta)<1$, then $\delta_m:=\delta_m(X, \Delta)<1$ for sufficiently large $m$. In this case one can find an $m$-basis type divisor $B_m$ such that $(X,\Delta+ \delta_mB_m)$
is strictly log canonical and admits an lc place $E_m$ over $X$. 
By \cite[Corollary 1.4.3]{BCHM10}, there is an extraction $g_m: Y_m\to X$ which only extracts $E_m$, i.e.
 \begin{align*}
 K_{Y_m}+\tilde{\Delta}+\delta_m\tilde{B}_m+E_m 
 = g_m^*\left(K_X+\Delta+ \delta_mB_m\right),
 \end{align*}
 where $\tilde{\bullet}$ is the birational transformation of $\bullet$. Note that $Y_m$ is of Fano type. We could run $-(K_{Y_m}+E_m)$-MMP to get a mimimal model $f_m: Y_m\dashrightarrow Y_m'$ such that $-(K_{Y_m}+{f_m}_*E_m)$ is nef and $Y_m'$ is also of Fano type. By \cite[Theorem 1.7]{Birkar19}, there exists a positive integer $N$ depending only on $d$ such that $E_m$ is an lc place of some $N$-complement of $K_X$. Therefore we have the following by applying \cite[Corollary 3.6]{BJ20}
 $$\delta(X,\Delta)=\lim_m \frac{A_{X,\Delta}(E_m)}{S_{X, \Delta}(E_m)},$$
 where $\{E_m\}_m$ is a sequence of lc places of some $N$-complements of $K_X$.
 
 Next we turn to the case $\delta(X,\Delta)=1$. By the same proof as \cite[Theorem 1.1]{ZZ21b}, for any rational $0<\epsilon_i\ll1$, one can find an effective $\bR$-divisor $\Delta_i\sim_\bR -K_X-\Delta$ such that 
 $$\delta(X,\Delta+\epsilon_i\Delta_i)<1.$$ We assume $\{\epsilon_i\}_i$ is a decreasing sequence of numbers tending to $0$. For each $\epsilon_i$, by the proof of the previous case, there exists a sequence of prime divisors $\{E_{i,m}\}_m$ which are lc places of some $N$-complements of $K_X$ such that
 $$\delta(X,\Delta+\epsilon_i \Delta_i)=\lim_m \frac{A_{X,\Delta+\epsilon_i\Delta_i}(E_{i,m})}{S_{X, \Delta+\epsilon_i\Delta_i}(E_{i,m})},$$ 
and $N$ depends only on $d$. Thus we obtain the following 
 $$\delta(X,\Delta)=\inf_{i,m}\frac{A_{X,\Delta}(E_{i,m})}{S_{X, \Delta}(E_{i,m})}. $$
 In both cases, we could find a sequence of lc places of some $N$-complements of $K_X$ to approximate $\delta(X,\Delta)$. The proof is finished.
\end{proof}

We emphasize that the positive integer $N$ in Proposition \ref{prop: N-complement} does not depend on the boundary $\Delta$. We will also need the following result, which is a variant of \cite[Proposition 4.1]{BLX22} to remove the dependence on the boundary.

\begin{proposition}\label{prop: S-constructible}
Let $\pi: (\mX, \mB)\to T$ be a flat family of weak log Fano pairs of dimension $d$ over a smooth base $T$ with $K_\mX+\mB$ being $\bQ$-Cartier. Let $\mD$ be an effective $\bQ$-divisor on $\mX$ such that the support of $\mD$ does not contain any fiber of $\pi$, $K_\mX+\mD\sim_\bQ 0/T$, and $(\mX, \mD)\to T$ is a family of log canonical pairs. Assume $(\mX, \mB+\mD)\to T$ admits a fiberwise log resolution $g: \mY\to \mX$ and $\mF$ is a toroidal divisor with respect to $\Exc(g)+\Supp(g_*^{-1}\mB)+\Supp(g_*^{-1}\mD)$ with $A_{\mX, \mD}(\mF)<1$. Then $A_{\mX_t, \mB_t}(\mF_t)$ and $S_{\mX_t, \mB_t}(\mF_t)$ are both independent of $t\in T$.
\end{proposition}

\begin{proof}
The idea of the proof is the same as that of \cite[Proposition 4.1]{BLX22}.
The fact that $A_{\mX_t, \mB_t}(\mF_t)$ is independent of $t\in T$ follows from the existence of fiberwise log resolution. We turn to the invariance of $S_{\mX_t, \mB_t}(\mF_t)$.

By repeatedly blowing up the center of $\mF$ on $\mY$, we may assume $\mF$ is a prime divisor on $\mY$. We fix $x\in \bQ^{+}$ and aim to show the following invariant 
$$\vol(g_t^*(-K_{\mX_t}-\mB_t)-x\mF_t)$$
is independent of $t\in T$.
Let $\Gamma_1, \Gamma_2$ be two effective $\bQ$-divisors on $\mY$ without common components in their supports such that
$$K_\mY+\Gamma_1=g^*(K_\mX+\mD)+\Gamma_2 $$
and $g_*\Gamma_1=\mD$. Note that the support of $\Gamma_1+\Gamma_2$ is relative snc over $T$ and write $b:=1-A_{\mX, \mD}(\mF)>0$. By inversion of adjunction we know that $(\mX, \mD)$ is log canonical. In particular, $\Gamma_1$ has coefficients in $[0, 1]$.
Since $-K_{\mX}-\mB$ is relatively semi-ample, we may apply Bertini Theorem to find an effective $\bQ$-divisor 
$$\mH\sim_\bQ -(b/x)(K_\mX+\mB)$$ 
such that $\Gamma_1+g^*\mH-b\mF$ has coefficients in the interval $[0, 1]$ and the support of $\Gamma_1+g^*\mH-b\mF$ is relative snc over $T$. Note that we have
$$K_\mY+\Gamma_1+g^*\mH-b\mF=g^*(K_\mX+\mD+\mH)-b\mF+\Gamma_2 .$$
Applying \cite[Theorem 1.8 (3)]{HMX13} we know that 
$$\vol\left(K_{\mY_t}+(\Gamma_1)_t+g_t^*\mH_t-b\mF_t\right)=\vol\left(g_t^*(K_{\mX_t}+\mD_t+\mH_t)-b\mF_t+(\Gamma_2)_t\right) $$
is independent of $t\in T$. Also note the following formula:
\begin{align*}
(x/b)\left(g_t^*(K_{\mX_t}+\mD_t+\mH_t)-b\mF_t+(\Gamma_2)_t\right)&\sim_\bQ  (x/b)g_t^*\mH_t-x\mF_t+(x/b)(\Gamma_2)_t\\
&\sim_\bQ -g_t^*(K_{\mX_t}+\mB_t)-x\mF_t+(x/b)(\Gamma_2)_t.
\end{align*}
Hence the following invariant is independent of $t\in T$:
$$\vol(-g_t^*(K_{\mX_t}+\mB_t)-x\mF_t+(x/b)(\Gamma_2)_t). $$
We claim 
$$\vol(-g_t^*(K_{\mX_t}+\mB_t)-x\mF_t+(x/b)(\Gamma_2)_t)=\vol(-g_t^*(K_{\mX_t}+\mB_t)-x\mF_t), $$
which finishes the proof. To conclude the claim, for any sufficiently divisible positive integer $m$, we take an effective Cartier divisor
$$G\in |m(-g_t^*(K_{\mX_t}+\mB_t)-x\mF_t+(x/b)(\Gamma_2)_t)|, $$
then we have
$$G+mx\mF_t \in |-g_t^*m(K_{\mX_t}+\mB_t)+m(x/b)(\Gamma_2)_t|.$$
Since $(\Gamma_2)_t$ is $g_t$-exceptional, we see
$$|-g_t^*m(K_{\mX_t}+\mB_t)+m(x/b)(\Gamma_2)_t|=|-g_t^*m(K_{\mX_t}+\mB_t)|+m(x/b)(\Gamma_2)_t,$$
which implies
$$G+mx\mF_t\geq m(x/b)(\Gamma_2)_t.$$
Recalling that $\mF_t$ and $(\Gamma_2)_t$ have no common components, we thus have
$$G\geq m(x/b)(\Gamma_2)_t,$$
and 
$$G-m(x/b)(\Gamma_2)_t\in |m(-g_t^*(K_{\mX_t}+\mB_t)-x\mF_t)|.$$
Above all, we conclude the claim and the proof is complete.
\end{proof}

\begin{remark}
Notation as in Proposition \ref{prop: S-constructible} and its proof. Suppose $\mL$ is a $\bQ$-line bundle on $\mX$ which is big and nef over $T$, then $\mL$ is semi-ample over $T$ (note that $\mX$ is Fano type over $T$). Replacing $-K_\mX-\mB$ with $\mL$ and applying the same argument as the proof of Proposition \ref{prop: S-constructible}, we easily see that the invariant
$\vol(g_t^*\mL_t-x\mF_t) $ does not depend on $t\in T$. In particular, the following invariant is independent of $t\in T$:
$$S_{\mL_t}(\mF_t):=\frac{1}{\vol(\mL_t)}\int_0^{\infty} \vol(g_t^*\mL_t-x\mF_t){\rm d}x.$$
\end{remark}

Let $X$ be a normal projective variety of dimension $d$ and $D_j, 1\leq j\leq k,$ are Weil divisors on $X$. Recall the following definition:
$$\Kss(X, \sum_{j=1}^k D_j):=\{(x_1,...,x_k)\in \bR_{\geq 0}^k\ |\ \textit{$(X, \sum_{j=1}^k x_jD_j)$ is a K-semistable log Fano pair}\}. $$

\begin{theorem}\label{thm: finiteness}
Let $\pi: \mX\to T$ be a flat family of normal projective varieties over a smooth base $T$, and $\mD_j\to T$ is a family of effective Weil divisors for each $1\leq j\leq k$. 

Suppose for any closed point $t\in T$, there is a vector $(c_1,...,c_k)\in \bR_{\geq 0}^{k}$ such that $(\mX_t, \sum_{j=1}^kc_j\mD_{j,t})$ is a K-semistable log Fano pair. 
Then the following set is finite:
$$\{\Kss(\mX_t, \sum_{j=1}^k\mD_{j, t})\ |\ \text{$t\in T$}\}.$$
\end{theorem}

\begin{proof}
By Proposition \ref{prop: LF stratification}, up to a constructible stratification of the base $T$, we may assume that the log Fano domain $\LF(\mX_t, \sum_{j=1}^k \mD_{j,t})$ does not depend on $t\in T$, and we denote this domain by $P$.

By further stratifying the base into smooth pieces, denoted by $T=\cup_{i=1}^r T_i$, we may assume that for every $i$, the family 
$$(\mX\times_T T_i, \sum_{j=1}^k\mD_j\times_T T_i)\to T_i$$ 
admits a fiberwise log resolution up to an \'etale base change. Put $\mX_i:=\mX\times_T T_i$ and $\pi_i: \mX_i\to T_i$.

Choosing $N$ as in Proposition \ref{prop: N-complement}, up to further stratifying each $T_i$ into smooth pieces, we may assume ${\pi_i}_*\mO(-NK_{\mX_i/T_i})$ is locally free for every piece (e.g. \cite[Theorem III.12.8, Corollary III.12.9]{Har77}).

Now we can assume that $(\mX, \sum_{j=1}^k\mD_j)\to T$ admits a fiberwise log resolution and  $\pi_*\mO_{\mX}(-NK_{\mX/T})$ is locally free. Denote $W:=\bP(\pi_*\mO_\mX(-NK_{\mX/T}))\to T, $
whose fiber over $t\in T$ is exactly $|-NK_{\mX_t}|$. Let $\mH\subset \mX\times_T W$ be the universal divisor associated to $\pi_*\mO_\mX(-NK_{\mX/T})$ (this means $\mH_t\subset \mX_t\times |-NK_{\mX_t}|$ is the universal divisor associated to $-NK_{\mX_t}$), and denote $\mD:=\frac{1}{N}\mH$. Consider the morphism $(\mX, \mD)\to W$. By the lower semicontinuity of log canonical thresholds, the locus
$$Z:=\{w\in W\ |\ \textit{$(\mX_w, \mD_w)$ is strictly log canonical}\} $$
is locally closed in $W$. Then the scheme $Z$ together with the $\bQ$-divisor $\mD_Z$ (which is obtained by pulling back $\mD$ under the morphism $\mX\times_T Z\to \mX\times_T W$) parametrize boundaries of the desired form.

For the morphism 
$$\pi: (\mX\times_T Z, \mD_Z+\sum_{j=1}^k \mD_j\times_T Z)\to Z,$$ 
we choose a constructible stratification of $Z$, denoted by $Z=\cup_{i=1}^rZ_i$, such that each $Z_i$ is smooth and there is an \'etale cover $Z_i'\to Z_i$ such that 
$$\pi_i: (\mX\times_T Z, \mD_Z+\sum_{j=1}^k \mD_j\times_T Z)\times _Z Z'_i\to Z'_i$$
admits a fiberwise log resolution. Note that the fibers of 
$$\pi_i: (\mX\times_T Z, \sum_{j=1}^k x_j\mD_j\times_T Z)\times _Z Z'_i\to Z'_i$$
are all log Fano pairs for any $(x_1,...,x_k)\in P$ by the assumption in the first paragraph of the proof. Up to a further constructible stratification of each $Z_i'$, we may assume 
$$K_{\mX\times_T Z_i'}+\sum_{j=1}^k x_j\mD_j\times_T Z_i'$$
is $\bR$-Cartier for any $(x_1,...,x_k)\in P$ and any index $i$.
For any vector $\vec{x}:=(x_1,...,x_k)\in P$, we define
$$d_i(\vec{x}):=\inf\Bigg\{\frac{A_{\mX_{z_i}, \sum_{j=1}^kx_j\mD_j\times_T z_i}(v)}{S_{\mX_{z_i}, \sum_{j=1}^kx_j\mD_j\times_T z_i}(v)}\ |\ \textit{$v\in \Val^*_{\mX_{z_i}}$ and $A_{\mX_{z_i}, \mD_{z_i}}(v)=0$}\Bigg\}, $$
where $z_i$ is a closed point on $Z_i'$. By Proposition \ref{prop: S-constructible} and \cite[Proposition 4.2]{BLX22}, the number $d_i(\vec{x})$ does not depend on the choice of $z_i\in Z_i'$ when $(x_1,...,x_k)\in P(\bQ)$,  thus $d_i(\vec{x})$ is a well defined function for any $\vec{x}\in P$ by continuity. We take the notation 
$$\tilde{\delta}(X, \Delta):=\min\{\delta(X, \Delta), 1\}$$
for a log Fano pair $(X, \Delta)$. By Proposition \ref{prop: N-complement}, we have
\begin{align*}
 \tilde{\delta}(\mX_t, \sum_{j=1}^k x_j\mD_{j,t})
=\ \big\{\min_i\{1, d_i(\vec{x})\}\ |\ \textit{$t$ is contained in the image of $Z_i'$ under $Z_i'\to T$}\big\}. 
\end{align*}
Denote $p: Z\to T$. We aim to construct a suitable finite stratification of $T$. Recall that $Z=\cup_{i}^r Z_i$. For each subset $\{i_1,...,i_s\}\subset \{1,2,...,r\}$, we define the following subset of $T$:
$$T_{i_1,...,i_s}:=\Big\{b\in T\ |\ \textit{$b\in p(Z_l)$ if and only if $l\in \{i_1,...,i_s\}$}\Big\}. $$
Clearly we have a finite decomposition of $T$:
$$T=\bigcup_{\{i_1,...,i_s\}\subset\{1,2,...,k\}} T_{i_1,...,i_s}.$$
By our construction, for any $t_1, t_2\in T_{i_1,...,i_s}$, we have
$$\tilde{\delta}(\mX_{t_1}, \sum_{j=1}^k x_j\mD_{j,t_1})=\tilde{\delta}(\mX_{t_2}, \sum_{j=1}^k x_j\mD_{j,t_2})=\min\{1, d_{i_1}(\vec{x}),...,d_{i_s}(\vec{x})\} $$
for any $\vec{x}\in P$. This implies 
$$\Kss(\mX_{t_1}, \sum_{j=1}^k\mD_{j, t_1})=\Kss(\mX_{t_2}, \sum_{j=1}^k\mD_{j, t_2}).$$ 
Thus there are only finitely many K-semistable domains for fibers of the morphism $(\mX, \sum_{j=1}^k\mD_j)\to T$.
\end{proof}

\section{Wall crossing for K-semistability}\label{sec: wall crossing}

Recalling that in \cite{Zhou23b}, we define the set $\mE:=\mE(d,k,v, I)$ of log pairs  for two positive integers $d$ and $k$, a positive number $v$,  and a positive integer $I$. More precisely, a log pair $(X, \sum_{i=j}^k D_j)$ is contained in $\mE$ if and only if it satisfies the following conditions:
\begin{enumerate}
\item $X$ is a Fano variety of dimension $d$ and $(-K_X)^d=v$;
\item $D_j$ is an effective $\bQ$-divisor satisfying $D_j\sim_\bQ -K_X$ for every $1\leq j\leq k$;
\item $I(K_X+D_j)\sim 0$ for every $1\leq j\leq k$;
\item there exists a rational vector $(c_1,...,c_k)\in \bQ_{\geq 0}^k$ such that $(X, \sum_{j=1}^k c_jD_j)$ is a K-semistable log Fano pair.
\end{enumerate}
It turns out that the set $\mE$ is log bounded and the set of K-semistable domains 
$$\{\Kss(X, \sum_{j=1}^kD_j)\ |\ \textit{$(X, \sum_{j=1}^kD_j)\in \mE$}\}$$ is finite (e.g. \cite[Theorem 1.4]{Zhou23b}). In the definition of $\mE$, the conditions (1) and (2) bring much convenience, and we call this setting to be proportional since all the boundary divisors are proportional to the anti-canonical divisor. In this section, we aim to prove the corresponding result in non-proportional setting.

We first define a set $\mG:=\mG(d,k,P)$ of couples for the given data $d, k, P$, where $d,k$ are two positive integers and  $P\subset [0,1]^k$ is a rational polytope. 
We say a couple $(X, \sum_{j=1}^kD_j)$ belongs to $\mG$ if and only if the following conditions are satisfied:
\begin{enumerate}
\item $X$ is a normal projective variety of dimension $d$;
\item $D_j, 1\leq j\leq k,$ are effective Weil divisors on $X$;
\item there exists a vector $(c_1,...,c_k)\in P$ such that $(X, \sum_{j=1}^kc_jD_j)$ is a K-semistable log Fano pair.
\end{enumerate}
We define the K-semistable domain (restricted to $P$) for each couple in $\mG$ as follows:
$$\Kss(X, \sum_{j=1}^k D_j)_P:=\{(x_1,...,x_k)\in P |\ \textit{$(X, \sum_{j=1}^k x_jD_j)$ is a K-semistable log Fano pair}\}.$$

\begin{lemma}\label{lem: locally closed}
Given $(X, \sum_{j=1}^kD_j)\in \mG:=\mG(d,k,P)$. Then $\Kss(X, \sum_{j=1}^k D_j)_P$ is a closed subset in $P\cap \LF(X, \sum_{j=1}^kD_j)$.
\end{lemma}

\begin{proof} 
Let $\{\vec{c}_i:=(c_{i1},...,c_{ik})\}_{i=1}^\infty$ be a sequence of vectors tending to $\vec{c}:=(c_1,...,c_k)$. Suppose $\{\vec{c}_i\}_{i=1}^\infty\cup \vec{c}\subset P\cap \LF(X, \sum_{j=1}^kD_j)$ and $\{\vec{c}_i\}_{i=1}^\infty\subset \Kss(X, \sum_{j=1}^k D_j)_P$. It suffices to show that $\vec{c}\in \Kss(X, \sum_{j=1}^k D_j)_P$. 

For any prime divisor $E$ over $X$, we have
$$A_{X, \sum_{j=1}^k c_{ij}D_j}(E)-S_{X, \sum_{j=1}^k c_{ij}D_j}(E) \geq 0$$
for any index $i$. By the continuity we have
$$A_{X, \sum_{j=1}^k c_{j}D_j}(E)-S_{X, \sum_{j=1}^k c_{j}D_j}(E) \geq 0.$$
Hence $(X, \sum_{j=1}^kc_jD_j)$ is a K-semistable log Fano pair by Definition \ref{def: Fujita-Li}.
\end{proof}

Since $\mG$ is not necessarily log bounded, we will consider a log bounded subset in the following result.

\begin{theorem}\label{thm: finiteness P}
Suppose $\widetilde{\mG}\subset \mG$ is a subset which is log bounded, then the following set is finite
$$\{\Kss(X, \sum_{j=1}^kD_j)_P | \textit{$(X, \sum_{j=1}^kD_j)\in \widetilde{\mG}$}\}.$$
\end{theorem}

\begin{proof}
Since $\widetilde{\mG}$ is log bounded, there exist a projective morphsim $\pi: \mX\to T$ of finite type schemes and Weil divisors $\mD_j $'s ($1\leq j\leq k$) on $\mX$ such that for any $(X, \sum_{j=1}^k D_j)\in \widetilde{\mG}$ we could find a closed point $t\in T$ such that 
$$(\mX_t, \sum_{j=1}^k\mD_{j,t})\cong (X, \sum_{j=1}^k D_j).$$ 
Up to a constructible stratification of $T$ into smooth pieces, we apply Proposition \ref{prop: algebraic locus} to get finitely many families $\pi_i: (\mX^{(i)}, \sum_j\mD^{(i)}_j)\to T_i$  satisfying the following conditions:
\begin{enumerate}
\item for each $i$, $\pi_i: \mX^{(i)}\to T_i$ is a flat family of normal projective varieties over the smooth base $T_i$; 
\item  for each $i$ and $j$, $\mD^{(i)}_j\to T_i$ is a family of effective Weil divisors;
\item for each $i$ and each $t\in T_i$, $(\mX^{(i)}_{t}, \sum_{j=1}^k c_j\mD^{(i)}_{j,t})$ is a K-semistable log Fano pair for some vector $(c_1,...,c_k)$;
\item for any $(X, \sum_{j=1}^kD_j)\in \widetilde{\mG}$, it is isomorphic to some fiber of the morphism $\pi_i: (\mX^{(i)}, \sum_{j=1}^k\mD^{(i)}_j)\to T_i$ for some $i$.
\end{enumerate}
Applying Theorem \ref{thm: finiteness}, we get the desired finiteness of K-semistable domains.
\end{proof}

As a direct corollary, we have the following result on finite chamber decomposition.

\begin{corollary}\label{cor: chamber}
Suppose $\widetilde{\mG}\subset\mG$ is a subset which is log bounded. Then there exists a finite disjoint union $P=\amalg_{i=0}^m\ I_i$ such that each $I_i$ is a locally closed subset of $P$ under Euclidean topology, and for any $(X, \sum_{j=1}^k D_j)\in \widetilde{\mG}$, the K-semistability of $(X, \sum_{j=1}^kc_jD_j)$ does not change as $(c_1,...,c_k)$ varies in $I_i$.
\end{corollary}

\begin{proof}
By Theorem \ref{thm: finiteness P} and Lemma \ref{lem: locally closed}, there are finitely many locally closed subsets of $P$ (under Euclidean topology), denoted by $K_1,...,K_l$, such that for any $(X, \sum_{j=1}^kD_j)\in \widetilde{\mG}$, we have $\Kss(X, \sum_{j=1}^k D_j)_P=K_r$ for some $1\leq r\leq l$. One could easily cook up a finite disjoint union of $\cup_{r=1}^l K_r$, denoted by $\cup_{r=1}^l K_r=\amalg_{i=1}^m\ I_i$, such that each $I_i$ is locally closed, and for any $(X, \sum_{j=1}^k D_j)\in \widetilde{\mG}$, the K-semistability of $(X, \sum_{j=1}^kc_jD_j)$ does not change as $(c_1,...,c_k)$ varies in $I_i$. Denote $I_0:=P\setminus (\cup_{r=1}^lK_r)$. The disjoint union $P=\amalg_{i=0}^m\ I_i$ is the required one.
\end{proof}

A priori, the chamber decomposition in Corollary \ref{cor: chamber} may be quite pathological. At least, it is not necessarily a polytope chamber decomposition as in \cite{Zhou23b}. We will treat this problem in Section \ref{sec: fcd}.

We give two examples of the set $\widetilde{\mG}$ to make sure log boundedness.

\begin{example}\label{example: P_1}
Fix two positive integers $d,k$, fix a polytope $P:=[\epsilon, 1]^k$ for some positive rational number $\epsilon$, and fix a positive real number $\epsilon_0$.
We define a set $\widetilde{\mG}\subset\mG(d,k,P)$ of couples for the given data $d, k, P, \epsilon_0$.
We say a couple $(X, \sum_{j=1}^kD_j)$ belongs to $\widetilde{\mG}$ if and only if the following conditions are satisfied:
\begin{enumerate}
\item $X$ is a normal projective variety of dimension $d$;
\item $D_j, 1\leq j\leq k,$ are prime Weil divisors on $X$;
\item there exists a vector $(c_1,...,c_k)\in P$ such that $(X, \sum_{j=1}^kc_jD_j)$ is a K-semistable log Fano pair with $\vol(-K_X-\sum_{j=1}^kc_jD_j)\geq \epsilon_0$.
\end{enumerate}
We claim that the set $\widetilde{\mG}$ is log bounded. To see this, 
we first show the set 
$$\{X\ |\  \text{$(X, \sum_{j=1}^kD_j)\in \widetilde{\mG}$}\}$$ 
is bounded. For any $(X, \sum_{j=1}^kD_j)\in \widetilde{\mG}$, there exists a vector $(c_1,...,c_k)\in P$ such that $(X, \sum_{j=1}^kc_jD_j)$ is a K-semistable log Fano pair with $\vol(-K_X-\sum_jc_jD_j)\geq \epsilon_0$. By Theorem \ref{cor: bdd}, we get the desired boundedness. Hence there exists a very ample line bundle $A$ on $X$ such that 
$$A^d\leq N_1\quad \text{and} \quad -K_X\cdot A^{d-1}\leq N_2,$$ 
where $N_1, N_2$ are two positive numbers depending only on $d, k, P, \epsilon_0$. Note that $-K_X-\sum_{j=1}^k c_jD_j$ is ample, we have
$$(\sum_{j=1}^k c_j D_j)\cdot A^{d-1}\leq -K_X\cdot A^{d-1}\leq N_2 .$$
Since $c_j\geq \epsilon$ for every $j$, we see that $D_j\cdot A^{d-1}\leq \frac{N_2}{\epsilon}$, which implies that $D_j$ is also bounded for every $j$. So $\widetilde{\mG}$ is log bounded as desired (e.g. \cite[Lemma 2.20]{Birkar19}). Given this log boundedness, we could apply Corollary \ref{cor: chamber} to obtain a chamber decomposition of $P$ (into a disjoint union of finitely many locally closed subsets under Euclidean topology) to control the variation of K-semistability of couples in $\widetilde{\mG}$.
\end{example}

\begin{example}\label{example: P_2}
This example is related to log Fano pairs of Maeda type. Fix two positive integers $d, k$, and denote $P:=[0,1]^k$. We define a set $\widetilde{\mG}\subset\mG(d,k,P)$ of couples depending on $d,k, P$. A couple $(X, \sum_{j=1}^kD_j)$ belongs to $\widetilde{\mG}$ if and only if the following conditions are satisfied:
\begin{enumerate}
\item $X$ is a normal projective variety of dimension $d$;
\item $D_j, 1\leq j\leq k,$ are prime Weil divisors on $X$;
\item $(X, \sum_{j=1}^kD_j)$ is a degeneration of log Fano manifolds;
\item $(X, \sum_{j=1}^kc_jD_j)$ is a K-semistable log Fano pair for some $(c_1,...,c_k)\in P.$
\end{enumerate}
Here, we say $(X, \sum_{j=1}^kD_j)$ is a log Fano manifold if it is log smooth with $-K_X-\sum_{j=1}^kD_j$ being ample; we say $(X, \sum_{j=1}^kD_j)$ is a degeneration of log Fano manifolds if there is a family of couples $(\mX, \sum_{j=1}^{k}\mD_j)\to  C\ni 0$ over a smooth pointed curve such that $(\mX_t, \sum_{j=1}^k\mD_{j,t})$ is a log Fano manifold for any $t\ne 0$ and $(\mX_0, \sum_{j=1}^k\mD_{j,0})\cong (X, \sum_{j=1}^kD_k)$. 

In this setting, it is not hard to see that $\vol(-K_X-\sum_{j=1}^kc_jD_j)$ is always bounded from below for any $(X, \sum_{j=1}^kD_j)\in \widetilde{\mG}$ and any $(c_1,...,c_k)\in P$. Similar as Example \ref{example: P_1}, we see
$\widetilde{\mG}$ is log bounded (the subset of $\widetilde{\mG}$ consisting of log Fano manifolds is proven to be log bounded in \cite{LZ23} by using another approach). 
\end{example}

\section{Boundedness of couples}\label{sec: bdd set}

Given two positive integers $d,k$ and a rational polytope $P\subset [0,1]^k$. We first define a set of couples $\mA:=\mA(d,k, P)$ depending on $d,k$ and $P$.  We say a couple $(X, \sum_{j=1}^kD_j)$ belongs to $\mA$ if and only if the following conditions are satisfied:
\begin{enumerate}
\item $X$ is a normal projective variety of dimension $d$;
\item $D_j$'s ($1\leq j\leq k$) are effective Weil divisors on $X$;
\item there exists a vector $(c_1,...,c_k)\in P$ such that $(X, \sum_{j=1}^kc_jD_j)$ is a weak log Fano pair.
\end{enumerate}

Fixing another two positive real numbers $\delta_0$ and $\epsilon$, we define two subsets of $\mA$ associated to $\delta_0$ and $\epsilon$, denoted by $\mA_i:=\mA_i(\delta_0, \epsilon)$ for $i=1,2$. We say $(X, \sum_{j=1}^kD_j)\in \mA$ belongs to $\mA_1$ (resp. $\mA_2$) if and only if there exists some vector $(c_1,...,c_k)\in P$ such that $(X, \sum_{j=1}^kc_jD_j)$ is a  log Fano pair (resp. weak log Fano pair) with 
$\delta(X, \sum_{j=1}^kc_jD_j)\geq \delta_0$, and $\vol(-K_X-\sum_{j=1}^kx_jD_j)\geq \epsilon$ for any $(x_1,...,x_k)\in P$.\footnote{Here the volume of an $\bR$-divisor of $X$ is well-defined since $X$ is of Fano type and hence we can take a small $\bQ$-factorization of $X$ and define volume there. It is easy to see that it does not depend on the choice of small $\bQ$-factorizations.}
It is clear that $\mA_1$ is contained in $\mA_2$. We have the following boundedness result.

\begin{proposition}\label{prop: Abdd}
Both $\mA_1$ and $\mA_2$ are log bounded.
\end{proposition}

\begin{proof}
Up to replacing $k$ with a smaller integer, we may assume there exists a point of $P$, denoted by $(a_1,...,a_k)$, such that $a_j\ne 0$ for every $j$.

It suffices to show that $\mA_2$ is log bounded. By Proposition \ref{prop: polytope LF} and \cite[Proposition 2.10]{LZ24}, for any couple $(X, \sum_{j=1}^kD_j)\in \mA_2$, there exists a vector $(c_1,...,c_k)\in P(\bQ)$ such that $(X, \sum_{j=1}^kc_jD_j)$ is a  weak log Fano pair with 
$$\delta(X, \sum_{j=1}^kc_jD_j)\geq \delta_0/2\quad \text{and}\quad \vol(-K_X-\sum_{j=1}^kc_jD_j)\geq \epsilon.$$ 
By Theorem \ref{cor: bdd} and the sentence below Theorem \ref{cor: bdd}, we know the following set of varieties is bounded:
$$\{X\ |\  \text{$(X, \sum_{j=1}^kD_j)\in \mA_2$}\}.$$
Thus there exists a very ample line bundle $A$ on $X$ with 
$$A^{d}\leq N_1\quad\text{and}\quad A^{d-1}\cdot (-K_X)\leq N_2,$$ 
where $N_1, N_2$ are two positive numbers which do not depend on the choice of $X$ in the above set.
On the other hand, since $-K_X-\sum_{j=1}^ka_jD_j$ is pseudo-effective, we have
$$A^{d-1}\cdot D_j\leq \frac{A^{d-1}\cdot (-K_X)}{a_j}\leq \frac{N_2}{a_j},$$
which implies that $D_j$ is also bounded for each $j$ and hence $\mA_2$ is log bounded.  
\end{proof}

When we consider the couples in $\mA_1$ or $\mA_2$, say $(X, \sum_{j=1}^d D_j)$, we cannot freely vary the coefficients of $(X, \sum_{j=1}^kc_jD_j)$ in $P$, since it is not necessarily log Fano or weak log Fano for some vectors in $P$. The following result will ensure us to vary the coefficients freely in each chamber under a finite rational polytope chamber decomposition of $P$.

\begin{proposition}\label{prop: perfect fcd}
There exists a finite chamber decomposition of $P$ (depending only on $\mA_1$ and $\mA_2$), denoted by $P=\cup_s P_s$, 
satisfying the following condition: every $P_s$ is a rational polytope and $P_s^\circ\cap P_{s'}^\circ=\emptyset$ if $P_s$ and $P_{s'}$ are adjacent;
moreover, for any chamber $P_s$ and any face $F$ of $P_s$, we have the following:
\begin{enumerate}
\item for any $(X,\sum_{j=1}^kD_j)\in \mA_1$, if $(X, \sum_{j=1}^kc_jD_j)$ is a log Fano pair for some vector $(c_1,...,c_k)\in F^\circ$, then $(X,\sum_{j=1}^kx_jD_j)$ is a log Fano pair for any vector $(x_1,...,x_k)\in F^\circ$;
\item for any $(X,\sum_{j=1}^kD_j)\in \mA_2$, if $(X, \sum_{j=1}^kc_jD_j)$ is a weak log Fano pair for some vector $(c_1,...,c_k)\in F^\circ$, then $(X,\sum_{j=1}^kx_jD_j)$ is a weak log Fano pair for any vector $(x_1,...,x_k)\in F^\circ$.
\end{enumerate}
\end{proposition}
\begin{proof}
Since $\mA_1$ and $\mA_2$ are both log bounded, there exist suitable families parametrizing them by the same argument as in the proof of Theorem \ref{thm: finiteness P}. Applying Proposition \ref{prop: LF stratification}, we could find two required finite rational polytope chamber decompositions of $P$ respectively for $\mA_1$ and $\mA_2$. A common refinement of the two chamber decompositions is the desired one.
\end{proof}

We will need the following result in the subsequent section.

\begin{proposition}\label{prop: uniform polarization}
Given $P, \mA_1, \mA_2$ as above and let $P=\cup_s P_s$ be the chamber decomposition satisfying the properties in Proposition \ref{prop: perfect fcd}. Let $\pi: \mX\to T$ be a flat family of normal projective varieties over a smooth base $T$, and $\mD_j\to T$ is a family of effective Weil divisors for each $1\leq j\leq k$. 

Suppose every fiber of $(\mX, \sum_{j=1}^k\mD_j)\to T$ is contained in $\mA_1$ and assume $P\subset \overline{\LF(\mX_t, \sum_{j=1}^k\mD_{j,t})}$ for any $t\in T$. Then we have the following conclusion: for any index $s$ and any face $F$ of $P_s$, if there exists a vector $(c_1,...,c_k)\in F^\circ$ such that $-K_{\mX/T}-\sum_{j=1}^kc_j\mD_j$ is $\bR$-Cartier and $(\mX_t, \sum_{j=1}^kc_j\mD_{j,t})$ is log Fano with $\delta(\mX_t, \sum_{j=1}^kc_j\mD_{j,t})\geq \delta_0$ for any $t\in T$, then $-K_{\mX/T}-\sum_{j=1}^kx_j\mD_j$ is $\bR$-Cartier and $(\mX_t, \sum_{j=1}^kx_j\mD_{j,t})$ is log Fano for any $(x_1,...,x_k)\in F^\circ$ and any $t\in T$.
\end{proposition}

\begin{proof}
By Proposition \ref{prop: perfect fcd}, $(\mX_t, \sum_{j=1}^kx_j\mD_{j,t})$ is log Fano for any $(x_1,...,x_k)\in F^\circ$ and any $t\in T$. It remains to show that $-K_{\mX/T}-\sum_{j=1}^kx_j\mD_j$ is $\bR$-Cartier for any $(x_1,...,x_k)\in F^\circ$ and it is enough to test rational $(x_1,...,x_k)\in F^\circ(\bQ)$.

First observe that $(\mX, \sum_{j=1}^kc_j\mD_j)$ is a klt (quasi-projective) log pair. Thus there exists a small $\bQ$-factorization $\phi$ as below by \cite{BCHM10}:
\begin{center}
	\begin{tikzcd}[column sep = 2em, row sep = 2em]
	 (\mX', \sum_{j=1}^kc_j\mD_j') \arrow[rd,"",swap] \arrow[rr,"\phi"]&& (\mX, \sum_{j=1}^kc_j\mD_j) \arrow[ld,"\pi"]\\
	 &T&.
	\end{tikzcd}
\end{center}
We first show that for each closed point $t\in T$, the fiber $(\mX'_t, \sum_{j=1}^k\mD'_{j,t})$ belongs to $\mA_2$.
To see this, we first note the following formula:
$$K_{\mX'/T}+\sum_{j=1}^kc_j\mD_j' =\phi^*(K_{\mX/T}+\sum_{j=1}^kc_j\mD_j).$$
Write $n=\dim T$. Taking $n$ sufficiently general hypersurfaces passing through $t\in T$, denoted by $H_1,...,H_n$,  we have
$$K_{\mX'/T}+\sum_{j=1}^kc_j\mD_j' +\sum_{r=1}^n\phi^*\pi^*H_r=\phi^*(K_{\mX/T}+\sum_{j=1}^kc_j\mD_j+\sum_{r=1}^n \pi^*H_r) .$$ 
Note that both sides are dlt by the generality of $H_r$'s. Applying adjunction we have the following formula $(\star)$:
$$K_{\mX'_t}+\sum_{j=1}^rc_j\mD'_{j,t}=\phi^*(K_{\mX_t}+\sum_{j=1}^kc_j\mD_{j,t}). $$
It is clear that $(\mX'_t, \sum_{j=1}^kc_j\mD'_{j,t})$ is a weak log Fano pair with $\delta(\mX'_t,  \sum_{j=1}^kc_j\mD'_{j,t})\geq \delta_0.$
On the other hand, $(\mX'_t, \sum_{j=1}^k\mD'_{j,t})$ is a small modification of $(\mX_t, \sum_{j=1}^k\mD_{j,t})$ for general $t\in T$. Applying \cite[Theorem 12.8, Corollary 12.9]{Har77} we see that 
$$\vol(-K_{\mX'_t}-\sum_{j=1}^ka_j\mD'_{j,t})\geq \epsilon$$ 
for any $t\in T$ and any $(a_1,...,a_k)\in P$.
This implies that $(\mX'_t, \sum_{j=1}^k\mD'_{j,t})$ indeed belongs to $\mA_2$ for any $t\in T$.
By Proposition \ref{prop: perfect fcd}, we know that $(\mX'_t,\sum_{j=1}^ka_j\mD'_{j,t})$ is a weak log Fano pair for any $t\in T$ and any vector $(a_1,...,a_k)\in F^\circ$.

Next we show that for any $(x_1,...,x_k)\in F^\circ$, the following crepant relation $(\star\star)$ holds:
$$K_{\mX'_t}+\sum_{j=1}^kx_j\mD'_{j,t}=\phi^*(K_{\mX_t}+\sum_{j=1}^kx_j\mD_{j,t}). $$
To see this, we denote $\vec{c}:=(c_1,...,c_k)$ and $\vec{x}:=(x_1,...,x_k)$. Then we choose another vector $\vec{y}:=(y_1,...,y_k)\in F^\circ$ lying on the line passing through $\vec{c}$ and $\vec{x}$ such that $\vec{c}$ stays in the middle of $\vec{x}$ and $\vec{y}$.  Note that both $(\mX'_t, \sum_{j=1}^kx_j\mD'_{j,t})$ and $(\mX'_t, \sum_{j=1}^ky_j\mD'_{j,t})$ are weak log Fano pairs. Applying non-negativity lemma (\cite[Lemma 3.39]{KM98}), we have
$$K_{\mX'_t}+\sum_{j=1}^kx_j\mD'_{j,t}-E_1=\phi^*(K_{\mX_t}+\sum_{j=1}^kx_j\mD_{j,t}), $$
and 
$$K_{\mX'_t}+\sum_{j=1}^ky_j\mD'_{j,t}-E_2=\phi^*(K_{\mX_t}+\sum_{j=1}^ky_j\mD_{j,t}). $$
where $E_1, E_2$ are both effective exceptional divisors. By a linear combination of the above two equations and comparing with the equation $(\star)$, we easily derive $E_1=E_2=0$. Thus we get the desired crepant relation $(\star\star)$.

Now for any given rational vector $(x_1,...,x_k)\in F^\circ(\bQ)$, we know that $(\mX'_t, \sum_{j=1}^kx_j\mD'_{j,t})$ is a crepant model of $(\mX_t, \sum_{j=1}^kx_j\mD_{j,t})$. Considering the following morphism by changing the coefficients:
\begin{center}
	\begin{tikzcd}[column sep = 2em, row sep = 2em]
	 (\mX', \sum_{j=1}^kx_j\mD_j') \arrow[rd,"",swap] \arrow[rr,"\phi"]&& (\mX, \sum_{j=1}^kx_j\mD_j) \arrow[ld,"\pi"]\\
	 &T&,
	\end{tikzcd}
\end{center}
we thus see that $\phi$ induces fiberwise anti-canonical models with coefficients $x_j$'s.
We aim to show that $(\mX, \sum_{j=1}^kx_j\mD_j)$ is the anti-canonical model of $(\mX', \sum_{j=1}^kx_j\mD'_j)$, which implies that $-K_{\mX/T}-\sum_{j=1}^kx_j\mD_j$ is indeed $\bQ$-Cartier. Suppose $(\mX'', \sum_{j=1}^kx_j\mD''_j)$ is the anti-canonical model of $(\mX', \sum_{j=1}^kx_j\mD'_j)$ as follows
\begin{center}
	\begin{tikzcd}[column sep = 2em, row sep = 2em]
	 (\mX', \sum_{j=1}^kx_j\mD_j') \arrow[rd,"\pi'",swap] \arrow[rr,"\phi'"]&& (\mX'', \sum_{j=1}^kx_j\mD_j'') \arrow[ld,""]\\
	 &T&,
	\end{tikzcd}
\end{center}
it suffices to show that $\phi'$ induces fiberwise anti-canonical models. To see this, we take any sufficiently divisible $m\in \bZ^+$ such that $-m(K_{\mX'/T}+\sum_{j=1}^kx_j\mD_j')$ is free over $T$ and it suffices to show the surjection of the following map $(\spadesuit)$:
$$\pi'_*\left(-m(K_{\mX'/T}+\sum_{j=1}^kx_j\mD_j')\right)\to H^0(\mX'_t, -m(K_{\mX_t'}+\sum_{j=1}^kx_j\mD_{j,t}')). $$
Let $H_1,...,H_n$ be sufficiently general hypersurfaces passing through $t$ as before and denote $H_0:=T$. We have the following surjection for each $0\leq r\leq n-1$ by  Kawamata-Viehweg vanishing:
\begin{align*}
&\pi'_*\left(-m(K_{\cap_{i=0}^r\pi'^*H_i/\cap_{i=0}^rH_i}+\sum_{j=1}^kx_j{\mD'_j}|_{\cap_{i=0}^r\pi'^*H_i})\right)\\
\to\ & \pi'_*\left(-m(K_{\cap_{i=0}^{r+1}\pi'^*H_i/\cap_{i=0}^{r+1}H_i}+\sum_{j=1}^kx_j{\mD'_j}|_{\cap_{i=0}^{r+1}\pi'^*H_i}) \right) .
\end{align*}
Composing these surjective maps we get the desired surjection of the map $(\spadesuit)$.
\end{proof}

As a direct corollary of Proposition \ref{prop: uniform polarization}, we have the following result.

\begin{corollary}\label{prop: perfect tc} 
Given $P, \mA_1, \mA_2$ as above and let $P=\cup_s P_s$ be the chamber decomposition satisfying the properties in Proposition \ref{prop: perfect fcd}. Take a couple $(X, \sum_{j=1}^kD_j)\in \mA_1$ and let $F$ be a face of a chamber $P_s$. Suppose $(X, \sum_{j=1}^kc_jD_j)$ is a log Fano pair with $\delta(X, \sum_{j=1}^kc_jD_j)\geq \delta_0$
for some vector $(c_1,...,c_k)\in F^\circ$, and $(\mX, \sum_{j=1}^kc_j\mD_j)\to \bA^1$ is a special test configuration of $(X, \sum_{j=1}^kc_jD_j)$ such that $\delta(\mX_0, \sum_{j=1}^kc_j\mD_{j,0})\geq \delta_0$. Then $(\mX, \sum_{j=1}^kx_j\mD_j)\to \bA^1$ is a special test configuration of $(X, \sum_{j=1}^kx_jD_j)$ for any vector $(x_1,...,x_k)\in F^\circ$.
\end{corollary}

\begin{proof}
First note that $(\mX_0, \sum_{j=1}^k\mD_{j,0})\in \mA_1$. Applying Proposition \ref{prop: perfect fcd}, we know that $(\mX_t, \sum_{j=1}^kx_j\mD_{j,t})$ is a log Fano pair for any $t\in \bA^1$ and any vector $(x_1,...,x_k)\in F^\circ$.
Thus to show that $(\mX, \sum_{j=1}^kx_j\mD_j)\to \bA^1$ is a special test configuration, it suffices to show that $-K_\mX-\sum_{j=1}^kx_j\mD_j$ is $\bQ$-Cartier for any $(x_1,...,x_k)\in F^\circ(\bQ)$.
This is implied by Proposition \ref{prop: uniform polarization}.
\end{proof}

\section{On the shape of the K-semistable domain}\label{sec: shape}

Unlike the proportional setting as in \cite{Zhou23b, Zhou23a}, the shape of the K-semistable domain may not be well behaved in non-proportional case. For instance, some non-polytope examples are worked out in \cite{Loginov23}. In this section, we will prove a property of K-semistable domain to make sure that it is not so much pathological.

\subsection{CM-line bundle}

In this subsection, we work in the following special setting $(\heartsuit)$:  
$\pi: (\mX, \sum_{j=1}^k\mD_j)\to T$ is a family of couples of relative dimension $d$ over a smooth base $T$, where $\mX$ is flat over $T$ and each $\mD_j$ is an effective Weil divisor which is either flat over $T$ or induces a well-defined family of cycles  (e.g. \cite[Section I.3]{Kollar96}); $P\subset [0,1]^k$ is a rational polytope satisfying that $\pi: (\mX, \mD(\vx))\to T$ is a family of log Fano pairs and $\mL(\vec{x}):=-K_{\mX/T}-\mD(\vx)$ is $\bR$-Cartier for any $\vec{x}=(x_1,...,x_k)\in P^\circ$, where we write $\mD(\vec{x}):=\sum_{j=1}^kx_j\mD_j$.

In the setting $(\heartsuit)$, $\mL(\vec{x})$ is indeed a polarization when $\vec{x}$  lies in $P^\circ$, but  not necessarily so when $\vec{x}$ lies on the boundary of $P$ (since it is only relatively nef in this case).

By the work of Mumford-Knudsen \cite{KM76}, for any $(x_1,...,x_k)\in P^\circ(\bQ)$, there exist (uniquely determined) $\bQ$-line bundles $\lambda_i(\vx) \ (i=0,1,...,d+1)$ on $T$ such that we have the following expansion for all sufficiently divisible $m\in \bN$:
$$\det \pi_*(\mL(\vx)^{\otimes m})= \lambda_{d+1}(\vx)^{\binom{m}{d+1}}\otimes\lambda_d(\vx)^{\binom{m}{d}}\otimes...\otimes\lambda_1(\vx)^{\binom{m}{1}}\otimes\lambda_0(\vx).$$
Following \cite[Definition 9.29]{Xu24}, we give below the definition on the CM bundle associated to $\pi: (\mX, \mD(\vx); \mL(\vx))\to T$.

\begin{definition}\label{def: CM}
For any given rational vector $\vec{c}:=(c_1,...,c_k)\in P^\circ(\bQ)$, the CM-line bundle for the family $\pi: (\mX,\mD(\vc);\mL(\vc))\to T$ is a $\bQ$-line bundle on $T$ defined as follows:
$$\lambda_{\CM,(\mX,\mD(\vc),\mL(\vc);\pi)}:=-\lambda_{d+1}(\vc). $$
\end{definition}

If the base $T$ is proper, by Riemann-Roch formula, e.g. \cite[Appendix]{CP21}, we have
$$\lambda_{\CM,(\mX,\mD(\vc),\mL(\vc);\pi)}= -\pi_*(\mL(\vc)^{d+1})=-\pi_*(-K_{\mX/T}-\mD(\vc))^{d+1}.$$

We have the following well-known property of the CM line bundle.

\begin{proposition}\label{prop: PT}
Notation as above. Take $\vc:=(c_1,...,c_k)\in P^\circ(\bQ)$. Suppose $\pi: (\mX,\mD(\vc);\mL(\vc))\to T$ is a special test configuration ($T=\bA^1$ in this case) of a log Fano pair $(X, \sum_{j=1}^kc_jD_j)$, where $D_j$'s ($1\leq j\leq k$) are effective Weil divisors on $X$. Then $\lambda_{\CM,(\mX,\mD(\vc),\mL(\vc);\pi)}$ is a $\bC^*$-linearized $\bQ$-line bundle on the base $\bA^1$ and we have
$$\Fut(\mX,\mD(\vc);\mL(\vc))=\frac{1}{(d+1)L(\vc)^d}\cdot \wt(\lambda_{\CM,(\mX,\mD(\vc),\mL(\vc);\pi)}),$$
where $L(\vc):=-K_X-\sum_{j=1}^k c_jD_j$.
\end{proposition}

\begin{proof}
Standard by \cite{PT06}.
\end{proof}

Note that in Definition \ref{def: CM}, the polarization is a $\bQ$-line bundle and the boundary has $\bQ$-coefficients. However, since we will also encounter ample $\bR$-line bundles as polarizations and the boundaries are also allowed to have $\bR$-coefficients, we aim to give the definition of the CM-line bundle associated to an $\bR$-polarization in the setting $(\heartsuit)$.

\begin{proposition}\label{prop: R-CM}
In the setting $(\heartsuit)$, where $T$ is not necessarily proper, there exist finitely many $\bQ$-line bundles $M_q$'s on $T$ and finitely many $\bQ$-polynomials $g_q(x_1,...,x_k)$'s such that
$$\lambda_{\CM,(\mX,\mD(\vx),\mL(\vec{x});\pi)}= \sum_q g_q(x_1,...,x_k) \cdot M_q $$
for any $\vx=(x_1,...,x_k)\in P^\circ(\bQ)$. Moreover, if $\pi: (\mX, \sum_{j=1}^k\mD_j)\to T$ admits a $G$-equivariant action for some group $G$, then $M_q$'s can be chosen $G$-equivariant.
\end{proposition}

\begin{proof}
For a rational vector $\vx=(x_1,...,x_k)\in P^\circ(\bQ)$,  by Definition \ref{def: CM}, we have 
$$\lambda_{\CM,(\mX,\mD(\vx),\mL(\vx);\pi)}:=-\lambda_{d+1}(\vx). $$
By \cite[Theorem 1]{PRS08}, we know that $\lambda_{d+1}(\vx)$ can be reformulated as the following Deligne pairing:
$$\lambda_{d+1}(\vx)= \langle\mL(\vx), ..., \mL(\vx)\rangle_{\mX/T}, $$
where there are $d+1$ terms (all of which are equal to $\mL(\vx)$) in the bracket. 
Let $\vec{a}_p:=(a_{p1},...,a_{pk})$'s be the finite vertices of $P$ (indexed by $p$), then we easily see:
$$ \mL(\vx)=-K_{\mX/T}-\mD(\vx)=\sum_{p} f_p(x_1,...,x_k)\cdot (-K_{\mX/T}-\mD(\va_p))=\sum_p f_p(x_1,...,x_k)\cdot \mL(\va_p),$$
where $f_p(x_1,...,x_k)$'s are finite $\bQ$-linear functions. Note that $\mL(\va_p)$ is a $\bQ$-line bundle for each index $p$ because $\vec{a}_p$ is a rational vector (recall that $P$ is a rational polytope). Thus we see
$$\lambda_{\CM,(\mX,\mD(\vx),\mL(\vx);\pi)}=- \langle\sum_p f_p(x_1,...,x_k)\cdot\mL(\va_p), ..., \sum_p f_p(x_1,...,x_k)\cdot\mL(\va_p)\rangle_{\mX/T}.$$
Recalling that the Deligne pairing is multilinear and symmetric (e.g. \cite[Section 1.1]{Zhang96}), the above formula could then be reformulated as
$$\lambda_{\CM,(\mX,\mD(\vx),\mL(\vec{x});\pi)}= \sum_q g_q(x_1,...,x_k) \cdot M_q $$
for some $\bQ$-polynomials $g_q(x_1,...,x_k)$'s and $\bQ$-line bundles $M_q$'s.

If $\pi: (\mX, \sum_{j=1}^k\mD_j)\to T$ admits a $G$-equivariant action for some group $G$, then all $\mL(\va_p)$'s are  $G$-equivariant $\bQ$-line bundles. Thus $M_q$'s are $G$-equivariant.
\end{proof}

\begin{definition}\label{def: R-CM}
In the setting $(\heartsuit)$, we define the CM line bundle for the family $\pi: (\mX, \mD(\vx); \mL(\vec{x}))\to T$ for any $(x_1,...,x_k)\in P$ as follows:
$$\lambda_{\CM,(\mX,\mD(\vx),\mL(\vec{x});\pi)}= \sum_q g_q(x_1,...,x_k) \cdot M_q, $$
where $M_q$'s and $g_q$'s are given by Proposition \ref{prop: R-CM}.
\end{definition}

\begin{example}\label{exam: R-CM}
In the setting $(\heartsuit)$, if $T$ is proper,  the CM line bundle for the family  $\pi: (\mX, \mD(\vx); \mL(\vec{x}))\to T$ for any $(x_1,...,x_k)\in P$ could be formulated as $(\star)$: %
$$\lambda_{\CM,(\mX,\mD(\vx),\mL(\vec{x});\pi)}:= -\pi_*(\mL(\vec{x})^{d+1})=-\pi_*(-K_{\mX/T}-\mD(\vx))^{d+1}.$$
In fact, when $\vec{x}\in P^\circ(\bQ)$, the above formulation is clearly compatible with Definition \ref{def: CM}.
Arguing as in the proof of Proposition \ref{prop: R-CM}, we could write
$$ \mL(\vx)=\sum_{p} f_p(x_1,...,x_k)\cdot \mL(\va_p),$$
where $f_p(x_1,...,x_k)$'s are finite linear functions. 
Thus the formula $(\star)$ says
$$\lambda_{\CM,(\mX,\mD(\vx),\mL(\vec{x});\pi)}= -\pi_*\left(\sum_{p} f_p(x_1,...,x_k)\cdot (-K_{\mX/T}-\mD(\va_p))\right)^{d+1},$$
which could then be reformulated as the following form:
$$\lambda_{\CM,(\mX,\mD(\vx),\mL(\vec{x});\pi)}= \sum_q g_q(x_1,...,x_k) \cdot M_q, $$
where $g_{q}(x_1,...,x_k)$'s are finite polynomials and $M_{q}$'s are finite 
$\bQ$-line bundles on $T$.
\end{example}

\begin{corollary}\label{cor: R-PT}
Let $(X, \sum_{j=1}^kc_jD_j)$ be a log Fano pair of dimension $d$, where $D_j$'s ($1\leq j\leq k$) are effective Weil divisors on $X$ and $c_j\in [0,1]$ (not necessarily rational). Suppose $\pi: (\mX, \sum_{j=1}^kc_j\mD_j;\mL(\vc))\to \bP^1$ is a special test configuration of $(X, \sum_{j=1}^kc_jD_j; L(\vc):=-K_X-\sum_{j=1}^kc_jD_j)$ (note here that the test configuration is taken to be compact by the natural compactification, see Definition \ref{def: tc}). Then we have
\begin{align*}
\Fut(\mX, \sum_{j=1}^kc_j\mD_j;\mL(\vc))\ &=\  \frac{1}{(d+1)L(\vc)^d}\cdot\wt(\lambda_{\CM,(\mX,\sum_{j=1}^kc_j\mD_j,\mL(\vc);\pi)})\\
&=\ \frac{1}{(d+1)L(\vc)^d}\cdot \wt(-\pi_*(-K_{\mX/\bP^1}-\sum_{j=1}^kc_j\mD_j)^{d+1}).
\end{align*}
\end{corollary}

\begin{proof}
If $(c_1,...,c_k)$ is a rational vector, the formula follows from Proposition \ref{prop: PT}. If $(c_1,...,c_k)$ is not rational, the formula is given by the continuity.
\end{proof}

\subsection{The shape of the K-semistable domain}

In this section, we will show that the K-semistable domain is semi-algebraic under the volume lower bound assumption.

\begin{lemma}\label{lem: GIT}
Given a log Fano pair $(X, \Delta)$ and a fixed real number $0<\delta'_0\leq 1$. Suppose $\delta(X, \Delta)\geq \delta_0'$. To test K-semistability of $(X, \Delta)$, it suffices to confirm the non-negativity of the generalized Futaki invariants of all special test configurations $(\mX, \Delta_\tc)\to \bA^1$ with $\delta(\mX_0, \Delta_{\tc,0})\geq \delta'_0$. 
\end{lemma}

\begin{proof}

We first assume $(X, \Delta)$ is a log pair with $\bQ$-coefficients.
Suppose $\Fut(\mX, \Delta_\tc)\geq 0$ for any special test configuration $(\mX, \Delta_\tc; -K_\mX-\Delta_\tc)\to \bA^1$ of $(X, \Delta;-K_X-\Delta)$ with $\delta(\mX_0, \Delta_{\tc,0})\geq \delta'_0$, we aim to show that $(X, \Delta)$ is K-semistable. If not, then $\delta(X, \Delta)<1$. By \cite{LXZ22}, there exists a prime divisor $E$ over $X$ such that $\delta(X, \Delta)=\frac{A_{X, \Delta}(E)}{S_{X, \Delta}(E)}$. By \cite{BLZ22}, $E$ induces a special test configuration $(\mX', \Delta_\tc')\to \bA^1$ with 
$$\delta(\mX_0', \Delta'_{\tc,0})=\delta(X, \Delta)\geq \delta'_0\quad\text{and} \quad \Fut(\mX', \Delta'_{\tc}; -K_{\mX'}-\Delta'_\tc)<0.$$ 
This leads to a contradiction.

If $(X, \Delta)$ is a log pair with $\bR$-coefficients, by Theorem \ref{thm: R-copy} (1), there still exists a prime divisor $E$ over $X$ computing $\delta(X, \Delta)$,  and this divisor leads to a contradiction by the same way as in the case of $\bQ$-coefficients. 
\end{proof}

Now we are ready to give a characterization of the K-semistable domain. 

\begin{proposition}\label{prop: kss property}
Let $X$ be a normal projective variety of dimension $d$, and $D_j$'s ($1\leq j\leq k$) are effective Weil divisors on $X$. Given a rational polytope $P\subset [0,1]^k$ and assume $P\subset \overline{\LF(X, \sum_{j=1}^kD_j)}$. Suppose there exists a positive real number $\epsilon$ such that $\vol(-K_X-\sum_{j=1}^kx_jD_j)\geq \epsilon$ for any $(x_1,...,x_k)\in P$. 
Then $\Kss(X, \sum_{j=1}^kD_j)_P$ can be characterized by the following way: there exists a finite polytope chamber decomposition $P=\cup_s P_s$, which satisfies that each $P_s$ is a rational polytope and $P_s^\circ\cap P_{s'}^\circ=\emptyset$ if $P_s$ and $P_{s'}$ are adjacent, such that for each index $s$ and each face $F$ of $P_s$, the vector $(x_1,...,x_k)\in F^\circ$ is contained in $\Kss(X, \sum_{j=1}^kD_j)_{F}$ if and only if $(x_1,...,x_k)$ satisfies the following finitely many inequalities indexed by $i$ and $l$:
$$G^{(i)}_{l}(x_1,...,x_k)\geq 0, $$
where $G^{(i)}_{l}(x_1,...,x_k)$'s are $\bQ$-polynomials depending on $s$ and $F$. 
\end{proposition}

\begin{proof}
We divide the proof into several steps.

\

\textbf{Step 0}. \textit{In this step, we define a positive real number $\delta'_0$}. For any $\vec{x}:=(x_1,...,x_k)\in P$, we denote $L_{\vec{x}}:=-K_X-\sum_{j=1}^k x_jD_j$. By assumption, $L_{\vec{x}}$ is a big and nef $\bR$-Cartier divisor for any $\vec{x}\in P$.
Suppose $(a_1,...,a_k)\in P$ satisfies that $(X, \sum_{j=1}^ka_jD_j)$ is a K-semistable log Fano pair, then there exists a positive rational number $\epsilon'$ depending only on $d, \epsilon$ such that $(X, \sum_{j=1}^ka_jD_j)$ is $\epsilon'$-lc (e.g. \cite[Corollary 2.22]{LZ24}). Let $Q$ be a rational polytope consisting of vectors $(x_1,...,x_k)\in P$ satisfying $(X, \sum_{j=1}^kx_jD_j)$ is $\epsilon'$-lc. 
Thus we have
$$\Kss(X, \sum_{j=1}^k D_j)_P\subset Q\cap \LF(X, \sum_{j=1}^kD_j).$$
In the subsequent, we replace $P$ with $Q$. 
By the continuity (e.g. \cite[Proposition 2.10]{LZ24}), there exists a positive real number $\delta'_0$ such that $\delta(X, \sum_{j=1}^kx_j D_j)\geq \delta'_0$ for any $(x_1,...,x_k)\in P$. 
We may assume $\delta_0'<1$, otherwise we have 
$$\Kss(X, \sum_{j=1}^k D_j)_P=P\cap \LF(X, \sum_{j=1}^kD_j).$$

\

\textbf{Step 1}.
\textit{In this step, we define a log bounded set $\mP$ and produce the finite polytope chamber decomposition}. Let $\delta'_0$ be the positive real number as in step 0, we define a set of couples $\mP$.  We say $(Y, \sum_{j=1}^kB_j)$ belongs to $\mP$ if and only if the following conditions are satisfied:
\begin{enumerate}
\item $Y$ is a normal projective variety of dimension $d$ and $B_j$'s ($1\leq j\leq k$) are effective Weil divisors on $Y$;
\item there exists a vector $(a_1,...,a_k)\in P$ such that $(Y, \sum_{j=1}^ka_jB_j)$ is a log Fano pair with $\delta(Y, \sum_{j=1}^ka_jB_j)\geq \delta'_0$;
\item $\vol(-K_Y-\sum_{j=1}^kx_jB_j)\geq \epsilon$ for any $(x_1,...,x_k)\in P$;
\item $(Y, \sum_{j=1}^kx_jB_j)$ is $\epsilon'$-lc for any $(x_1,...,x_k)\in P$ (see step 0 for $\epsilon'$).
\end{enumerate}
By a similar argument as in the proof of Proposition \ref{prop: Abdd}, we know that $\mP$ is log bounded. By Lemma \ref{lem: bdd delta}, there exists a positive real number $\delta_0$ (which could be chosen to be  smaller than $\delta_0'$) such that $\delta(Y, \sum_{j=1}^kx_jB_j)\geq \delta_0$ for any $(Y, \sum_{j=1}^kB_j)\in \mP$ and any $(x_1,...,x_k)\in P$ making $(Y, \sum_{j=1}^kx_jB_j)$ a log Fano pair. 

Let $\mA_1:=\mA_1(\delta_0, \epsilon)$ be defined as in Section \ref{sec: bdd set} and it is clear that $\mP$ is a subset of $\mA_1$.  Let $P=\cup_s P_s$ be a finite rational polytope chamber decomposition satisfying the properties in Proposition \ref{prop: perfect fcd}.
By Proposition \ref{prop: perfect fcd},  for any $(Y, \sum_{j=1}^kB_j)\in \mA_1$ and  any face $F$ of $P_s$, if $(Y, \sum_{j=1}^ka_jB_j)$ is a log Fano pair for some $(a_1,...,a_k)\in F^\circ$, then $(Y, \sum_{j=1}^kx_jB_j)$ is a log Fano pair for any $(x_1,...,x_k)\in F^\circ$.  Note that 
$$\Kss(X, \sum_{j=1}^kD_j)_P=\cup_{F^\circ} \Kss(X, \sum_{j=1}^kD_j)_{F^\circ},$$
where $F$ runs through all faces of all chambers $P_s$'s and
$$\Kss(X, \sum_{j=1}^k D_j)_{F^\circ}:=\{(x_1,...,x_k)\in F^\circ |\ \textit{$(X, \sum_{j=1}^k x_jD_j)$ is a K-semistable log Fano pair}\}.$$
Thus it suffices to figure out the structure of each $\Kss(X, \sum_{j=1}^kD_j)_{F^\circ}$.
In the subsequent steps, we replace $P$ with $F^\circ$. 
Correspondingly, we replace $\mP$ with a subset for which the following condition is satisfied (to replace the second condition in the definition of $\mP$):

(2') there exists a vector $(a_1,...,a_k)\in F^\circ$ such that $(Y, \sum_{j=1}^ka_jB_j)$ is a log Fano pair with $\delta(Y, \sum_{j=1}^ka_jB_j)\geq \delta'_0$.

\

\textbf{Step 2}.
\textit{In this step, we find suitable families parametrizing $\mP$}. Fix a rational vector $\vec{c}:=(c_1,...,c_k)\in F^\circ(\bQ)$. By step 1,  the following set
$$\mP_{\vec{c}}:=\{(Y, \sum_{j=1}^kc_jB_j)\ |\  \text{$(Y, \sum_{j=1}^kB_j)\in \mP$}\} $$
is a log bounded family of log Fano pairs. 
Thus there exists a sufficiently divisible integer $m$ (which does not depend on the choice of $(Y, \sum_{j=1}^kc_jB_j)\in \mP_{\vec{c}}$) such that $-m(K_Y+\sum_{j=1}^kc_jB_j)$ is a very ample line bundle.
By a standard argument on Hilbert scheme and applying Proposition \ref{prop: algebraic locus}, there exist finitely many projective spaces $\bP^{N_i}$ and finitely many families indexed by $i$
$$\pi_i: (\mX^{(i)}, \sum_{j=1}^k\mD^{(i)}_j)\to Z_i$$ 
such that  the following conditions are satisfied:
\begin{enumerate}
\item for each $i$, $\pi_i: \mX^{(i)}\to Z_i$ is a proper flat morphism over a finite type base, and $\mD^{(i)}_j$'s are families of effective Weil divisors such that $(\mX_z^{(i)}, \sum_{j=1}^kc_j\mD^{(i)}_{j,z})$ is log Fano for any $z\in Z_i$; 
\item for each $i$ and each closed point $z\in Z_i$, the fiber $(\mX_z^{(i)}, \sum_{j=1}^k\mD^{(i)}_{j,z})$ is embedded into $\bP^{N_i}$ and $\pi_i$ is the universal family embedded into $\bP^{N_i}\times Z_i$ as follows:
\begin{center}
	\begin{tikzcd}[column sep = 2em, row sep = 2em]
	(\mX^{(i)}, \sum_{j=1}^k\mD^{(i)}_j)  \arrow[rd,"\pi_i",swap] \arrow[rr,"\phi_i"]&& \bP^{N_i}\times Z_i \arrow[ld,""]\\
	 &Z_i&;
	\end{tikzcd}
\end{center}
\item $-m\left(K_{\mX^{(i)}/Z_i}+\sum_{j=1}^kc_j\mD_j^{(i)}\right)$ is linearly equivalent to the pull-back of $\mO_{\bP_{Z_i}^{N_i}}(1)$;
\item  for any $(Y, \sum_{j=1}^kB_j)\in \mP$, there exists some index $i$ such that $(Y, \sum_{j=1}^kB_j)$ is isomorphic to a fiber of $\pi_i$; 
\item for each $i$, the base $Z_i$ admits a $G_i$-action, where $G_i:=\SL(N_i+1)$.
\end{enumerate}

By the construction of $\delta'_0$ and $\mP$ (in step 0 and step 1), we clearly see that $(X, \sum_{j=1}^kD_j)$ is contained in  $\mP$. Up to removing the redundant families, we may assume
there exists a closed point $z^{(i)}_0\in Z_i$ corresponding to $(X, \sum_{j=1}^kD_j)$ for each index $i$. Since we only concern about $z^{(i)}_0\in Z_i$ for each family $\pi_i$, up to a base change, we may replace $Z_i$ with $\overline{G_i\cdot z^{(i)}_0}\cap Z_i$. Up to an equivariant resolution and a further base change, we may assume each $Z_i$ is smooth.

\

\textbf{Step 3}. 
\textit{In this step, we define a polarization depending on the coefficients}. By step 2, we see that $-K_{\mX^{(i)}/Z_i}-\sum_{j=1}^kc_j\mD^{(i)}_j$ is $\bQ$-Cartier and relatively ample for the family $(\mX^{(i)}, \sum_{j=1}^kc_j\mD^{(i)}_j)\to Z_i$. 
By \cite{BLX22, Xu20}, up to a replacement of $Z_i$ with a $G_i$-invariant open subset, we may assume 
$$\delta(\mX^{(i)}_t, \sum_{j=1}^kc_j\mD^{(i)}_{j,t})\geq \delta_0$$ 
for every $t\in Z_i$ (note that $\delta_0$ is given in step 1). By the construction of the family $(\mX^{(i)}, \sum_{j=1}^kc_j\mD^{(i)}_j)\to Z_i$, it is clear that $(\mX^{(i)}_t, \sum_{j=1}^k\mD^{(i)}_{j,t})$ is contained in $\mA_1$ for every $t\in Z_i$.
Thus we have $F^\circ\subset \LF(\mX^{(i)}_t, \sum_{j=1}^k\mD^{(i)}_{j,t})$ for every $t\in Z_i$ by the assumption in step 1. Applying Proposition \ref{prop: uniform polarization}, we know that
$$\mL(\vec{x})_i:=-K_{\mX^{(i)}/Z_i}-\sum_{j=1}^kx_j\mD^{(i)}_j$$
is $\bR$-Cartier and relatively ample for any $\vec{x}:=(x_1,...,x_k)\in F^\circ$. 

\

\textbf{Step 4}.
\textit{In this step, we compute the CM line bundle}.
By step 3, we know that $\mL(\vec{x})_i$ is a relatively ample $\bR$-line bundle on $\mX^{(i)}$ for any $\vec{x}:=(x_1,...,x_k)\in F^\circ$. Denote by:
$$\Lambda(\vec{x})_i:=\lambda_{\CM,(\mX^{(i)},\sum_{j=1}^kx_j\mD^{(i)}_j,\mL(\vec{x})_i;\pi_i)}.$$
By Proposition \ref{prop: R-CM} and Definition \ref{def: R-CM},  $\Lambda(\vec{x})_i$ can be reformulated as follows: 
$$\Lambda(\vec{x})_i=\sum_{q}g^{(i)}_{q}(x_1,...,x_k) M^{(i)}_{q}, $$
where $g^{(i)}_{q}(x_1,...,x_k)$'s are finite polynomials and $M^{(i)}_{q}$'s are finite $G_i$-linearized line bundles on $Z_i$.

\

\textbf{Step 5}.
\textit{In this step, we apply Lemma \ref{lem: GIT} to confirm K-semistability}. For any $(x_1,...,x_k)\in F^\circ$ and any special test configuration of $(X, \sum_{j=1}^kx_jD_j)$ whose central fiber, denoted by $(Y, \sum_{j=1}^kx_jB_j)$, admits delta invariant $\geq \delta'_0$, the couple $(Y, \sum_{j=1}^kB_j)$ is clearly contained in $\mP$ by the construction of $\mP$ in step 1.  Again by step 1,  we have  $\delta(Y, \sum_{j=1}^kx_jB_j)\geq \delta_0$ for any $(x_1,...,x_k)\in F^\circ$.
This implies that there exists an index $i$ such that $(Y, \sum_{j=1}^kB_j)$ is a fiber of the family $(\mX^{(i)}, \sum_{j=1}^k\mD^{(i)}_j)\to Z_i$ constructed in step 3.

For any $z\in Z_i$ and any vector $\vec{x}:=(x_1,...,x_k)\in F^\circ$, let $\rho: \bC^*\to G_i$ be a one parameter subgroup such that $\lim_{t\to 0} \rho(t)\cdot z\in Z_i$, then it induces a special test configuration of $(\mX^{(i)}_z, \sum_{j=1}^kx_j\mD^{(i)}_{j,z};-K_{\mX^{(i)}_z}-\sum_{j=1}^kx_j\mD^{(i)}_{j,z})$, denoted by
$$(\fX_\rho, \sum_{j=1}^kx_j\fD_{\rho, j};-K_{\fX_\rho}-\sum_{j=1}^kx_j\fD_{\rho, j})\to \bA^1. $$ 
By Corollary \ref{cor: R-PT}, we have
\begin{align*}
\Fut(z,\vec{x}, \rho):=&\ \Fut(\fX_\rho, \sum_{j=1}^kx_j\fD_{\rho, j};-K_{\fX_\rho}-\sum_{j=1}^kx_j\fD_{\rho, j})\\
=&\ \frac{1}{d+1}\cdot \frac{\mu^{\Lambda(\vec{x})_i}(z, \rho)}{(-K_{\mX^{(i)}_z}-\sum_{j=1}^kx_j\mD^{(i)}_{j,z})^d}, 
\end{align*}
where $\mu^{\Lambda(\vec{x})_i}(z, \rho)$ is the weight of the $\bC^*$-action (induced by $\rho$) on $\Lambda(\vec{x})_i$.
Let $M^{(i)}$ be a sufficiently ample $G_i$-linearized line bundle on $Z_i$ such that 
$$L^{(i)}_{q}:=M^{(i)}_{q}+M^{(i)}$$ 
is ample for every $q$, then by step 4 we have
$$\Fut(z, \vec{x}, \rho)= \frac{\sum_q g^{(i)}_{q}(x_1,...,x_k)\cdot\left(\mu^{L^{(i)}_{q}}(z, \rho)-\mu^{M^{(i)}}(z, \rho)\right)}{(d+1)\cdot (-K_{\mX^{(i)}_z}-\sum_{j=1}^kx_j\mD^{(i)}_{j,z})^d}.$$
Recall that $z^{(i)}_0\in Z_i$ is the closed point corresponding to $(X, \sum_{j=1}^kD_j)$.
For any $(x_1,...,x_k)\in F^\circ$, by Lemma \ref{lem: GIT}, we see that $(X, \sum_{j=1}^kx_jD_j)$ is K-semistable if and only if 
$$\Fut(z^{(i)}_0, \vec{x}, \rho)\geq 0 $$
for any index $i$ and any one parameter subgroup $\rho:\bC^*\to G_i$ with $\lim_{t\to 0}\rho(t)\cdot z^{(i)}_0\in Z_i$.

\

\textbf{Step 6}. 
\textit{In this step, we conduct an analysis on $\mu^{L^{(i)}_{q}}(z_0^{(i)}, \rho)-\mu^{M^{(i)}}(z_0^{(i)}, \rho)$ and finish the proof}.
By step 5, for any $(x_1,...,x_k)\in F^\circ$, the log Fano pair $(X, \sum_{j=1}^kx_jD_j)$ is K-semistable if and only if
$$\sum_q g^{(i)}_{q}(x_1,...,x_k)\cdot\left(\mu^{L^{(i)}_{q}}(z^{(i)}_0, \rho)-\mu^{M^{(i)}}(z^{(i)}_0, \rho)\right)\geq 0 $$
for any index $i$ and any one parameter subgroup $\rho:\bC^*\to G_i$ with $\lim_{t\to 0}\rho(t)\cdot z^{(i)}_0\in Z_i$.
Fix a maximal torus $\bT_i\subset G_i$. Applying \cite[Lemma A.3]{LWX19}, we know that there exists a constructible decomposition $Z_i=\amalg_I S_{iI}^{\bT_i}$ such that for any $z\in S_{iI}^{\bT_i}$,  $\mu^{L^{(i)}_p}(z, \rho)$ and $\mu^{M^{(i)}}(z, \rho)$ are rational piecewise linear functions on $\Hom_\bQ(\bC^*, \bT_i)$ which are independent of the choice of $z\in S_{iI}^{\bT_i}$. We denote these two functions by $\mu_{iI}^{L^{(i)}_p}(\rho)$ and $\mu_{iI}^{M^{(i)}}(\rho)$ respectively. Note that any one parameter subgroup $\rho$ of $G_i$ is conjugate via some $g\in G_i$ to a one parameter subgroup $g\rho g^{-1}$ of $\bT_i$ (e.g. \cite[Definition A.1]{LWX19}), and by \cite[Lemma A.2]{LWX19} we have
$$\mu^{L^{(i)}_p}(z, \rho)=\mu^{L^{(i)}_p}(gz, g\rho g^{-1})\quad \text{and}\quad  \mu^{M^{(i)}}(z, \rho)=\mu^{M^{(i)}}(gz, g\rho g^{-1}).$$ 
Thus we see that $(X, \sum_{j=1}^kx_jD_j)$ is K-semistable if and only if
$$\sum_q g^{(i)}_q(x_1,...,x_k)\cdot\left(\mu_{iI}^{L^{(i)}_q}(\rho)-\mu_{iI}^{M^{(i)}}(\rho)\right)\geq 0 $$
for any index $i$ and any one parameter subgroup $\rho\in \Hom_\bQ(\bC^*, \bT_i)$ with $\lim_{t\to 0}\rho(t)\cdot z^{(i)}_0\in Z_i$ and any $I$ with $z^{(i)}_0\in S_{iI}^{G_i}$, where $S_{iI}^{G_i}:=G_i\cdot S_{iI}^{\bT_i}$ is a constructible subset of $Z_i$ by Chevalley’s Lemma \cite[Exercise II.3.19]{Har77}. Assume $\bT_i$ has rank $n_i$. Since $\mu_{iI}^{L^{(i)}_q}(\rho)$ and $\mu_{iI}^{M^{(i)}}(\rho)$ are rational piecewise linear functions on $\Hom_\bQ(\bC^*, \bT_i)$, we may write
$$h^{(i)}_q(y_1,...,y_{n_i}):= \mu_{iI}^{L^{(i)}_q}(\rho_{\vec{y}})-\mu_{iI}^{M^{(i)}}(\rho_{\vec{y}})$$
for $\vec{y}:=(y_1,...,y_{n_i})$, where $\rho_{\vec{y}}\in \Hom_\bQ(\bC^*, \bT_i)$ corresponds to $(y_1,...,y_{n_i})\in \bQ^{\oplus n_i}$ and $h^{(i)}_q(y_1,...,y_{n_i})$ is a rational piecewise linear function on $ \bQ^{\oplus n_i}$.
As we only consider those $\rho_{\vec{y}}\in \Hom_\bQ(\bC^*, \bT_i)$ with the property $\lim_{t\to 0}\rho_{\vec{y}}(t)\cdot z^{(i)}_0\in Z_i$, we always restrict to a subset $J_i$ of $\bQ^{\oplus n_i}$ when we consider $h^{(i)}_q(y_1,...,y_{n_i})$.

Above all, restricted to $F^\circ$, the K-semistable domain $\Kss(X, \sum_{j=1}^kD_j)_{F^\circ}$ is cut out by the following equations (indexed by $i$) for any $(y_1,...,y_{n_i})$ contained in the subset $J_i$ of $\bQ^{\oplus n_i}$:
$$ \sum_q g^{(i)}_q(x_1,...,x_k)\cdot h^{(i)}_q(y_1,...,y_{n_i})\geq 0,$$
where $g^{(i)}_q(x_1,...,x_k)$'s are polynomials and $h^{(i)}_q(y_1,...,y_{n_i})$'s
 are rational piecewise linear functions on $\bQ^{\oplus n_i}$. Since $h^{(i)}_q(y_1,...,y_{n_i})$'s are piecewise linear (e.g. \cite[Lemma A.3]{LWX19}), it is not hard to see that the above equation is determined by finitely many points in $J_i\subset \bQ^{\oplus n_i}$ indexed by $l$, denoted by $(y_{l1}^{(i)},...,y_{ln_i}^{(i)})$. Write
 $$G^{(i)}_l(x_1,...,x_k):= \sum_q g^{(i)}_q(x_1,...,x_k)\cdot h^{(i)}_q(y_{l1}^{(i)},...,y_{ln_i}^{(i)}).$$ 
 Then we see that $(x_1,...,x_k)\in F^\circ$ is contained in $\Kss(X, \sum_{j=1}^kD_j)_{F^\circ}$ if and only if $(x_1,...,x_k)$ satisfies the following finite equations indexed by $i$ and $l$:
$$G^{(i)}_{l}(x_1,...,x_k)\geq 0.$$
The proof is complete.
\end{proof}

The following lemma is applied in step 1 of the above proof.

\begin{lemma}\label{lem: bdd delta}
Given two positive integers $d$ and $k$, a rational polytope $P\subset [0,1]^k$, and two positive real numbers $\epsilon, \delta_0'$. Let $\mP$ be a set of couples as defined in step 1 of the proof of Proposition \ref{prop: kss property}. Then there exists a positive real number $\delta_0$ depending only on the given data $d,k,P, \epsilon, \delta_0'$ such that $\delta(Y, \sum_{j=1}^kx_jB_j)\geq \delta_0$ for any $(Y, \sum_{j=1}^kB_j)\in \mP$ and any $(x_1,...,x_k)\in P$ making $(Y, \sum_{j=1}^kx_jB_j)$ a log Fano pair.
\end{lemma}

\begin{proof}
As we have seen in step 1 of the proof of Proposition \ref{prop: kss property}, $\mP$ is log bounded. By Proposition \ref{prop: algebraic locus}, we could find a family of couples, denoted by $(\mX, \sum_{j=1}^k\mD_j)\to T$, satisfying the following conditions:
\begin{enumerate}
\item the base $T$ is a finite disjoint union of smooth varieties;
\item $\mX\to T$ is a flat family of normal projective varieties and and $\mD_j\to T$ is a family of effective Weil divisors for each $1\leq j\leq k$;
\item for any $t\in T$, $(\mX_t, \sum_{j=1}^kc_j\mD_{j,t})$ is log Fano for some numbers $c_j$'s;
\item the set $\mP$ is contained in the set of fibers of $(\mX, \sum_{j=1}^k\mD_j)\to T$.
\end{enumerate}
By Proposition \ref{prop: LF stratification}, there exists a constructible stratification $T=\amalg_i T_i$ such that the log Fano domain (restricted to $P$) of $(\mX_t, \sum_{j=1}^k\mD_{j,t})$ does not depend on $t\in T_i$. We denote this restricted log Fano domain to be $P_i$. Up to a further constructible stratification of each $T_i$, we may assume each $T_i$ is smooth. Let $(\mX_i, \sum_{j=1}^k\mD_{ij})\to T_i$ be the pull-back of $(\mX, \sum_{j=1}^k\mD_j)\to T$ along $T_i\to T$. Then we see the fibers of $(\mX_i, \sum_{j=1}^k\mD_{ij})\to T_i$ share the same log Fano domain restricted to $P$, which is $P_i$ by our notation.

It suffices to show that there exists a positive real number $\delta_0$ depending only on the given data $d,k,P,\epsilon,\delta_0'$ such that 
$\delta(\mX_{i,t}, \sum_{j=1}^kx_j\mD_{ij, t})\geq \delta_0$
for any $(x_1,...,x_k)\in P_i$ and any $t\in T_i$. 
Replace $(\mX, \sum_{j=1}^k\mD_j)\to T$ (resp. $P$) with $(\mX_i, \sum_{j=1}^k\mD_{ij})\to T_i$ (resp. $P_i$). By the proof of Theorem \ref{thm: finiteness}, there exists a finite decomposition $T=\amalg_r U_r$ such that the following invariant
$$\tilde{\delta}(\mX_t, \sum_{j=1}^k x_j\mD_{j,t}):=\min\{1,\ \delta(\mX_t, \sum_{j=1}^k x_j\mD_{j,t})\} $$
does not depend on $t\in U_r$ for any $(x_1,...,x_k)\in P$. For each $r$, we arbitrarily choose $t_r\in U_r$. By the continuity of delta invariants (e.g. \cite[Proposition 2.10]{LZ24}), there exists a positive real number $\mu_r<1$ such that $\delta(\mX_{t_r}, \sum_{j=1}^kx_j\mD_{j,t_r})\geq \mu_r$ for any $(x_1,...,x_k)\in P$. Define
$$\delta_0:= \min_{r}\{\mu_r\}. $$
It is clear that $\delta_0$ is the desired positive real number. 
\end{proof}

Proposition \ref{prop: kss property} gives a clear description of the shape of the K-semistable domain. Roughly speaking, if we take $P$ to be $\overline{\LF(X, \sum_{j=1}^kD_j)}$, then there exists a finite rational polytope chamber decomposition of $P$ such that the K-semistable domain restricted to each face of each chamber is cut out by some algebraic functions. 
To characterize the shape of the K-semistable domain in formal language, we introduce the following concept of semi-algebraic sets. 

\begin{definition}\label{def: semi-algebraic}
A subset of $\bR^k$ is called \textit{a basic semi-algebraic set} if it is defined by finite polynomial equalities and finite polynomial inequalities. \textit{A semi-algebraic set} is defined to be the finite union of basic semi-algebraic sets.
\end{definition}

A polynomial equality is of the form $F(x_1,...,x_k)=0$ and a polynomial inequality is of the form $G(x_1,...,x_k)>0$, where $F, G$ are polynomials.
We have the following basic properties for semi-algebraic sets.

\begin{proposition}\label{prop: semi-algebraic property}
Let $Z_1, Z_2$ be two semi-algebraic sets in $\bR^k$, then we have the following conclusions:
\begin{enumerate}
\item $Z_1\cup Z_2$ is semi-algebraic;
\item $Z_1\cap Z_2$ is semi-algebraic;
\item the complement of $Z_1$ (resp. $Z_2$) is semi-algebraic;
\item the closure (under Euclidean topology) of $Z_1$ (resp. $Z_2$) is semi-algebraic;
\item $Z_1$ (resp. $Z_2$) can be decomposed into disjoint union of a finite number of semi-algebraic sets, and each of them is homeomorphic to an open hypercube $(0,1)^r$ for some integer $r\in \bN$ (with $(0,1)^0$ being a point).
\end{enumerate}
\end{proposition}

\begin{proof}
The first four statements are clear by the definition of semi-algebraic sets. The last statement is referred to \cite[Theorem 2.3.6]{BCR98}.
\end{proof}

\begin{corollary}\label{cor: semi-algebraic kss}
Let $X$ be a normal projective variety of dimension $d$, and $D_j$'s ($1\leq j\leq k$) are effective Weil divisors on $X$. Given a rational polytope $P\subset [0,1]^k$ and assume $P\subset \overline{\LF(X, \sum_{j=1}^kD_j)}$. Suppose there exists a positive real number $\epsilon$ such that $\vol(-K_X-\sum_{j=1}^kx_jD_j)\geq \epsilon$ for any $(x_1,...,x_k)\in P$. 
Then $\Kss(X, \sum_{j=1}^kD_j)_P$ is a semi-algebraic set which can be decomposed into a disjoint union of a finite number of semi-algebraic sets, and each of them is homeomorphic to an open hypercube $(0,1)^r$ for some integer $r\in \bN$ (with $(0,1)^0$ being a point).
\end{corollary}

\begin{proof}
The proof is a combination of Propositions \ref{prop: kss property} and \ref{prop: semi-algebraic property}.
\end{proof}

\begin{remark}
One may notice that there is a condition on volume in Corollary \ref{cor: semi-algebraic kss}, i.e. $\vol(-K_X-\sum_{j=1}^kx_jD_j)\geq \epsilon$ for any $(x_1,...,x_k)\in P$. This condition is natural in many settings (e.g. Example \ref{example: P_2} and Example \ref{example: p1p2}) and it is essential in the proof of Proposition \ref{prop: kss property}.  But it is expected to have the same conclusion (on the shape of the K-semistable domain) without the volume condition. This seems difficult and will be explored in future works.
\end{remark}

\section{Wall crossing for K-moduli}\label{sec: fcd}

In this section, we aim to establish the wall crossing theory for K-moduli in non-proportional setting. Throughout the section, we fix two positive integers $d,k$ and a rational polytope $P\subset [0,1]^k$. Let us first recall the set of couples $\mG:=\mG(d,k, P)$ depending on $d,k$ and $P$.  We say a couple $(X, \sum_{j=1}^kD_j)$ belongs to $\mG$ if and only if the following conditions are satisfied:
\begin{enumerate}
\item $X$ is a normal projective variety of dimension $d$;
\item $D_j$'s ($1\leq j\leq k$) are effective Weil divisors on $X$;
\item there exists a vector $(c_1,...,c_k)\in P$ such that $(X, \sum_{j=1}^kc_jD_j)$ is a K-semistable log Fano pair.
\end{enumerate}
Recall that the restricted K-semistable domain for $(X, \sum_{j=1}^kD_j)\in \mG$ is defined as follows: 
$$\Kss(X, \sum_{j=1}^kD_j)_P:=\{(x_1,...,x_k)\in P | \textit{$(X, \sum_{j=1}^kx_jD_j)$ is a K-semistable log Fano pair}\}.$$

\subsection{Semi-algebraic chamber decomposition}\label{subsec: sacd}
As $\mG$ is not necessarily log bounded, we will fix a log bounded subset $\widetilde{\mG}$ in this subsection.

\begin{theorem}\label{thm: fcd}
Suppose $\widetilde{\mG}\subset \mG$ is a subset which is log bounded and there exists a positive real number $\epsilon$ such that $\vol(-K_X-\sum_{j=1}^kx_jD_j)\geq \epsilon$ for any $(X, \sum_{j=1}^kD_j)\in \widetilde{\mG}$ and any $(x_1,...,x_k)\in P$. Then there exists a  decomposition of $P$ into finite disjoint pieces (depending only on $\widetilde{\mG}$), denoted by
$$P=\amalg_{i=0}^m A_i $$
such that the following conditions are satisfied:
\begin{enumerate}
\item each $A_i$ is a semi-algebraic set;
\item the K-semistability of $(X, \sum_{j=1}^kx_jD_j)$ does not change as $(x_1,...,x_k)$ varies in $A_i$ for any $0\leq i\leq m$ and any $(X, \sum_{j=1}^kD_j)\in \widetilde{\mG}$.
\end{enumerate}
\end{theorem}

\begin{proof}
By Theorem \ref{thm: finiteness P}, there exist finitely many locally closed subsets (indexed by $r$) $K_r\subset P$ which could appear as the restricted K-semistable domains.  More precisely, for any $(X, \sum_{j=1}^kD_j)\in \widetilde{\mG}$, there exists some $r$ such that $\Kss(X, \sum_{j=1}^kD_j)_P=K_r$. 
By Corollary \ref{cor: semi-algebraic kss}, each $K_r$ is semi-algebraic. Applying Proposition \ref{prop: semi-algebraic property}, one could easily cook up a desired disjoint union of semi-algebraic sets.
\end{proof}

In the case of one boundary divisor, we have the following result.

\begin{corollary}\label{cor: fcd}
Suppose $k=1$ and $P:=[a,b]\subset [0,1]$. Suppose $\widetilde{\mG}\subset \mG$ is a subset which is log bounded and there exists a positive real number $\epsilon$ such that $\vol(-K_X-bD)\geq \epsilon$ for any $(X, D)\in \widetilde{\mG}$. Then there exist finitely many algebraic numbers (depending only on $\widetilde{\mG}$)
$$a=w_0<w_1<w_2<...<w_p<w_{p+1}=b $$
such that the K-semistability of $(X, xD)$ does not change as $x$ varies in $(w_i, w_{i+1})$ for any $0\leq i\leq p$ and any $(X, D)\in \widetilde{\mG}$.
\end{corollary}

\begin{proof}
It is directly implied by Theorem \ref{thm: fcd} and Proposition \ref{prop: kss property}.
\end{proof}

In the proportional situation, the chamber decomposition in Corollary \ref{cor: fcd} would easily lead to a wall crossing theory for K-moduli (e.g. \cite{ADL19, Zhou23}), as every K-semistable chamber is also a K-polystable chamber. More precisely, we have the following result.

\begin{proposition}\label{prop: linear futaki}
Let $X$ be a Fano variety and $D$ an effective $\bQ$-divisor on $X$ satisfying $D\sim_\bQ -K_X$. Suppose $(X, cD)$ is a K-semistable log Fano pair for any $c\in (a,b)\cap \bQ$ and $(X, c_0D)$ is a K-polystable log Fano pair for some $c_0\in (a,b)\cap \bQ$, then $(X, cD)$ is K-polystable for any $c\in (a,b)\cap \bQ$.
\end{proposition}

\begin{proof}
Taking any rational number $c$ satisfying $a<c<c_0$, we show that $(X, cD)$ is K-polystable. Note that $(X, cD)$ is K-semistable by assumption. Suppose it is not K-polystable, then there exists a non-product type test configuration $(\mX, c\mD; -K_\mX-c\mD)\to \bA^1$ which degenerates $(X, cD)$ into a K-polystable log Fano pair $(\mX_0, c\mD_0)$. Note that 
$$(\mX, a\mD; -K_\mX-a\mD)\to \bA^1\quad \text{and}\quad (\mX, c_0\mD; -K_\mX-c_0\mD)\to \bA^1$$
are also test configurations of $(X, aD)$ and  $(X, c_0D)$
respectively. By interpolation of generalized Futaki invariants (e.g. \cite[Lemma 2.7]{Zhou23}), we have
$$\Fut(\mX, c\mD)=\lambda\cdot \Fut(\mX, a\mD)+(1-\lambda)\cdot \Fut(\mX, c_0\mD)=0$$
for $\lambda=\frac{c_0-c}{c_0-a}$. Since 
$$\Fut(\mX, a\mD)\geq 0\quad\text{ and }
\quad \Fut(\mX, c_0\mD)\geq 0$$ 
by K-semistability, we must have 
$$\Fut(\mX, a\mD)= 0\quad\text{ and }\quad\Fut(\mX, c_0\mD)= 0.$$ 
However, since $(X, c_0D)$ is K-polystable, the equality $\Fut(\mX, c_0\mD)= 0$ implies that $(\mX, c\mD)\to \bA^1$ is a test configuration of product type (e.g. \cite{LWX21}), leading to a contradiction. The contradiction tells that $(X, cD)$ is K-polystable for $c\in (a, c_0)\cap \bQ$. One can similarly conclude it for $c\in (c_0, b)\cap \bQ$.
\end{proof}

However, in non-proportional setting, there is no reason to still have the above proposition without the interpolation property as applied in the proof. The natural idea is to find a refined chamber decomposition of the one obtained in Corollary \ref{cor: fcd} to make sure that each chamber is a K-polystable chamber and the decomposition is still finite. Before achieving this, 
we first define a larger set $\widehat{\mG}$ than $\widetilde{\mG}$ to make it complete in moduli sense.

\begin{definition}\label{def: G-complete}
We say a couple $(Y, \sum_{j=1}^kB_j)\in \mG$ is contained in $\widehat{\mG}$ if and only if there exists a family of couples $(\mY, \sum_{j=1}^k\mB_j)\to C\ni 0$ over a smooth pointed curve such that the following conditions are satisfied:
\begin{enumerate}
\item $(\mY_t, \sum_{j=1}^k \mB_{j,t})\in \widetilde{\mG}$ for $t\ne 0$;
\item $(\mY_0, \sum_{j=1}^k \mB_{j,0})\cong (Y, \sum_{j=1}^kB_j)$;
\item for some $(c_1,...,c_k)\in P$, $(\mY, \sum_{j=1}^k c_j\mB_j)\to C$ is a family of K-semistable log Fano pairs.
\end{enumerate}
\end{definition}
In other words, $\widehat{\mG}$ is obtained by adding to $\widetilde{\mG}$ all the K-semistable degenerations of K-semistable log Fano pairs arising from $\widetilde{\mG}$.

\begin{proposition}\label{prop: G-complete}
Let $\widetilde{\mG}\subset \mG$ be a subset which is log bounded. Suppose there exists a positive real number $\epsilon$ such that $\vol(-K_X-\sum_{j=1}^kx_jD_j)\geq \epsilon$ for any $(X, \sum_{j=1}^kD_j)\in \widetilde{\mG}$ and any $(x_1,...,x_k)\in P$, then the set $\widehat{\mG}$ is log bounded.
\end{proposition}

\begin{proof}
 It is implied by Proposition \ref{prop: Abdd}.
 \end{proof}

\begin{theorem}\label{thm: kps fcd}
Suppose $\widetilde{\mG}\subset \mG$ is a subset which is log bounded and there exists a positive real number $\epsilon$ such that $\vol(-K_X-\sum_{j=1}^kx_jD_j)\geq \epsilon$ for any $(X, \sum_{j=1}^kD_j)\in \widetilde{\mG}$ and any $(x_1,...,x_k)\in P$. Then there exists a  decomposition of $P$ into finite disjoint subsets (depending only on $\widetilde{\mG}$), denoted by
$$P=\amalg_{i'=1}^{m'} A'_{i'} $$
such that the following conditions are satisfied:
\begin{enumerate}
\item each $A'_{i'}$ is a semi-algebraic set homeomorphic to $(0,1)^r$ for some $r\in \bN$ (with $(0,1)^0$ being a point);
\item the K-polystability of $(X, \sum_{j=1}^kx_jD_j)$ does not change as $(x_1,...,x_k)$ varies in $A'_{i'}$ for any $0\leq i'\leq m'$ and any $(X, \sum_{j=1}^kD_j)\in \widetilde{\mG}$.
\end{enumerate}
\end{theorem}

\begin{proof}
By Theorem \ref{thm: fcd}, there exists a  decomposition of $P$ into finite disjoint pieces (depending only on $\widetilde{\mG}$), denoted by
$$P=\amalg_{i=1}^{m} A_{i} $$
such that each $A_i$ is semi-algebraic and the K-semistability of $(X, \sum_{j=1}^kx_jD_j)$ does not change as $(x_1,...,x_k)$ varies in $A_i$ for any $0\leq i\leq m$ and any $(X, \sum_{j=1}^kD_j)\in \widetilde{\mG}$. 
In the sequel, we will automatically apply Proposition \ref{prop: semi-algebraic property} to refine the above decomposition. Replacing $\widetilde{\mG}$ with $\widehat{\mG}$, up to a refinement of the above decomposition, we may assume the following condition by Proposition \ref{prop: perfect fcd}: \textit{for any $0\leq i\leq m$ and any $(X,\sum_{j=1}^kD_j)\in \widetilde{\mG}$, if $(X,\sum_{j=1}^kx_jD_j)$ is a log Fano pair for some $(x_1,...,x_k)\in A_i$, then $(X,\sum_{j=1}^kx_jD_j)$ is a log Fano pair for any $(x_1,...,x_k)\in A_i$}.

To find the suitable further refinement of the decomposition, we first define a set $\mP:=\mP(d, k,P, \epsilon)$ of couples by the following way. We say a couple $(X', \sum_{j=1}^kD_j')$ belongs to $\mP$ if and only if the following conditions are satisfied: 
\begin{enumerate}
\item $X'$ is a normal projective variety of dimension $d$;
\item $D'_j$'s ($1\leq j\leq k$) are effective Weil  divisors on $X'$;
\item $(X', \sum_{j=1}^kc_jD_j')$ is a K-semistable weak log Fano pair for some $(c_1,...,c_k)\in P$;
\item $\vol(-K_{X'}-\sum_{j=1}^kx_jD_j')\geq \epsilon$ for any $(x_1,...,x_k)\in P$.
\end{enumerate}
By Proposition \ref{prop: Abdd}, $\mP$ is log bounded and $\widetilde{\mG}$ is a subset of $\mP$. Applying Propositions \ref{prop: perfect fcd}, up to a further refinement of 
the decomposition of $P$,
we may assume the following condition: \textit{for any $0\leq i\leq m$ and any $(X', \sum_{j=1}^kD_j')\in \mP$, if $(X', \sum_{j=1}^kx_jD_j')$ is a weak log Fano pair for some $(x_1,...,x_k)\in A_i$, then $(X', \sum_{j=1}^kx_jD_j')$ is a weak log Fano pair for any $(x_1,...,x_k)\in A_i$.}

We aim to show that, for any $(X, \sum_{j=1}^kD_j)\in \widetilde{\mG}$, the K-polystability of $(X, \sum_{j=1}^kx_jD_j)$ does not change as $(x_1,...,x_k)$ varies in $A_i$ for any $0\leq i\leq m$. Suppose $(X, \sum_{j=1}^kc_jD_j)$ is K-polystable for some $(c_1,...,c_k)\in A_i$, then $(X, \sum_{j=1}^kx_jD_j)$ is K-semistable for any $(x_1,...,x_k)\in A_i$ by assumption. If $(X, \sum_{j=1}^kc_j'D_j)$ is K-semistable but not K-polystable for some $(c_1',...,c_k')\in A_i$, then 
there exists a non-product type special test configuration
$$(\mX, \sum_{j=1}^kc_j'\mD_j)\to \bA^1, $$
which degenerates $(X, \sum_{j=1}^kc_j'D_j)$ into its K-polystable degeneration $(\mX_0, \sum_{j=1}^kc_j'\mD_{j,0})$ (e.g. Theorem \ref{thm: R-copy} (2)). On the other hand, by Corollary \ref{prop: perfect tc}, the following degeneration by changing the coefficients from $(c'_1,...,c_k')$ to $(c_1,...,c_k)$
$$(\mX, \sum_{j=1}^kc_j\mD_j)\to \bA^1 $$
is still a special test configuration by degenerating the K-polystable log Fano pair $(X, \sum_{j=1}^kc_jD_j)$ into another K-semistable log Fano pair $(\mX_0, \sum_{j=1}^kc_j\mD_{j,0})$. This implies that $(\mX, \sum_{j=1}^kc_j\mD_j)\to \bA^1$ is of product type, leading to a contradiction. The contradiction tells that the K-polystability of $(X, \sum_{j=1}^kx_jD_j)$ does not change as $(x_1,...,x_k)$ varies in $A_i$ for any $0\leq i\leq m$.

Again by Proposition \ref{prop: semi-algebraic property}, we may further decompose each $A_i$ into finite disjoint union of semi-algebraic sets where each of them is homeomorphic to $(0,1)^r$ for some $r\in \bN$.
\end{proof}

\begin{corollary}\label{cor: kps fcd}
Suppose $k=1$ and $P:=[a,b]\subset [0,1]$. Suppose $\widetilde{\mG}\subset \mG$ is a subset which is log bounded and there exists a positive real number $\epsilon$ such that $\vol(-K_X-bD)\geq \epsilon$ for any $(X, D)\in \widetilde{\mG}$. Then there exist finitely many algebraic numbers (depending only on $\widetilde{\mG}$)
$$a=w'_0<w'_1<w'_2<...<w'_q<w'_{q+1}=b $$
such that the K-polystability of $(X, xD)$ does not change as $x$ varies in $(w'_i, w'_{i+1})$ for any $0\leq i\leq q$ and any $(X, D)\in \widetilde{\mG}$.
\end{corollary}

\begin{proof}
It is directly implied by Theorem \ref{thm: kps fcd} and Proposition \ref{prop: kss property}.
\end{proof}

\begin{remark}
We emphasize two points towards Corollary \ref{cor: kps fcd} in this remark. First, the decomposition in Corollary \ref{cor: kps fcd} may not be the minimal one, i.e. the number of $w_i'$'s is not necessarily minimal to satisfy the invariance of K-polystability. However, one can always remove the redundant terms to get the minimal one. Second, the walls $w_i'$'s are not necessarily rational as in the proportional setting (although they are always algebraic), e.g. \cite{ADL19, Zhou23}. For instance, an irrational example is worked out in \cite{DJ+24} recently. 
\end{remark}

\subsection{Wall crossing}
Given the log bounded subset $\widetilde{\mG}\subset \mG$ and its moduli completion $\widehat{\mG}\subset \mG$ as in Section \ref{subsec: sacd}.  By \cite{LXZ22, LZ24}, for each  $\vec{x}:=(x_1,...,x_k)\in P$, there exists a finite type Artin stack $\mM^{\Kss}_{\widehat{\mG},\vec{x}}$ parametrizing all K-semistable log Fano pairs $(X, \sum_{j=1}^kx_jD_j)$, where $(X, \sum_{j=1}^kD_j)\in \widehat{\mG}$. Moreover, the stack descends to a separated and proper good moduli space $M^{\Kps}_{\widehat{\mG},\vec{x}}$, which parametrizes K-polystable objects. Let $\mM^{\Kss, \red}_{\widehat{\mG},\vec{x}}$ be the moduli functor associated to reduced schemes and  $M^{\Kps,\red}_{\widehat{\mG}, \vec{x}}$ the corresponding good moduli space.  Recalling in \cite[Section 5.3]{LZ24}, for a reduced scheme $S$, we define $\mM^{\Kss, \red}_{\widehat{\mG},\vec{x}}(S)$ to be the set of families $(\mX, \sum_{j=1}^k x_j\mD_j)\to S$ satisfying:
\begin{enumerate}
\item $\mX\to S$ is proper and flat;
\item $\mD_j$ is a family of Weil divisors on $\mX$ for each $1\leq j\leq k$;
\item $-K_{\mX/S}-\sum_{j=1}^k x_j\mD_j$ is $\bR$-Cartier;
\item for each $s\in S$, the fiber $(\mX_s, \sum_{j=1}^kx_j\mD_{j,s})$ is a K-semistable log Fano pair with $(\mX_s, \sum_{j=1}^k\mD_{j,s})\in \widehat{\mG}$.
\end{enumerate}
Write $\phi_{\vec{x}}: \mM^{\Kss,\red}_{\widehat{\mG},\vec{x}}\to M^{\Kps,\red}_{\widehat{\mG},\vec{x}}$. We have the following result.

\begin{theorem}\label{thm: moduli invariance}
Let $A_{i'}'$ be a chamber in the chamber decomposition in Theorem \ref{thm: kps fcd}. Then $\phi_{\vec{x}}: \mM^{\Kss,\red}_{\widehat{\mG},\vec{x}}\to M^{\Kps,\red}_{\widehat{\mG},\vec{x}}$ does not change as $\vec{x}$ varies in $A_{i'}'$.
\end{theorem}

\begin{proof}
As we only consider reduced bases, by Proposition \ref{prop: uniform polarization} and \cite[Corollary 4.35]{Kollar23}, the third condition in the definition of $\mM^{\Kss, \red}_{\widehat{\mG},\vec{x}}(S)$ is preserved as we vary $\vec{x}$ in $A_{i'}'$. Thus Theorem \ref{thm: kps fcd} tells that $\phi_{\vec{x}}$ does not change as $\vec{x}$ varies in $A'_{i'}$.
\end{proof}

We also point out the difference with the proportional setting when crossing the wall.
Let $\vec{c}:=(c_1,...,c_k), \vec{c'}:=(c_1',...,c_k')$ be two vectors satisfying $\vec{c}\in A'_{i'}$ and $\vec{c'}\in \partial\bar{A'_{i'}}$ (i,e, $\vec{c'}$ is a  vector lying on the boundary of the closure of $A'_{i'}$). In the proportional setting, we have an open embedding $\mM^{\Kss,\red}_{\widehat{\mG},\vec{c}} \to \mM^{\Kss,\red}_{\widehat{\mG},\vec{c'}}$ (e.g. \cite{ADL19, Zhou23, Zhou23b}). However, in non-proportional setting, it is not necessarily an open embedding as the log pairs may not be log Fano (but keep being weak log Fano) after changing the coefficients. For example, suppose $(c_1',...,c_k')$ is rational and we take a family $(\mX, \sum_{j=1}^kc_j\mD_j)\to S$ in $\mM^{\Kss,\red}_{\widehat{\mG},\vec{c}}(S)$, where $S$ is a reduced scheme. Changing the coefficients from $\vec{c}$ to $\vec{c'}$, it is clear that $(\mX, \sum_{j=1}^kc'_j\mD_j)\to S$ is a family of K-semistable weak log Fano pairs by the construction of the chamber decomposition in Theorem \ref{thm: kps fcd}. If $S$ is smooth, by the proof of Proposition \ref{prop: uniform polarization}, we see that the anti-canonical model of $(\mX, \sum_{j=1}^kc'_j\mD_j)\to S$ over $S$, denoted by 
\begin{center}
	\begin{tikzcd}[column sep = 2em, row sep = 2em]
	 (\mX, \sum_{j=1}^kc'_j\mD_j) \arrow[rd,"\pi",swap] \arrow[rr,""]&& (\mX^\ac, \sum_{j=1}^kc'_j\mD^\ac_j) \arrow[ld,"\pi'"]\\
	 &S&,
	\end{tikzcd}
\end{center}
is a \textit{simultaneous  anti-canonical model}, i.e. restricting to any fiber, $(\mX_s^\ac, \sum_{j=1}^kc'_j\mD^\ac_{j,s})$ is the anti-canonical model of $(\mX_s, \sum_{j=1}^kc'_j\mD_{j,s})$. Thus we have 
$$[(\mX^\ac, \sum_{j=1}^kc'_j\mD^\ac_j)\to S]\in \mM^{\Kss,\red}_{\widehat{\mG},\vec{c'}}(S).$$
If $S$ is not necessarily smooth, we have the following result.

\begin{proposition}\label{prop: ac-model}
Notation as above. Given $[(\mX, \sum_{j=1}^kc_j\mD_j)\to S]\in \mM^{\Kss,\red}_{\widehat{\mG},\vec{c}}(S)$, where $S$ is a reduced scheme. Suppose $\vec{c'}=(c_1',...,c_k')$ is rational. Then there exists a unique element $[(\mY, \sum_{j=1}^kc'_j\mB_j)\to S]\in \mM^{\Kss,\red}_{\widehat{\mG},\vec{c'}}(S)$ such that $(\mY, \sum_{j=1}^kc'_j\mB_j)\to S$ is a simultaneous anti-canonical model of $(\mX, \sum_{j=1}^kc'_j\mD_j)\to S$.
\end{proposition}


\begin{proof}
Since $-K_{\mX/S}-\sum_{j=1}^kc_j\mD_j$ is $\bR$-Cartier and relatively ample, one can easily find an ample $\bR$-divisor $\mH$ on $\mX$ such that 
$$\mH\sim_{\bR} -K_{\mX/S}-\sum_{j=1}^kc_j\mD_j$$ 
over $S$ and 
$(\mX, \sum_{j=1}^kc_j\mD_j+\mH)\to S$ is a family of klt log Calabi-Yau pairs. By Proposition \ref{prop: uniform polarization} and \cite[Corollary 4.35]{Kollar23}, $-K_{\mX/S}-\sum_{j=1}^kc'_j\mD_j$ is $\bQ$-Cartier and there exists a sufficiently divisible integer $m\in \bZ_+$ such that $-m(K_{\mX_s}+\sum_{j=1}^k c_j'\mD_{j,s})$ is a free line bundle for any $s\in S$. Note that $\pi: \mX\to S$ is flat and 
$$H^l(\mX_s, -m(K_{\mX_s}+\sum_{j=1}^k c_j'\mD_{j,s}))=0$$ 
for any $l>0$ (by Kawamata-Viehweg vanishing). Thus $\pi_*(-m(K_{\mX/S}+\sum_{j=1}^kc'_j\mD_j))$ is a vector bundle on $S$ and $-m(K_{\mX/S}+\sum_{j=1}^kc'_j\mD_j)$ is free over $S$. Therefore, one can find an effective $\bQ$-Cartier divisor $\mH'$ such that 
$$\mH'\sim_\bQ -K_{\mX/S}-\sum_{j=1}^kc_j'\mD_j$$ 
over $S$ and $(\mX, \sum_{j=1}^kc_j\mD_j+\mH+\mH')\to S$ is a locally stable family (e.g. \cite[Definition 4.7]{Kollar23}) of klt log pairs such that $K_{\mX/S}+\sum_{j=1}^kc_j\mD_{j}+\mH+\mH'$ is semi-ample over $S$. By \cite[Lemma 4.9]{MZ23}, $(\mX, \sum_{j=1}^kc_j\mD_j+\mH+\mH')\to S$ admits an ample model over $S$ which is a stable family (e.g. \cite[Definition 4.7]{Kollar23}).
Thus $(\mX, \sum_{j=1}^kc'_j\mD_j)\to S$ admits a simultaneous anti-canonical model by \cite[Lemma 4.4]{MZ23}, denoted by $(\mY, \sum_{j=1}^kc'_j\mB_j)\to S$. As a family of log Fano pairs which is also locally stable, it is not hard to confirm $[(\mY, \sum_{j=1}^kc'_j\mB_j)\to S]\in \mM^{\Kss,\red}_{\widehat{\mG},\vec{c'}}(S)$ (e.g. \cite[Lemma 10.58]{Kollar23}). The uniqueness is clear.
\end{proof}

We denote $[(\mY, \sum_{j=1}^kc'_j\mB_j)\to S]$ in Proposition \ref{prop: ac-model} by $[(\mX^\ac, \sum_{j=1}^kc'_j\mD^\ac_j)\to S]$.
To describe the wall crossing diagram, we fix a semi-algebraic chamber $A_i$ in Theorem \ref{thm: kps fcd} and denote $\partial \bar{A_i}$ to be the boundary of the closure of $A_i$. Fix a vector $\fa_0:=(a_{01},...,a_{0k})\in \partial \bar{A_i}$ and a vector $\fa=(a_1,...,a_k)\in A_i$, which are not necessarily rational. 

\begin{theorem}\label{thm: diagram}
Notation as above, there exists a wall crossing diagram as follows:
\begin{center}
	\begin{tikzcd}[column sep = 2em, row sep = 1.5em]
	 \mM^{\Kss,\red}_{\widehat{\mG}, \fa} \arrow[d,"\phi_\fa", swap] \arrow[rr,"\Psi_{\fa, \fa_0}"]&& \mM^{\Kss,\red}_{\widehat{\mG}, \fa_0} \arrow[d,"\phi_{\fa_0}",swap]\\
	 M^{\Kps,\red}_{\widehat{\mG}, \fa}\arrow[rr,"\psi_{\fa, \fa_0}"]&& M^{\Kps,\red}_{\widehat{\mG}, \fa_0}
	\end{tikzcd}
\end{center}
where 
$\psi_{\fa, \fa_0}$ is an induced proper morphism by the universal property of good moduli spaces. If the $\bQ$-span of $\fa_0$ contains the chamber $A_i$, then $\Psi_{\fa, \fa_0}$ is an open embedding.
\end{theorem}

\begin{proof}
We may assume $\mM^{\Kss,\red}_{\widehat{\mG}, \fa}$ is not empty. We first show the existence of the diagram, for which it is enough to establish the morphism $\Psi_{\fa, \fa_0}$. Given an element $[(\mX, \sum_{j=1}^k a_j\mD_j)\to T]\in \mM^{\Kss,\red}_{\widehat{\mG}, \fa}(T)$. If $\fa_0$ is rational, by Proposition \ref{prop: ac-model}, we define 
$$\Psi_{\fa, \fa_0}([(\mX, \sum_{j=1}^k a_{j}\mD_j)\to T]):=[(\mX^\ac, \sum_{j=1}^k a_{0j}\mD^\ac_j)\to T] \in \mM^{\Kss,\red}_{\widehat{\mG}, \fa_0}(T).$$
Suppose $\fa_0$ is irrational. By Proposition \ref{prop: uniform polarization} and the construction of the chamber decomposition in Theorem \ref{thm: kps fcd}, there exists a rational polytope $Q$ containing $\fa$ with $\fa_0\in \partial{Q}$ such that $(\mX, \sum_{j=1}^k x_j\mD_j)\to T$ is an $\bR$-Gorenstein family of log Fano pairs for any $(x_1,...,x_k)\in Q^\circ$. In particular, the face of 
$Q$ containing $\fa_0$ as an interior point, denoted by $F$ (which is also a rational polytope), satisfies the following condition: $(\mX, \sum_{j=1}^k y_j\mD_j)\to T$ is an $\bR$-Gorenstein family of weak log Fano pairs for any $(y_1,...,y_k)$ contained in $F^\circ$. By Shokurov's polytope decomposition of  canonical models (e.g. \cite[Corollary 1.1.5]{BCHM10}), there exists a smaller rational polytope $F'\subset F^\circ$ containing $\fa_0$ such that $(\mX, \sum_{j=1}^k y_j\mD_j)\to T$ admits the same anti-canonical model for any $(y_1,...,y_k)\in F'^\circ(\bQ)$, denoted by $(\mX^\ac, \sum_{j=1}^k y_j\mD^\ac_j)\to T$. By the same argument as in the last paragraph of the proof of Proposition \ref{prop: uniform polarization} and the proof of Proposition \ref{prop: ac-model}, we see that $(\mX^\ac, \sum_{j=1}^k y_j\mD^\ac_j)\to T$ is a simultaneous anti-canonical model of $(\mX, \sum_{j=1}^k y_j\mD_j)\to T$ for any $(y_1,...,y_k)\in F'^\circ(\bQ)$. In particular, $(\mX^\ac, \sum_{j=1}^k \fa_{0j}\mD^\ac_j)\to T$ is a simultaneous anti-canonical model of $(\mX, \sum_{j=1}^k \fa_{0j}\mD_j)\to T$. Therefore, we can similarly define
$$\Psi_{\fa, \fa_0}([(\mX, \sum_{j=1}^k a_{j}\mD_j)\to T]):=[(\mX^\ac, \sum_{j=1}^k a_{0j}\mD^\ac_j)\to T] \in \mM^{\Kss,\red}_{\widehat{\mG}, \fa_0}(T).$$

We already establish $\Psi_{\fa,\fa_0}$, and now we turn to the last statement. Assume the $\bQ$-span of $\fa_0$ contains $A_i$. 
By Proposition \ref{prop: polytope LF} and the construction of the chamber decomposition in Theorem \ref{thm: kps fcd}, we see that for any $[(\mX, \sum_{j=1}^k a_j\mD_j)\to T]\in \mM^{\Kss,\red}_{\widehat{\mG}, \fa}(T)$,  the family $(\mX, \sum_{j=1}^k a_{0j}\mD_j)\to T$ is still a family of K-semistable log Fano pairs (not just weak log Fano pairs) after changing the coefficients. This means 
$[(\mX, \sum_{j=1}^k a_{0j}\mD_j)\to T]\in \mM^{\Kss,\red}_{\widehat{\mG}, \fa_0}(T)$ and we do not need the step of taking the simultaneous anti-canonical model as in the previous case. Thus $\Psi_{\fa, \fa_0}$ is itself an open embedding by \cite[Section 5]{LZ24} in this case. The proof is complete.
\end{proof}

\begin{remark}
    We note that wall crossing for log Fano pairs obtained by taking products of proportional log Fano pairs was studied in \cite{Pap24} based on the product theorem of K-stability \cite{zhuang19}.
\end{remark}

We end this section with the following example, which is explicitly explored in \cite{DJ+24}.

\begin{example}\label{example: p1p2}
Take $d=3, k=1$ and $P=[0, 1/2]$. We construct a log bounded subset $\widetilde{\mG}\subset \mG(d,k,P)$ as follows. Denote $X:=\bP^1\times \bP^2$ and $L:=\mO(2,2)$. Let $\pi: (X\times |L|, \mD)\to |L|$ be the universal family associated to $|L|$. We define $\widetilde{\mG}$ to be the set consisting of log smooth fibers of $\pi$. Observe the following: (1) $(X, cD)$ is K-semistable for any $(X, D)\in \widetilde{\mG}$ and any $c\in P^\circ(\bQ)$ (e.g. \cite[Theorem 2]{DJ+24}); (2) $\vol(-K_X-cD)\geq 12$ for any $(X, D)\in \widetilde{\mG}$ and any $c\in P$. By Proposition \ref{prop: G-complete}, the set $\widehat{\mG}$ is also log bounded. By Corollary \ref{cor: kps fcd}, there exists a finite polytope chamber decomposition of $P$ to control the variation of K-semi(poly)stability for couples in $\widehat{\mG}$. In fact, the implicit computation in \cite{DJ+24} shows that there is only one wall number in $P$, which is even irrational. 
\end{example}

\section{GIT-stability  vs  K-stability before the first wall}\label{sec: GIT=K}

As an application of the general theory established before, we explore the relationship between K-stability and GIT-stability in non-proportional setting. In this section, we fix $X$ to be a K-polystable Fano variety of dimension $d$ which is $\bQ$-factorial terminal, and
 fix $|L|$ to be the linear system of an effective line bundle $L$ on $X$. We also fix a rational polytope $P:=[0,b]$ such that $-K_X-bL$ is ample. Denote $G:=\Aut(X)$. Then $G$ is reductive by \cite{ABHLX20} and naturally acts on $|L|$. As in the proportional setting (e.g. \cite{Zhou21a}), it is natural to figure out the relationship between the GIT-stability of $D\in |L|$ (under $G$-action) and the K-stability of $(X, \epsilon D)$ for $0<\epsilon\ll 1$. We aim to address this problem in this section.

We say $(Y, B)$ is a degeneration of $(X, |L|)$ if there exists a family of couples $(\mY, \mB)\to  C\ni 0$ over a smooth pointed curve satisfying the following conditions:
\begin{enumerate}
\item for any closed point $t\in C$ with $t\ne 0$, we have $\mY_t\cong X$ and $\mB_t\in |L|$;
\item $(\mY_0, \mB_0)\cong (Y, B)$.
\end{enumerate}
For convenience, we just say $(Y, B)$ is produced by $(\mY, \mB)\to  C\ni 0$. Note that $Y$ may not be $\bQ$-Gorenstein, although $X$ is. We define a set of couples (denoted by $\kA$) as follows: we say a couple $(Y, B)$ belongs to $\kA$ if and only if  the following conditions are satisfied:
\begin{enumerate}
\item $(Y, B)$ is a degeneration of $(X, |L|)$;
\item $(Y, cB)$ is a K-semistable log Fano pair for some $c\in P$. \end{enumerate}

We also define an even larger set of couples (denoted by $\kA'$) as follows. We say a couple $(Y, B)$ is contained in $\kA'$ if and only if the following conditions are satisfied:
\begin{enumerate}
\item $Y$ is a normal projective variety of dimension $d$ and $B$ is an effective Weil divisor on $Y$;
\item $\vol(-K_Y-xB)\geq \vol(-K_X-bL)$ for any $x\in P$;
\item $(Y, cB)$ is a K-semistable weak log Fano pair for some $c\in P$.
\end{enumerate}

\begin{proposition}\label{prop: first gap}
Notation as above, we have the following conclusions:
\begin{enumerate}
\item both $\kA$ and $\kA'$ are log bounded; 
\item there exists a positive real number $0<\epsilon_0<b$ (which only depends on $X, |L|, P$) satisfying the following condition: for any $(Y, B)\in \kA$, if $(Y, cB)$ is log Fano for some $c\in (0,\epsilon_0)$, then $(Y, xB)$ is log Fano for any $x\in (0,\epsilon_0)$; for any $(Y, B)\in \kA'$, if $(Y, cB)$ is a weak log Fano pair for some $c\in (0,\epsilon_0)$, then $(Y, xB)$ is a weak log Fano pair for any $x\in (0,\epsilon_0)$;
\item moreover, the following condition is also satisfied: for any $(Y, B)\in \kA$, if $(Y, cB)$ is a K-semistable log Fano pair for some $c\in (0,\epsilon_0)$, then $(Y, xB)$ is a K-semistable log Fano pair  for any $x\in (0,\epsilon_0)$.
\end{enumerate}
\end{proposition}

\begin{proof}
The first statement is a consequence of Proposition \ref{prop: Abdd}.   The second statement is a consequence of Proposition \ref{prop: perfect fcd}. The third statement is a consequence of Corollary \ref{cor: fcd}.
\end{proof}

For any $a\in (0, b)$, we define a subset of $\kA$, denoted by $\kA_a$, satisfying the following condition: a couple $(Y, B)\in \kA$ is contained in $\kA_a$ if and only if there exists a family of couples $(\mY, \mB)\to C\ni 0$ over a smooth pointed curve (which produces $(Y, B)$) such that $(\mY, c\mB)\to C$ is a family of K-semistable log Fano pairs and $-K_\mY-c\mB$ is $\bQ$-Cartier for some $c\in (0, a)\cap \bQ$. Denote $\mP_{a}:=\{Y| \textit{$(Y, B)\in \kA_{a}$}\}$. We have the following result.

\begin{proposition}\label{prop: fano set}
Let $\epsilon_0$ be as in Proposition \ref{prop: first gap}, then
$\mP_{\epsilon_0}$ is a set of Fano varieties. Moreover, every variety in $\mP_{\epsilon_0}$ is isomorphic to $X$.
\end{proposition}

\begin{proof}
Arbitrarily choosing a couple $(Y, B)\in \kA_{\epsilon_0}$,  there exists a family of K-semistable log Fano pairs over a smooth pointed curve (which produces $(Y, B)$), denoted by $(\mY, c\mB)\to C\ni 0$ for some $c\in (0,\epsilon_0)\cap \bQ$, such that $K_\mY+c\mB$ is $\bQ$-Cartier. 

We aim to show that $Y$ is a Fano variety. By Proposition \ref{prop: first gap}, $(Y, xB)$ is a K-semistable log Fano pair for any $x\in (0, \epsilon_0)\cap \bQ$. Thus $Y$ is a K-semistable weak Fano variety since $-K_Y$ is big and nef. By Proposition \ref{prop: uniform polarization}, we see that $K_\mY+x\mB$ is $\bQ$-Cartier for any $x\in (0, \epsilon_0)\cap \bQ$, thus $\mY\to C$ is a family of weak Fano varieties with $K_\mY$ being $\bQ$-Cartier. Let $\mY'\to C$ be the anti-canonical model of $\mY\to C$ and applying the same argument as in the proof of Proposition \ref{prop: uniform polarization}, we see that $\mY'\to C$ actually carries fiberwise anti-canonical models of $\mY\to C$. Therefore, the following two morphisms
$$\mY\to C\quad \text{and}\quad \mY'\to C $$
are isomorphic over $C\setminus 0$ and birational over the distinguished point $0$.
Denote $\phi: \mY\to \mY'/C$. Then we have
$$K_{\mY}=\phi^* K_{\mY'}\quad \text{and}\quad K_{\mY_0}=\phi^*K_{\mY'_0}, $$
where we still denote $\phi: \mY_0\to \mY'_0$. Note that $\mY_0'$ is a K-semistable Fano variety, as it is the anti-canonical model of $Y\cong\mY_0$. Applying \cite[Theorem 3.4]{Zhou21a} we see that $\mY_0'\cong X$, and in particular $\mY_0'$ is a $\bQ$-factorial terminal Fano variety by assumption. This tells that $\mY_0\to \mY_0'$ is an isomorphism, which implies 
$$Y\cong \mY_0\cong \mY_0'\cong X.$$
The proof is complete.
\end{proof}

We are ready to prove the following theorem, which is analogue to \cite[Theorem 1.1]{Zhou21a} in non-proportional setting (see also \cite[Theorem 1.4]{ADL19}). For a special case of the following theorem where $X$ is a product of projective spaces, see \cite[Theorem 2.10]{DJ+24}.

\begin{theorem}\label{thm: GIT=K}
Let $\epsilon_0$ be as in Proposition \ref{prop: first gap} and suppose $\kA_{\epsilon_0}$ is not empty. Then K-(semi/poly)stability of $(X, c D)$ is identical to GIT-(semi/poly)stability of $D\in |L|$ under $G$-action for any $c\in (0,\epsilon_0)\cap \bQ$.
\end{theorem}

Let $\mD\subset X\times |L|$ be the universal divisor and $ (X\times |L|, \mD)\to |L|$ the universal family. Then the condition that $\kA_{\epsilon_0}$ is non-empty just means that $(X\times |L|, c\mD)\to |L|$ admits a fiber which is a K-semistable log Fano pair for some $c\in (0,\epsilon_0)\cap \bQ$ (thus for any $c\in (0,\epsilon_0)\cap \bQ$ by Proposition \ref{prop: first gap}).

\begin{proof}
Denote $\bP:=|L|$. For any $c\in (0,\epsilon_0)\cap \bQ$, we first compute the CM-line bundle for the universal family $\pi: (X\times \bP, c\mD)\to \bP$ with the polarization 
$$\mL_c:=-K_{X\times \bP/\bP}-c\mD.$$ 
Since the base is of Picard number one, it suffices to assume $\bP=\bP^1$ up to a base change (note that CM-line bundle is functorial by \cite[Proposition 3.8]{CP21}). By Definition \ref{def: CM}, we have the following computation:
$$\lambda_{\CM, (X\times \bP, c\mD, \mL_c; \pi)}=-\pi_*(-K_{X\times \bP/\bP}-c\mD)^{d+1}=(2-c)(d+1)\cdot\vol(-K_X-cL)\cdot \mO_\bP(1), $$
which is a positive ($\bQ$)-line bundle on the base.

For any $c\in (0,\epsilon_0)\cap \bQ$, let $(\mY, c\mB)\to C\ni 0$ be a family of K-semistable  log Fano pairs over a smooth pointed curve such that
\begin{enumerate}
\item $(\mY_t, \mB_t)$ is isomorphic to some fiber of $(X\times |L|, \mD)\to |L|$ for $t\ne 0$;
\item $K_\mY+c\mB$ is $\bQ$-Cartier.
\end{enumerate}
Denote $(Y, B):=(\mY_0, \mB_0)$. It is clear that $(Y, B)\in \kA_{\epsilon_0}$. By Proposition \ref{prop: fano set}, we see $Y\cong X$. By now, we have an isotrivial family $\mY\to C$ where every fiber is isomorphic to $X$. Let $ C\to |L|$ be the natural extending morphism induced by $C\setminus 0\to |L|$, then we have another isotrivial (actually trivial) family $\mY'\to C$ via the pull-back of $X\times |L|\to |L|$ along $C\to |L|$. By Lemma \ref{lem: etale torsor}, up to shrinking $C$ around $0\in C$ and taking an \'etale base change, we may assume both two families
$$\mY\to C\quad \text{and}\quad \mY'\to C $$
are trivial families and hence they are isomorphic over $C$. Denote $(\mY', \mB')\to C$ to be pull-back of the universal family $(X\times |L|, \mD)\to |L|$ via the base change along $C\to |L|$. Then we have two families of couples
$$(\mY, \mB)\to C\quad \text{and}\quad (\mY', \mB')\to C$$
such that they are isomorphic over $C\setminus 0$, i.e.
$$(\mY, \mB)\times_C (C\setminus 0) \cong (\mY', \mB')\times_C (C\setminus 0).$$
On the other hand, since $\mY$ and $\mY'$ are isomorphic over $C$, we  derive that $\mB-\mB'\sim 0$ (via identifying $\mY$ with $\mY'$) and hence $\mB_0\in |L|$. This tells that $(Y, B)$ is actually isomorphic to a fiber of the universal family $(X\times |L|, \mD)\to |L|$.

The rest of the proof is similar to the proof of \cite[Theorem 4.3]{Zhou21a}. We provide some details for the readers' convenience. Let $D$ be an element in $|L|$ such that $(X, cD)$ is K-(semi/poly) stable for some $c\in (0, \epsilon_0)\cap \bQ$, then $D$ is GIT-(semi/poly)stable by \cite[Theorem 2.22]{ADL19}. It remains to show the converse statement. Suppose $D\in |L|$ is GIT-(semi/poly)stable under $G$-action. Since $\kA_{\epsilon_0}$ is non-empty, by \cite{BLX22, Xu20}, there exists a family of couples 
$$(\mY'', \mB'')\to C''\ni 0''$$ 
over a smooth pointed curve satisfying the following conditions:
\begin{enumerate}
\item $(\mY'', \mB'')\to C''$ is the pull-back of $(X\times |L|, \mD)\to |L|$ along some morphism $C''\to |L|$;
\item $(\mY''_{0''}, \mB''_{0''})\cong (X, D)$;
\item $(\mY'', c\mB'')\times_ {C''} (C''\setminus 0'')\to C''\setminus 0'' $
is a family of K-semistable log Fano pairs.
\end{enumerate}
By the properness of K-moduli (e.g. \cite{LXZ22, Xu21}), up to a finite base change of $C''$, one could find another family $(\mY''', c\mB''')\to C''$ such that
\begin{enumerate}
\item $(\mY''', c\mB''')\times_{C''}(C''\setminus 0'')\cong (\mY'', c\mB'')\times_{C''}(C''\setminus 0'')$, and
\item $(\mY_{0''}''', c\mB_{0''}''')$ is a K-semistable log Fano pair.\end{enumerate}
By the argument in the previous paragraph, we know that $(\mY_{0''}''', \mB_{0''}''')$ is isomorphic to a fiber of the universal family $(X\times |L|, \mD)\to |L|$ and thus $B'''_{0''}\in |L|$ is GIT-semistable. 
By the separatedness of GIT moduli space, we know that $D$ and $B'''_{0''}$ lie on the same $S$-equivalence class under $G$-action. Thus $(X, cD)$ is K-(semi/poly)stable.
\end{proof}


\begin{lemma}\label{lem: etale torsor}
Let $W\to T$ be a flat projective morphism over a smooth base $T$ such that $-K_{W/T}$ is $\bQ$-Cartier and relatively ample. Suppose every geometric fiber is isomorphic to a given Fano variety $F$. Then the Isom scheme $\Isom(W, F\times T)\to T$ is an \'etale $\Aut(F)$-torsor.
\end{lemma}

\begin{proof}
Choose a sufficiently divisible positive integer $m$ such that $|-mK_F|$ is a very ample linear system. Then there is an embedding $F\to \bP^{N}$ induced by $|-mK_F|$ for some $N$. Thus $\Aut(F)$ is a linear algebraic subgroup of $\PGL(N+1)$ and we have the following \'etale $\Aut(F)$-torsor:
$$\PGL(N+1)\to \PGL(N+1)/\Aut(F). $$
Note that every fiber of $W\to T$ is isomorphic to $F$. We see that the Isom scheme $\Isom(W, F\times T)\to T$ is induced by the base change along some morphism 
$$T\to \PGL(N+1)/\Aut(F),$$ 
which implies that the Isom scheme is an \'etale $\Aut(F)$-torsor.
\end{proof}

\begin{remark}\label{rmk: GIT=K}
Suppose $\kA_{\epsilon_0}$ is non-empty. Then for any $c\in (0, \epsilon_0)\cap \bQ$ there exists a finite type Artin stack (denoted by $\mM^{\Kss}_c$) parametrizing K-semistable log Fano pairs $(Y, cB)$, where $(Y, B)\in \kA_{\epsilon_0}$; moreover, the stack descends to a separated projective scheme (denoted by $M^{\Kps}_c$) parametrizing K-polystable objects (e.g. \cite{Xu21, LXZ22}). Theorem \ref{thm: GIT=K} tells that the K-moduli coincides with GIT-moduli. More precisely, we have the following commutative diagram:
\begin{center}
	\begin{tikzcd}[column sep = 2em, row sep = 1.5em]
	 \mM^{\Kss}_c\arrow[rr,""] \arrow[d,"",swap]&& \mM_{X, |L|}^{\GIT}\arrow[d,"",swap]\\
	 M^{\Kps}_{c}\arrow[rr,""]&& M^{\GIT}_{X, |L|}
	 	\end{tikzcd}
\end{center}
where $\mM_{X, |L|}^{\GIT}$ (resp. $M^{\GIT}_{X, |L|}$) is the GIT stack (resp. GIT space) of $|L|$ under the $G$-action, and the two horizontal morphisms are isomorphic.

\end{remark}

\bibliography{reference.bib}
\end{document}